

\ifx\shlhetal\undefinedcontrolsequence\let\shlhetal\relax\fi

\def\fmtname{AmS-TeX}

\def\fmtversion{2.1}
\catcode`\@=11
\ifx\amstexloaded@\relax\catcode`\@=\active
  \endinput\else\let\amstexloaded@\relax\fi
\newlinechar=`\^^J
\def\W@{\immediate\write\sixt@@n}
\def\CR@{\W@{^^J\fmtname - Version \fmtversion^^J%
COPYRIGHT 1985, 1990, 1991 - AMERICAN MATHEMATICAL SOCIETY^^J%
Use of this macro package is not restricted provided^^J%
each use is acknowledged upon publication.^^J}}
\CR@ \everyjob{\CR@}
\message{Loading definitions for}
\message{misc utility macros,}
\toksdef\toks@@=2
\long\def\rightappend@#1\to#2{\toks@{\\{#1}}\toks@@
 =\expandafter{#2}\xdef#2{\the\toks@@\the\toks@}\toks@{}\toks@@{}}
\def\alloclist@{}
\newif\ifalloc@
\def\showallocations{{\def\\{\immediate\write\m@ne}\alloclist@}\alloc@true}
\def\alloc@#1#2#3#4#5{\global\advance\count1#1by\@ne
 \ch@ck#1#4#2\allocationnumber=\count1#1
 \global#3#5=\allocationnumber
 \edef\next@{\string#5=\string#2\the\allocationnumber}%
 \expandafter\rightappend@\next@\to\alloclist@}
\newcount\count@@
\newcount\count@@@
\def\FN@{\futurelet\next}
\def\DN@{\def\next@}
\def\DNii@{\def\nextii@}
\def\RIfM@{\relax\ifmmode}
\def\RIfMIfI@{\relax\ifmmode\ifinner}
\def\setboxz@h{\setbox\z@\hbox}
\def\wdz@{\wd\z@}
\def\boxz@{\box\z@}
\def\setbox@ne{\setbox\@ne}
\def\wd@ne{\wd\@ne}
\def\iterate{\body\expandafter\iterate\else\fi}
\def\err@#1{\errmessage{AmS-TeX error: #1}}
\newhelp\defaulthelp@{Sorry, I already gave what help I could...^^J
Maybe you should try asking a human?^^J
An error might have occurred before I noticed any problems.^^J
``If all else fails, read the instructions.''}
\def\Err@{\errhelp\defaulthelp@\err@}
\def\eat@#1{}
\def\in@#1#2{\def\in@@##1#1##2##3\in@@{\ifx\in@##2\in@false\else\in@true\fi}%
 \in@@#2#1\in@\in@@}
\newif\ifin@
\def\space@.{\futurelet\space@\relax}
\space@. %
\newhelp\athelp@
{Only certain combinations beginning with @ make sense to me.^^J
Perhaps you wanted \string\@\space for a printed @?^^J
I've ignored the character or group after @.}
{\catcode`\~=\active 
 \lccode`\~=`\@ \lowercase{\gdef~{\FN@\at@}}}
\def\at@{\let\next@\at@@
 \ifcat\noexpand\next a\else\ifcat\noexpand\next0\else
 \ifcat\noexpand\next\relax\else
   \let\next\at@@@\fi\fi\fi
 \next@}
\def\at@@#1{\expandafter
 \ifx\csname\space @\string#1\endcsname\relax
  \expandafter\at@@@ \else
  \csname\space @\string#1\expandafter\endcsname\fi}
\def\at@@@#1{\errhelp\athelp@ \err@{\Invalid@@ @}}
\def\atdef@#1{\expandafter\def\csname\space @\string#1\endcsname}
\newhelp\defahelp@{If you typed \string\define\space cs instead of
\string\define\string\cs\space^^J
I've substituted an inaccessible control sequence so that your^^J
definition will be completed without mixing me up too badly.^^J
If you typed \string\define{\string\cs} the inaccessible control sequence^^J
was defined to be \string\cs, and the rest of your^^J
definition appears as input.}
\newhelp\defbhelp@{I've ignored your definition, because it might^^J
conflict with other uses that are important to me.}
\def\define{\FN@\define@}
\def\define@{\ifcat\noexpand\next\relax
 \expandafter\define@@\else\errhelp\defahelp@                               
 \err@{\string\define\space must be followed by a control
 sequence}\expandafter\def\expandafter\nextii@\fi}                          
\def\undefined@@@@@@@@@@{}
\def\preloaded@@@@@@@@@@{}
\def\next@@@@@@@@@@{}
\def\define@@#1{\ifx#1\relax\errhelp\defbhelp@                              
 \err@{\string#1\space is already defined}\DN@{\DNii@}\else
 \expandafter\ifx\csname\expandafter\eat@\string                            
 #1@@@@@@@@@@\endcsname\undefined@@@@@@@@@@\errhelp\defbhelp@
 \err@{\string#1\space can't be defined}\DN@{\DNii@}\else
 \expandafter\ifx\csname\expandafter\eat@\string#1\endcsname\relax          
 \global\let#1\undefined\DN@{\def#1}\else\errhelp\defbhelp@
 \err@{\string#1\space is already defined}\DN@{\DNii@}\fi
 \fi\fi\next@}

\def\predefine#1#2{\let#1#2}
\def\undefine#1{\let#1\undefined}
\message{page layout,}
\newdimen\captionwidth@
\captionwidth@\hsize
\advance\captionwidth@-1.5in
\def\pagewidth#1{\hsize#1\relax
 \captionwidth@\hsize\advance\captionwidth@-1.5in}
\def\pageheight#1{\vsize#1\relax}
\def\hcorrection#1{\advance\hoffset#1\relax}
\def\vcorrection#1{\advance\voffset#1\relax}
\message{accents/punctuation,}

\let\graveaccent\`
\let\acuteaccent\'
\let\tildeaccent\~
\let\hataccent\^
\let\underscore\_
\let\B\=
\let\D\.
\let\ic@\/
\def\/{\unskip\ic@}
\def\textfonti{\the\textfont\@ne}
\def\t#1#2{{\edef\next@{\the\font}\textfonti\accent"7F \next@#1#2}}
\def~{\unskip\nobreak\ \ignorespaces}
\def\.{.\spacefactor\@m}
\atdef@;{\leavevmode\null;}
\atdef@:{\leavevmode\null:}
\atdef@?{\leavevmode\null?}
\edef\@{\string @}
\def\relaxnext@{\let\next\relax}
\atdef@-{\relaxnext@\leavevmode
 \DN@{\ifx\next-\DN@-{\FN@\nextii@}\else
  \DN@{\leavevmode\hbox{-}}\fi\next@}%
 \DNii@{\ifx\next-\DN@-{\leavevmode\hbox{---}}\else
  \DN@{\leavevmode\hbox{--}}\fi\next@}%
 \FN@\next@}
\def\srdr@{\kern.16667em}
\def\drsr@{\kern.02778em}
\def\sldl@{\drsr@}
\def\dlsl@{\srdr@}
\atdef@"{\unskip\relaxnext@
 \DN@{\ifx\next\space@\DN@. {\FN@\nextii@}\else
  \DN@.{\FN@\nextii@}\fi\next@.}%
 \DNii@{\ifx\next`\DN@`{\FN@\nextiii@}\else
  \ifx\next\lq\DN@\lq{\FN@\nextiii@}\else
  \DN@####1{\FN@\nextiv@}\fi\fi\next@}%
 \def\nextiii@{\ifx\next`\DN@`{\sldl@``}\else\ifx\next\lq
  \DN@\lq{\sldl@``}\else\DN@{\dlsl@`}\fi\fi\next@}%
 \def\nextiv@{\ifx\next'\DN@'{\srdr@''}\else
  \ifx\next\rq\DN@\rq{\srdr@''}\else\DN@{\drsr@'}\fi\fi\next@}%
 \FN@\next@}

\def\textfontii{\the\textfont\tw@}
\def\lbrace@{\delimiter"4266308 }
\def\rbrace@{\delimiter"5267309 }
\def\{{\RIfM@\lbrace@\else{\textfontii f}\spacefactor\@m\fi}
\def\}{\RIfM@\rbrace@\else
 \let\@sf\empty\ifhmode\edef\@sf{\spacefactor\the\spacefactor}\fi
 {\textfontii g}\@sf\relax\fi}
\let\lbrace\{
\let\rbrace\}
\def\AmSTeX{{\textfontii A\kern-.1667em%
  \lower.5ex\hbox{M}\kern-.125emS}-\TeX}
\message{line and page breaks,}
\def\vmodeerr@#1{\Err@{\string#1\space not allowed between paragraphs}}
\def\mathmodeerr@#1{\Err@{\string#1\space not allowed in math mode}}
\def\linebreak{\RIfM@\mathmodeerr@\linebreak\else
 \ifhmode\unskip\unkern\break\else\vmodeerr@\linebreak\fi\fi}

\newskip\saveskip@
\def\allowlinebreak{\RIfM@\mathmodeerr@\allowlinebreak\else
 \ifhmode\saveskip@\lastskip\unskip
 \allowbreak\ifdim\saveskip@>\z@\hskip\saveskip@\fi
 \else\vmodeerr@\allowlinebreak\fi\fi}
\def\nolinebreak{\RIfM@\mathmodeerr@\nolinebreak\else
 \ifhmode\saveskip@\lastskip\unskip
 \nobreak\ifdim\saveskip@>\z@\hskip\saveskip@\fi
 \else\vmodeerr@\nolinebreak\fi\fi}
\def\newline{\relaxnext@
 \DN@{\RIfM@\expandafter\mathmodeerr@\expandafter\newline\else
  \ifhmode\ifx\next\par\else
  \expandafter\unskip\expandafter\null\expandafter\hfill\expandafter\break\fi
  \else
  \expandafter\vmodeerr@\expandafter\newline\fi\fi}%
 \FN@\next@}
\def\dmatherr@#1{\Err@{\string#1\space not allowed in display math mode}}
\def\nondmatherr@#1{\Err@{\string#1\space not allowed in non-display math
 mode}}
\def\onlydmatherr@#1{\Err@{\string#1\space allowed only in display math mode}}
\def\nonmatherr@#1{\Err@{\string#1\space allowed only in math mode}}
\def\mathbreak{\RIfMIfI@\break\else
 \dmatherr@\mathbreak\fi\else\nonmatherr@\mathbreak\fi}
\def\nomathbreak{\RIfMIfI@\nobreak\else
 \dmatherr@\nomathbreak\fi\else\nonmatherr@\nomathbreak\fi}
\def\allowmathbreak{\RIfMIfI@\allowbreak\else
 \dmatherr@\allowmathbreak\fi\else\nonmatherr@\allowmathbreak\fi}
\def\pagebreak{\RIfM@
 \ifinner\nondmatherr@\pagebreak\else\postdisplaypenalty-\@M\fi
 \else\ifvmode\removelastskip\break\else\vadjust{\break}\fi\fi}
\def\nopagebreak{\RIfM@
 \ifinner\nondmatherr@\nopagebreak\else\postdisplaypenalty\@M\fi
 \else\ifvmode\nobreak\else\vadjust{\nobreak}\fi\fi}
\def\nonvmodeerr@#1{\Err@{\string#1\space not allowed within a paragraph
 or in math}}
\def\vnonvmode@#1#2{\relaxnext@\DNii@{\ifx\next\par\DN@{#1}\else
 \DN@{#2}\fi\next@}%
 \ifvmode\DN@{#1}\else
 \DN@{\FN@\nextii@}\fi\next@}
\def\newpage{\vnonvmode@{\vfill\break}{\nonvmodeerr@\newpage}}
\def\smallpagebreak{\vnonvmode@\smallbreak{\nonvmodeerr@\smallpagebreak}}
\def\medpagebreak{\vnonvmode@\medbreak{\nonvmodeerr@\medpagebreak}}
\def\bigpagebreak{\vnonvmode@\bigbreak{\nonvmodeerr@\bigpagebreak}}
\def\NoBlackBoxes{\global\overfullrule\z@}
\def\BlackBoxes{\global\overfullrule5\p@}
\def\Invalid@#1{\def#1{\Err@{\Invalid@@\string#1}}}
\def\Invalid@@{Invalid use of }
\message{figures,}
\Invalid@\caption
\Invalid@\captionwidth
\newdimen\smallcaptionwidth@
\def\topspace{\mid@false\ins@}
\def\midspace{\mid@true\ins@}
\newif\ifmid@
\def\captionfont@{}
\def\ins@#1{\relaxnext@\allowbreak
 \smallcaptionwidth@\captionwidth@\gdef\thespace@{#1}%
 \DN@{\ifx\next\space@\DN@. {\FN@\nextii@}\else
  \DN@.{\FN@\nextii@}\fi\next@.}%
 \DNii@{\ifx\next\caption\DN@\caption{\FN@\nextiii@}%
  \else\let\next@\nextiv@\fi\next@}%
 \def\nextiv@{\vnonvmode@
  {\ifmid@\expandafter\midinsert\else\expandafter\topinsert\fi
   \vbox to\thespace@{}\endinsert}
  {\ifmid@\nonvmodeerr@\midspace\else\nonvmodeerr@\topspace\fi}}%
 \def\nextiii@{\ifx\next\captionwidth\expandafter\nextv@
  \else\expandafter\nextvi@\fi}%
 \def\nextv@\captionwidth##1##2{\smallcaptionwidth@##1\relax\nextvi@{##2}}%
 \def\nextvi@##1{\def\thecaption@{\captionfont@##1}%
  \DN@{\ifx\next\space@\DN@. {\FN@\nextvii@}\else
   \DN@.{\FN@\nextvii@}\fi\next@.}%
  \FN@\next@}%
 \def\nextvii@{\vnonvmode@
  {\ifmid@\expandafter\midinsert\else
  \expandafter\topinsert\fi\vbox to\thespace@{}\nobreak\smallskip
  \setboxz@h{\noindent\ignorespaces\thecaption@\unskip}%
  \ifdim\wdz@>\smallcaptionwidth@\centerline{\vbox{\hsize\smallcaptionwidth@
   \noindent\ignorespaces\thecaption@\unskip}}%
  \else\centerline{\boxz@}\fi\endinsert}
  {\ifmid@\nonvmodeerr@\midspace
  \else\nonvmodeerr@\topspace\fi}}%
 \FN@\next@}
\message{comments,}
\def\newcodes@{\catcode`\\12\catcode`\{12\catcode`\}12\catcode`\#12%
 \catcode`\%12\relax}
\def\oldcodes@{\catcode`\\0\catcode`\{1\catcode`\}2\catcode`\#6%
 \catcode`\%14\relax}
\def\comment{\newcodes@\endlinechar=10 \comment@}
{\lccode`\0=`\\
\lowercase{\gdef\comment@#1^^J{\comment@@#10endcomment\comment@@@}%
\gdef\comment@@#10endcomment{\FN@\comment@@@}%
\gdef\comment@@@#1\comment@@@{\ifx\next\comment@@@\let\next\comment@
 \else\def\next{\oldcodes@\endlinechar=`\^^M\relax}%
 \fi\next}}}
\def\pr@m@s{\ifx'\next\DN@##1{\prim@s}\else\let\next@\egroup\fi\next@}
\def\prime{{\null\prime@\null}}
\mathchardef\prime@="0230
\let\dsize\displaystyle

\let\ssize\scriptstyle

\message{math spacing,}
\def\,{\RIfM@\mskip\thinmuskip\relax\else\kern.16667em\fi}
\def\!{\RIfM@\mskip-\thinmuskip\relax\else\kern-.16667em\fi}
\let\thinspace\,
\let\negthinspace\!
\def\medspace{\RIfM@\mskip\medmuskip\relax\else\kern.222222em\fi}
\def\negmedspace{\RIfM@\mskip-\medmuskip\relax\else\kern-.222222em\fi}
\def\thickspace{\RIfM@\mskip\thickmuskip\relax\else\kern.27777em\fi}
\let\;\thickspace
\def\negthickspace{\RIfM@\mskip-\thickmuskip\relax\else
 \kern-.27777em\fi}
\atdef@,{\RIfM@\mskip.1\thinmuskip\else\leavevmode\null,\fi}
\atdef@!{\RIfM@\mskip-.1\thinmuskip\else\leavevmode\null!\fi}
\atdef@.{\RIfM@&&\else\leavevmode.\spacefactor3000 \fi}
\def\and{\DOTSB\;\mathchar"3026 \;}

\message{fractions,}
\def\frac#1#2{{#1\over#2}}

\newdimen\ex@
\ex@.2326ex
\Invalid@\thickness
\def\thickfrac{\relaxnext@
 \DN@{\ifx\next\thickness\let\next@\nextii@\else
 \DN@{\nextii@\thickness1}\fi\next@}%
 \DNii@\thickness##1##2##3{{##2\above##1\ex@##3}}%
 \FN@\next@}

\def\thickfracwithdelims#1#2{\relaxnext@\def\ldelim@{#1}\def\rdelim@{#2}%
 \DN@{\ifx\next\thickness\let\next@\nextii@\else
 \DN@{\nextii@\thickness1}\fi\next@}%
 \DNii@\thickness##1##2##3{{##2\abovewithdelims
 \ldelim@\rdelim@##1\ex@##3}}%
 \FN@\next@}

\def\:{\nobreak\hskip.1111em\mathpunct{}\nonscript\mkern-\thinmuskip{:}\hskip
 .3333emplus.0555em\relax}
\def\snug{\unskip\kern-\mathsurround}
\message{smash commands,}
\def\topsmash{\top@true\bot@false\smash@}
\def\botsmash{\top@false\bot@true\smash@}
\newif\iftop@
\newif\ifbot@
\def\smash{\top@true\bot@true\smash@}
\def\smash@{\RIfM@\expandafter\mathpalette\expandafter\mathsm@sh\else
 \expandafter\makesm@sh\fi}
\def\finsm@sh{\iftop@\ht\z@\z@\fi\ifbot@\dp\z@\z@\fi\leavevmode\boxz@}
\message{large operator symbols,}
\def\LimitsOnSums{\global\let\slimits@\displaylimits}
\def\NoLimitsOnSums{\global\let\slimits@\nolimits}
\LimitsOnSums
\mathchardef\coprod@="1360       \def\coprod{\DOTSB\coprod@\slimits@}
\mathchardef\bigvee@="1357       \def\bigvee{\DOTSB\bigvee@\slimits@}
\mathchardef\bigwedge@="1356     \def\bigwedge{\DOTSB\bigwedge@\slimits@}
\mathchardef\biguplus@="1355     \def\biguplus{\DOTSB\biguplus@\slimits@}
\mathchardef\bigcap@="1354       \def\bigcap{\DOTSB\bigcap@\slimits@}
\mathchardef\bigcup@="1353       \def\bigcup{\DOTSB\bigcup@\slimits@}
\mathchardef\prod@="1351         \def\prod{\DOTSB\prod@\slimits@}
\mathchardef\sum@="1350          \def\sum{\DOTSB\sum@\slimits@}
\mathchardef\bigotimes@="134E    \def\bigotimes{\DOTSB\bigotimes@\slimits@}
\mathchardef\bigoplus@="134C     \def\bigoplus{\DOTSB\bigoplus@\slimits@}
\mathchardef\bigodot@="134A      \def\bigodot{\DOTSB\bigodot@\slimits@}
\mathchardef\bigsqcup@="1346     \def\bigsqcup{\DOTSB\bigsqcup@\slimits@}
\message{integrals,}
\def\LimitsOnInts{\global\let\ilimits@\displaylimits}
\def\NoLimitsOnInts{\global\let\ilimits@\nolimits}
\NoLimitsOnInts
\def\int{\DOTSI\intop\ilimits@}
\def\oint{\DOTSI\ointop\ilimits@}
\def\intic@{\mathchoice{\hskip.5em}{\hskip.4em}{\hskip.4em}{\hskip.4em}}
\def\negintic@{\mathchoice
 {\hskip-.5em}{\hskip-.4em}{\hskip-.4em}{\hskip-.4em}}
\def\intkern@{\mathchoice{\!\!\!}{\!\!}{\!\!}{\!\!}}
\def\intdots@{\mathchoice{\plaincdots@}
 {{\cdotp}\mkern1.5mu{\cdotp}\mkern1.5mu{\cdotp}}
 {{\cdotp}\mkern1mu{\cdotp}\mkern1mu{\cdotp}}
 {{\cdotp}\mkern1mu{\cdotp}\mkern1mu{\cdotp}}}
\newcount\intno@
\def\iint{\DOTSI\intno@\tw@\FN@\ints@}
\def\iiint{\DOTSI\intno@\thr@@\FN@\ints@}
\def\iiiint{\DOTSI\intno@4 \FN@\ints@}
\def\idotsint{\DOTSI\intno@\z@\FN@\ints@}
\def\ints@{\findlimits@\ints@@}
\newif\iflimtoken@
\newif\iflimits@
\def\findlimits@{\limtoken@true\ifx\next\limits\limits@true
 \else\ifx\next\nolimits\limits@false\else
 \limtoken@false\ifx\ilimits@\nolimits\limits@false\else
 \ifinner\limits@false\else\limits@true\fi\fi\fi\fi}
\def\multint@{\int\ifnum\intno@=\z@\intdots@                                
 \else\intkern@\fi                                                          
 \ifnum\intno@>\tw@\int\intkern@\fi                                         
 \ifnum\intno@>\thr@@\int\intkern@\fi                                       
 \int}                                                                      
\def\multintlimits@{\intop\ifnum\intno@=\z@\intdots@\else\intkern@\fi
 \ifnum\intno@>\tw@\intop\intkern@\fi
 \ifnum\intno@>\thr@@\intop\intkern@\fi\intop}
\def\ints@@{\iflimtoken@                                                    
 \def\ints@@@{\iflimits@\negintic@\mathop{\intic@\multintlimits@}\limits    
  \else\multint@\nolimits\fi                                                
  \eat@}                                                                    
 \else                                                                      
 \def\ints@@@{\iflimits@\negintic@
  \mathop{\intic@\multintlimits@}\limits\else
  \multint@\nolimits\fi}\fi\ints@@@}
\def\LimitsOnNames{\global\let\nlimits@\displaylimits}
\def\NoLimitsOnNames{\global\let\nlimits@\nolimits@}
\LimitsOnNames
\def\nolimits@{\relaxnext@
 \DN@{\ifx\next\limits\DN@\limits{\nolimits}\else
  \let\next@\nolimits\fi\next@}%
 \FN@\next@}
\message{operator names,}
\def\newmcodes@{\mathcode`\'"27\mathcode`\*"2A\mathcode`\."613A%
 \mathcode`\-"2D\mathcode`\/"2F\mathcode`\:"603A }
\def\operatorname#1{\mathop{\newmcodes@\kern\z@\fam\z@#1}\nolimits@}
\def\operatornamewithlimits#1{\mathop{\newmcodes@\kern\z@\fam\z@#1}\nlimits@}
\def\qopname@#1{\mathop{\fam\z@#1}\nolimits@}
\def\qopnamewl@#1{\mathop{\fam\z@#1}\nlimits@}
\def\arccos{\qopname@{arccos}}
\def\arcsin{\qopname@{arcsin}}
\def\arctan{\qopname@{arctan}}
\def\arg{\qopname@{arg}}
\def\cos{\qopname@{cos}}
\def\cosh{\qopname@{cosh}}
\def\cot{\qopname@{cot}}
\def\coth{\qopname@{coth}}
\def\csc{\qopname@{csc}}
\def\deg{\qopname@{deg}}
\def\det{\qopnamewl@{det}}
\def\dim{\qopname@{dim}}
\def\exp{\qopname@{exp}}
\def\gcd{\qopnamewl@{gcd}}
\def\hom{\qopname@{hom}}
\def\inf{\qopnamewl@{inf}}
\def\injlim{\qopnamewl@{inj\,lim}}
\def\ker{\qopname@{ker}}
\def\lg{\qopname@{lg}}
\def\lim{\qopnamewl@{lim}}
\def\liminf{\qopnamewl@{lim\,inf}}
\def\limsup{\qopnamewl@{lim\,sup}}
\def\ln{\qopname@{ln}}
\def\log{\qopname@{log}}
\def\max{\qopnamewl@{max}}
\def\min{\qopnamewl@{min}}
\def\Pr{\qopnamewl@{Pr}}
\def\projlim{\qopnamewl@{proj\,lim}}
\def\sec{\qopname@{sec}}
\def\sin{\qopname@{sin}}
\def\sinh{\qopname@{sinh}}
\def\sup{\qopnamewl@{sup}}
\def\tan{\qopname@{tan}}
\def\tanh{\qopname@{tanh}}
\def\varinjlim{\mathop{\vtop{\ialign{##\crcr
 \hfil\rm lim\hfil\crcr\noalign{\nointerlineskip}\rightarrowfill\crcr
 \noalign{\nointerlineskip\kern-\ex@}\crcr}}}}
\def\varprojlim{\mathop{\vtop{\ialign{##\crcr
 \hfil\rm lim\hfil\crcr\noalign{\nointerlineskip}\leftarrowfill\crcr
 \noalign{\nointerlineskip\kern-\ex@}\crcr}}}}
\def\varliminf{\mathop{\underline{\vrule height\z@ depth.2exwidth\z@
 \hbox{\rm lim}}}}

\newdimen\buffer@
\buffer@\fontdimen13 \tenex
\newdimen\buffer
\buffer\buffer@

\def\ResetBuffer{\fontdimen13 \tenex\buffer@\global\buffer\buffer@}
\def\shave#1{\mathop{\hbox{$\m@th\fontdimen13 \tenex\z@                     
 \displaystyle{#1}$}}\fontdimen13 \tenex\buffer}

\message{multilevel sub/superscripts,}
\Invalid@\\
\def\Let@{\relax\iffalse{\fi\let\\=\cr\iffalse}\fi}
\Invalid@\vspace
\def\vspace@{\def\vspace##1{\crcr\noalign{\vskip##1\relax}}}
\def\multilimits@{\bgroup\vspace@\Let@
 \baselineskip\fontdimen10 \scriptfont\tw@
 \advance\baselineskip\fontdimen12 \scriptfont\tw@
 \lineskip\thr@@\fontdimen8 \scriptfont\thr@@
 \lineskiplimit\lineskip
 \vbox\bgroup\ialign\bgroup\hfil$\m@th\scriptstyle{##}$\hfil\crcr}
\def\Sb{_\multilimits@}
\def\endSb{\crcr\egroup\egroup\egroup}
\def\Sp{^\multilimits@}

\def\spreadlines#1{\RIfMIfI@\onlydmatherr@\spreadlines\else
 \openup#1\relax\fi\else\onlydmatherr@\spreadlines\fi}
\def\Mathstrut@{\copy\Mathstrutbox@}
\newbox\Mathstrutbox@
\setbox\Mathstrutbox@\null
\setboxz@h{$\m@th($}
\ht\Mathstrutbox@\ht\z@
\dp\Mathstrutbox@\dp\z@
\message{matrices,}
\newdimen\spreadmlines@
\def\spreadmatrixlines#1{\RIfMIfI@
 \onlydmatherr@\spreadmatrixlines\else
 \spreadmlines@#1\relax\fi\else\onlydmatherr@\spreadmatrixlines\fi}
\def\matrix{\null\,\vcenter\bgroup\Let@\vspace@
 \normalbaselines\openup\spreadmlines@\ialign
 \bgroup\hfil$\m@th##$\hfil&&\quad\hfil$\m@th##$\hfil\crcr
 \Mathstrut@\crcr\noalign{\kern-\baselineskip}}
\def\endmatrix{\crcr\Mathstrut@\crcr\noalign{\kern-\baselineskip}\egroup
 \egroup\,}
\def\format{\crcr\egroup\iffalse{\fi\ifnum`}=0 \fi\format@}
\newtoks\hashtoks@
\hashtoks@{#}
\def\format@#1\\{\def\preamble@{#1}%
 \def\l{$\m@th\the\hashtoks@$\hfil}%
 \def\c{\hfil$\m@th\the\hashtoks@$\hfil}%
 \def\r{\hfil$\m@th\the\hashtoks@$}%
 \edef\preamble@@{\preamble@}\ifnum`{=0 \fi\iffalse}\fi
 \ialign\bgroup\span\preamble@@\crcr}
\def\smallmatrix{\null\,\vcenter\bgroup\vspace@\Let@
 \baselineskip9\ex@\lineskip\ex@
 \ialign\bgroup\hfil$\m@th\scriptstyle{##}$\hfil&&\thickspace\hfil
 $\m@th\scriptstyle{##}$\hfil\crcr}
\def\endsmallmatrix{\crcr\egroup\egroup\,}

\newmuskip\dotsspace@
\dotsspace@1.5mu
\def\strip@#1 {#1}
\def\spacehdots#1\for#2{\multispan{#2}\xleaders
 \hbox{$\m@th\mkern\strip@#1 \dotsspace@.\mkern\strip@#1 \dotsspace@$}\hfill}
\def\hdotsfor#1{\spacehdots\@ne\for{#1}}
\def\multispan@#1{\omit\mscount#1\unskip\loop\ifnum\mscount>\@ne\sp@n\repeat}
\def\spaceinnerhdots#1\for#2\after#3{\multispan@{\strip@#2 }#3\xleaders
 \hbox{$\m@th\mkern\strip@#1 \dotsspace@.\mkern\strip@#1 \dotsspace@$}\hfill}
\def\innerhdotsfor#1\after#2{\spaceinnerhdots\@ne\for#1\after{#2}}
\def\cases{\bgroup\spreadmlines@\jot\left\{\,\matrix\format\l&\quad\l\\}
\def\endcases{\endmatrix\right.\egroup}
\message{multiline displays,}
\newif\ifinany@
\newif\ifinalign@
\newif\ifingather@
\def\strut@{\copy\strutbox@}
\newbox\strutbox@
\setbox\strutbox@\hbox{\vrule height8\p@ depth3\p@ width\z@}
\def\topaligned{\null\,\vtop\aligned@}
\def\botaligned{\null\,\vbox\aligned@}
\def\aligned{\null\,\vcenter\aligned@}
\def\aligned@{\bgroup\vspace@\Let@
 \ifinany@\else\openup\jot\fi\ialign
 \bgroup\hfil\strut@$\m@th\displaystyle{##}$&
 $\m@th\displaystyle{{}##}$\hfil\crcr}
\def\endaligned{\crcr\egroup\egroup}

\def\alignedat#1{\null\,\vcenter\bgroup\doat@{#1}\vspace@\Let@
 \ifinany@\else\openup\jot\fi\ialign\bgroup\span\preamble@@\crcr}
\newcount\atcount@
\def\doat@#1{\toks@{\hfil\strut@$\m@th
 \displaystyle{\the\hashtoks@}$&$\m@th\displaystyle
 {{}\the\hashtoks@}$\hfil}
 \atcount@#1\relax\advance\atcount@\m@ne                                    
 \loop\ifnum\atcount@>\z@\toks@=\expandafter{\the\toks@&\hfil$\m@th
 \displaystyle{\the\hashtoks@}$&$\m@th
 \displaystyle{{}\the\hashtoks@}$\hfil}\advance
  \atcount@\m@ne\repeat                                                     
 \xdef\preamble@{\the\toks@}\xdef\preamble@@{\preamble@}}

\def\gathered{\null\,\vcenter\bgroup\vspace@\Let@
 \ifinany@\else\openup\jot\fi\ialign
 \bgroup\hfil\strut@$\m@th\displaystyle{##}$\hfil\crcr}
\def\endgathered{\crcr\egroup\egroup}
\newif\iftagsleft@
\def\TagsOnLeft{\global\tagsleft@true}
\def\TagsOnRight{\global\tagsleft@false}
\TagsOnLeft
\newif\ifmathtags@
\def\TagsAsMath{\global\mathtags@true}
\def\TagsAsText{\global\mathtags@false}
\TagsAsText
\def\tagform@#1{\hbox{\rm(\ignorespaces#1\unskip)}}
\def\thetag{\leavevmode\tagform@}
\def\tag#1$${\iftagsleft@\leqno\else\eqno\fi                                
 \maketag@#1\maketag@                                                       
 $$}                                                                        
\def\maketag@{\FN@\maketag@@}
\def\maketag@@{\ifx\next"\expandafter\maketag@@@\else\expandafter\maketag@@@@
 \fi}
\def\maketag@@@"#1"#2\maketag@{\hbox{\rm#1}}                                
\def\maketag@@@@#1\maketag@{\ifmathtags@\tagform@{$\m@th#1$}\else
 \tagform@{#1}\fi}
\interdisplaylinepenalty\@M
\def\allowdisplaybreaks{\RIfMIfI@
 \onlydmatherr@\allowdisplaybreaks\else
 \interdisplaylinepenalty\z@\fi\else\onlydmatherr@\allowdisplaybreaks\fi}
\Invalid@\allowdisplaybreak
\Invalid@\displaybreak
\Invalid@\intertext
\def\allowdisplaybreak@{\def\allowdisplaybreak{\crcr\noalign{\allowbreak}}}
\def\displaybreak@{\def\displaybreak{\crcr\noalign{\break}}}
\def\intertext@{\def\intertext##1{\crcr\noalign{%
 \penalty\postdisplaypenalty \vskip\belowdisplayskip
 \vbox{\normalbaselines\noindent##1}%
 \penalty\predisplaypenalty \vskip\abovedisplayskip}}}
\newskip\centering@
\centering@\z@ plus\@m\p@
\def\align{\relax\ifingather@\DN@{\csname align (in
  \string\gather)\endcsname}\else
 \ifmmode\ifinner\DN@{\onlydmatherr@\align}\else
  \let\next@\align@\fi
 \else\DN@{\onlydmatherr@\align}\fi\fi\next@}
\newhelp\andhelp@
{An extra & here is so disastrous that you should probably exit^^J
and fix things up.}
\newif\iftag@
\newcount\and@
\def\align@{\inalign@true\inany@true
 \vspace@\allowdisplaybreak@\displaybreak@\intertext@
 \def\tag{\global\tag@true\ifnum\and@=\z@\DN@{&&}\else
          \DN@{&}\fi\next@}%
 \iftagsleft@\DN@{\csname align \endcsname}\else
  \DN@{\csname align \space\endcsname}\fi\next@}
\def\Tag@{\iftag@\else\errhelp\andhelp@\err@{Extra & on this line}\fi}
\newdimen\lwidth@
\newdimen\rwidth@
\newdimen\maxlwidth@
\newdimen\maxrwidth@
\newdimen\totwidth@
\def\measure@#1\endalign{\lwidth@\z@\rwidth@\z@\maxlwidth@\z@\maxrwidth@\z@
 \global\and@\z@                                                            
 \setbox@ne\vbox                                                            
  {\everycr{\noalign{\global\tag@false\global\and@\z@}}\Let@                
  \halign{\setboxz@h{$\m@th\displaystyle{\@lign##}$}
   \global\lwidth@\wdz@                                                     
   \ifdim\lwidth@>\maxlwidth@\global\maxlwidth@\lwidth@\fi                  
   \global\advance\and@\@ne                                                 
   &\setboxz@h{$\m@th\displaystyle{{}\@lign##}$}\global\rwidth@\wdz@        
   \ifdim\rwidth@>\maxrwidth@\global\maxrwidth@\rwidth@\fi                  
   \global\advance\and@\@ne                                                
   &\Tag@
   \eat@{##}\crcr#1\crcr}}
 \totwidth@\maxlwidth@\advance\totwidth@\maxrwidth@}                       
\def\displ@y@{\global\dt@ptrue\openup\jot
 \everycr{\noalign{\global\tag@false\global\and@\z@\ifdt@p\global\dt@pfalse
 \vskip-\lineskiplimit\vskip\normallineskiplimit\else
 \penalty\interdisplaylinepenalty\fi}}}
\def\black@#1{\noalign{\ifdim#1>\displaywidth
 \dimen@\prevdepth\nointerlineskip                                          
 \vskip-\ht\strutbox@\vskip-\dp\strutbox@                                   
 \vbox{\noindent\hbox to#1{\strut@\hfill}}
 \prevdepth\dimen@                                                          
 \fi}}
\expandafter\def\csname align \space\endcsname#1\endalign
 {\measure@#1\endalign\global\and@\z@                                       
 \ifingather@\everycr{\noalign{\global\and@\z@}}\else\displ@y@\fi           
 \Let@\tabskip\centering@                                                   
 \halign to\displaywidth
  {\hfil\strut@\setboxz@h{$\m@th\displaystyle{\@lign##}$}
  \global\lwidth@\wdz@\boxz@\global\advance\and@\@ne                        
  \tabskip\z@skip                                                           
  &\setboxz@h{$\m@th\displaystyle{{}\@lign##}$}
  \global\rwidth@\wdz@\boxz@\hfill\global\advance\and@\@ne                  
  \tabskip\centering@                                                       
  &\setboxz@h{\@lign\strut@\maketag@##\maketag@}
  \dimen@\displaywidth\advance\dimen@-\totwidth@
  \divide\dimen@\tw@\advance\dimen@\maxrwidth@\advance\dimen@-\rwidth@     
  \ifdim\dimen@<\tw@\wdz@\llap{\vtop{\normalbaselines\null\boxz@}}
  \else\llap{\boxz@}\fi                                                    
  \tabskip\z@skip                                                          
  \crcr#1\crcr                                                             
  \black@\totwidth@}}                                                      
\newdimen\lineht@
\expandafter\def\csname align \endcsname#1\endalign{\measure@#1\endalign
 \global\and@\z@
 \ifdim\totwidth@>\displaywidth\let\displaywidth@\totwidth@\else
  \let\displaywidth@\displaywidth\fi                                        
 \ifingather@\everycr{\noalign{\global\and@\z@}}\else\displ@y@\fi
 \Let@\tabskip\centering@\halign to\displaywidth
  {\hfil\strut@\setboxz@h{$\m@th\displaystyle{\@lign##}$}%
  \global\lwidth@\wdz@\global\lineht@\ht\z@                                 
  \boxz@\global\advance\and@\@ne
  \tabskip\z@skip&\setboxz@h{$\m@th\displaystyle{{}\@lign##}$}%
  \global\rwidth@\wdz@\ifdim\ht\z@>\lineht@\global\lineht@\ht\z@\fi         
  \boxz@\hfil\global\advance\and@\@ne
  \tabskip\centering@&\kern-\displaywidth@                                  
  \setboxz@h{\@lign\strut@\maketag@##\maketag@}%
  \dimen@\displaywidth\advance\dimen@-\totwidth@
  \divide\dimen@\tw@\advance\dimen@\maxlwidth@\advance\dimen@-\lwidth@
  \ifdim\dimen@<\tw@\wdz@
   \rlap{\vbox{\normalbaselines\boxz@\vbox to\lineht@{}}}\else
   \rlap{\boxz@}\fi
  \tabskip\displaywidth@\crcr#1\crcr\black@\totwidth@}}
\expandafter\def\csname align (in \string\gather)\endcsname
  #1\endalign{\vcenter{\align@#1\endalign}}
\Invalid@\endalign
\newif\ifxat@
\def\alignat{\RIfMIfI@\DN@{\onlydmatherr@\alignat}\else
 \DN@{\csname alignat \endcsname}\fi\else
 \DN@{\onlydmatherr@\alignat}\fi\next@}
\newif\ifmeasuring@
\newbox\savealignat@
\expandafter\def\csname alignat \endcsname#1#2\endalignat                   
 {\inany@true\xat@false
 \def\tag{\global\tag@true\count@#1\relax\multiply\count@\tw@
  \xdef\tag@{}\loop\ifnum\count@>\and@\xdef\tag@{&\tag@}\advance\count@\m@ne
  \repeat\tag@}%
 \vspace@\allowdisplaybreak@\displaybreak@\intertext@
 \displ@y@\measuring@true                                                   
 \setbox\savealignat@\hbox{$\m@th\displaystyle\Let@
  \attag@{#1}
  \vbox{\halign{\span\preamble@@\crcr#2\crcr}}$}%
 \measuring@false                                                           
 \Let@\attag@{#1}
 \tabskip\centering@\halign to\displaywidth
  {\span\preamble@@\crcr#2\crcr                                             
  \black@{\wd\savealignat@}}}                                               
\Invalid@\endalignat
\def\xalignat{\RIfMIfI@
 \DN@{\onlydmatherr@\xalignat}\else
 \DN@{\csname xalignat \endcsname}\fi\else
 \DN@{\onlydmatherr@\xalignat}\fi\next@}
\expandafter\def\csname xalignat \endcsname#1#2\endxalignat
 {\inany@true\xat@true
 \def\tag{\global\tag@true\def\tag@{}\count@#1\relax\multiply\count@\tw@
  \loop\ifnum\count@>\and@\xdef\tag@{&\tag@}\advance\count@\m@ne\repeat\tag@}%
 \vspace@\allowdisplaybreak@\displaybreak@\intertext@
 \displ@y@\measuring@true\setbox\savealignat@\hbox{$\m@th\displaystyle\Let@
 \attag@{#1}\vbox{\halign{\span\preamble@@\crcr#2\crcr}}$}%
 \measuring@false\Let@
 \attag@{#1}\tabskip\centering@\halign to\displaywidth
 {\span\preamble@@\crcr#2\crcr\black@{\wd\savealignat@}}}
\def\attag@#1{\let\Maketag@\maketag@\let\TAG@\Tag@                          
 \let\Tag@=0\let\maketag@=0
 \ifmeasuring@\def\llap@##1{\setboxz@h{##1}\hbox to\tw@\wdz@{}}%
  \def\rlap@##1{\setboxz@h{##1}\hbox to\tw@\wdz@{}}\else
  \let\llap@\llap\let\rlap@\rlap\fi                                         
 \toks@{\hfil\strut@$\m@th\displaystyle{\@lign\the\hashtoks@}$\tabskip\z@skip
  \global\advance\and@\@ne&$\m@th\displaystyle{{}\@lign\the\hashtoks@}$\hfil
  \ifxat@\tabskip\centering@\fi\global\advance\and@\@ne}
 \iftagsleft@
  \toks@@{\tabskip\centering@&\Tag@\kern-\displaywidth
   \rlap@{\@lign\maketag@\the\hashtoks@\maketag@}%
   \global\advance\and@\@ne\tabskip\displaywidth}\else
  \toks@@{\tabskip\centering@&\Tag@\llap@{\@lign\maketag@
   \the\hashtoks@\maketag@}\global\advance\and@\@ne\tabskip\z@skip}\fi      
 \atcount@#1\relax\advance\atcount@\m@ne
 \loop\ifnum\atcount@>\z@
 \toks@=\expandafter{\the\toks@&\hfil$\m@th\displaystyle{\@lign
  \the\hashtoks@}$\global\advance\and@\@ne
  \tabskip\z@skip&$\m@th\displaystyle{{}\@lign\the\hashtoks@}$\hfil\ifxat@
  \tabskip\centering@\fi\global\advance\and@\@ne}\advance\atcount@\m@ne
 \repeat                                                                    
 \xdef\preamble@{\the\toks@\the\toks@@}
 \xdef\preamble@@{\preamble@}
 \let\maketag@\Maketag@\let\Tag@\TAG@}                                      
\Invalid@\endxalignat
\def\xxalignat{\RIfMIfI@
 \DN@{\onlydmatherr@\xxalignat}\else\DN@{\csname xxalignat
  \endcsname}\fi\else
 \DN@{\onlydmatherr@\xxalignat}\fi\next@}
\expandafter\def\csname xxalignat \endcsname#1#2\endxxalignat{\inany@true
 \vspace@\allowdisplaybreak@\displaybreak@\intertext@
 \displ@y\setbox\savealignat@\hbox{$\m@th\displaystyle\Let@
 \xxattag@{#1}\vbox{\halign{\span\preamble@@\crcr#2\crcr}}$}%
 \Let@\xxattag@{#1}\tabskip\z@skip\halign to\displaywidth
 {\span\preamble@@\crcr#2\crcr\black@{\wd\savealignat@}}}
\def\xxattag@#1{\toks@{\tabskip\z@skip\hfil\strut@
 $\m@th\displaystyle{\the\hashtoks@}$&%
 $\m@th\displaystyle{{}\the\hashtoks@}$\hfil\tabskip\centering@&}%
 \atcount@#1\relax\advance\atcount@\m@ne\loop\ifnum\atcount@>\z@
 \toks@=\expandafter{\the\toks@&\hfil$\m@th\displaystyle{\the\hashtoks@}$%
  \tabskip\z@skip&$\m@th\displaystyle{{}\the\hashtoks@}$\hfil
  \tabskip\centering@}\advance\atcount@\m@ne\repeat
 \xdef\preamble@{\the\toks@\tabskip\z@skip}\xdef\preamble@@{\preamble@}}
\Invalid@\endxxalignat
\newdimen\gwidth@
\newdimen\gmaxwidth@
\def\gmeasure@#1\endgather{\gwidth@\z@\gmaxwidth@\z@\setbox@ne\vbox{\Let@
 \halign{\setboxz@h{$\m@th\displaystyle{##}$}\global\gwidth@\wdz@
 \ifdim\gwidth@>\gmaxwidth@\global\gmaxwidth@\gwidth@\fi
 &\eat@{##}\crcr#1\crcr}}}
\def\gather{\RIfMIfI@\DN@{\onlydmatherr@\gather}\else
 \ingather@true\inany@true\def\tag{&}%
 \vspace@\allowdisplaybreak@\displaybreak@\intertext@
 \displ@y\Let@
 \iftagsleft@\DN@{\csname gather \endcsname}\else
  \DN@{\csname gather \space\endcsname}\fi\fi
 \else\DN@{\onlydmatherr@\gather}\fi\next@}
\expandafter\def\csname gather \space\endcsname#1\endgather
 {\gmeasure@#1\endgather\tabskip\centering@
 \halign to\displaywidth{\hfil\strut@\setboxz@h{$\m@th\displaystyle{##}$}%
 \global\gwidth@\wdz@\boxz@\hfil&
 \setboxz@h{\strut@{\maketag@##\maketag@}}%
 \dimen@\displaywidth\advance\dimen@-\gwidth@
 \ifdim\dimen@>\tw@\wdz@\llap{\boxz@}\else
 \llap{\vtop{\normalbaselines\null\boxz@}}\fi
 \tabskip\z@skip\crcr#1\crcr\black@\gmaxwidth@}}
\newdimen\glineht@
\expandafter\def\csname gather \endcsname#1\endgather{\gmeasure@#1\endgather
 \ifdim\gmaxwidth@>\displaywidth\let\gdisplaywidth@\gmaxwidth@\else
 \let\gdisplaywidth@\displaywidth\fi\tabskip\centering@\halign to\displaywidth
 {\hfil\strut@\setboxz@h{$\m@th\displaystyle{##}$}%
 \global\gwidth@\wdz@\global\glineht@\ht\z@\boxz@\hfil&\kern-\gdisplaywidth@
 \setboxz@h{\strut@{\maketag@##\maketag@}}%
 \dimen@\displaywidth\advance\dimen@-\gwidth@
 \ifdim\dimen@>\tw@\wdz@\rlap{\boxz@}\else
 \rlap{\vbox{\normalbaselines\boxz@\vbox to\glineht@{}}}\fi
 \tabskip\gdisplaywidth@\crcr#1\crcr\black@\gmaxwidth@}}
\newif\ifctagsplit@
\def\CenteredTagsOnSplits{\global\ctagsplit@true}
\def\TopOrBottomTagsOnSplits{\global\ctagsplit@false}
\TopOrBottomTagsOnSplits
\def\split{\relax\ifinany@\let\next@\insplit@\else
 \ifmmode\ifinner\def\next@{\onlydmatherr@\split}\else
 \let\next@\outsplit@\fi\else
 \def\next@{\onlydmatherr@\split}\fi\fi\next@}
\def\insplit@{\global\setbox\z@\vbox\bgroup\vspace@\Let@\ialign\bgroup
 \hfil\strut@$\m@th\displaystyle{##}$&$\m@th\displaystyle{{}##}$\hfill\crcr}
\def\endsplit{\crcr\egroup\egroup\iftagsleft@\expandafter\lendsplit@\else
 \expandafter\rendsplit@\fi}
\def\rendsplit@{\global\setbox9 \vbox
 {\unvcopy\z@\global\setbox8 \lastbox\unskip}
 \setbox@ne\hbox{\unhcopy8 \unskip\global\setbox\tw@\lastbox
 \unskip\global\setbox\thr@@\lastbox}
 \global\setbox7 \hbox{\unhbox\tw@\unskip}
 \ifinalign@\ifctagsplit@                                                   
  \gdef\split@{\hbox to\wd\thr@@{}&
   \vcenter{\vbox{\moveleft\wd\thr@@\boxz@}}}
 \else\gdef\split@{&\vbox{\moveleft\wd\thr@@\box9}\crcr
  \box\thr@@&\box7}\fi                                                      
 \else                                                                      
  \ifctagsplit@\gdef\split@{\vcenter{\boxz@}}\else
  \gdef\split@{\box9\crcr\hbox{\box\thr@@\box7}}\fi
 \fi
 \split@}                                                                   
\def\lendsplit@{\global\setbox9\vtop{\unvcopy\z@}
 \setbox@ne\vbox{\unvcopy\z@\global\setbox8\lastbox}
 \setbox@ne\hbox{\unhcopy8\unskip\setbox\tw@\lastbox
  \unskip\global\setbox\thr@@\lastbox}
 \ifinalign@\ifctagsplit@                                                   
  \gdef\split@{\hbox to\wd\thr@@{}&
  \vcenter{\vbox{\moveleft\wd\thr@@\box9}}}
  \else                                                                     
  \gdef\split@{\hbox to\wd\thr@@{}&\vbox{\moveleft\wd\thr@@\box9}}\fi
 \else
  \ifctagsplit@\gdef\split@{\vcenter{\box9}}\else
  \gdef\split@{\box9}\fi
 \fi\split@}
\def\outsplit@#1$${\align\insplit@#1\endalign$$}
\newdimen\multlinegap@
\multlinegap@1em
\newdimen\multlinetaggap@
\multlinetaggap@1em
\def\MultlineGap#1{\global\multlinegap@#1\relax}
\def\multlinegap#1{\RIfMIfI@\onlydmatherr@\multlinegap\else
 \multlinegap@#1\relax\fi\else\onlydmatherr@\multlinegap\fi}
\def\nomultlinegap{\multlinegap{\z@}}
\def\multline{\RIfMIfI@
 \DN@{\onlydmatherr@\multline}\else
 \DN@{\multline@}\fi\else
 \DN@{\onlydmatherr@\multline}\fi\next@}
\newif\iftagin@
\def\tagin@#1{\tagin@false\in@\tag{#1}\ifin@\tagin@true\fi}
\def\multline@#1$${\inany@true\vspace@\allowdisplaybreak@\displaybreak@
 \tagin@{#1}\iftagsleft@\DN@{\multline@l#1$$}\else
 \DN@{\multline@r#1$$}\fi\next@}
\newdimen\mwidth@
\def\rmmeasure@#1\endmultline{%
 \def\shoveleft##1{##1}\def\shoveright##1{##1}
 \setbox@ne\vbox{\Let@\halign{\setboxz@h
  {$\m@th\@lign\displaystyle{}##$}\global\mwidth@\wdz@
  \crcr#1\crcr}}}
\newdimen\mlineht@
\newif\ifzerocr@
\newif\ifonecr@
\def\lmmeasure@#1\endmultline{\global\zerocr@true\global\onecr@false
 \everycr{\noalign{\ifonecr@\global\onecr@false\fi
  \ifzerocr@\global\zerocr@false\global\onecr@true\fi}}
  \def\shoveleft##1{##1}\def\shoveright##1{##1}%
 \setbox@ne\vbox{\Let@\halign{\setboxz@h
  {$\m@th\@lign\displaystyle{}##$}\ifonecr@\global\mwidth@\wdz@
  \global\mlineht@\ht\z@\fi\crcr#1\crcr}}}
\newbox\mtagbox@
\newdimen\ltwidth@
\newdimen\rtwidth@
\def\multline@l#1$${\iftagin@\DN@{\lmultline@@#1$$}\else
 \DN@{\setbox\mtagbox@\null\ltwidth@\z@\rtwidth@\z@
  \lmultline@@@#1$$}\fi\next@}
\def\lmultline@@#1\endmultline\tag#2$${%
 \setbox\mtagbox@\hbox{\maketag@#2\maketag@}
 \lmmeasure@#1\endmultline\dimen@\mwidth@\advance\dimen@\wd\mtagbox@
 \advance\dimen@\multlinetaggap@                                            
 \ifdim\dimen@>\displaywidth\ltwidth@\z@\else\ltwidth@\wd\mtagbox@\fi       
 \lmultline@@@#1\endmultline$$}
\def\lmultline@@@{\displ@y
 \def\shoveright##1{##1\hfilneg\hskip\multlinegap@}%
 \def\shoveleft##1{\setboxz@h{$\m@th\displaystyle{}##1$}%
  \setbox@ne\hbox{$\m@th\displaystyle##1$}%
  \hfilneg
  \iftagin@
   \ifdim\ltwidth@>\z@\hskip\ltwidth@\hskip\multlinetaggap@\fi
  \else\hskip\multlinegap@\fi\hskip.5\wd@ne\hskip-.5\wdz@##1}
  \halign\bgroup\Let@\hbox to\displaywidth
   {\strut@$\m@th\displaystyle\hfil{}##\hfil$}\crcr
   \hfilneg                                                                 
   \iftagin@                                                                
    \ifdim\ltwidth@>\z@                                                     
     \box\mtagbox@\hskip\multlinetaggap@                                    
    \else
     \rlap{\vbox{\normalbaselines\hbox{\strut@\box\mtagbox@}%
     \vbox to\mlineht@{}}}\fi                                               
   \else\hskip\multlinegap@\fi}                                             
\def\multline@r#1$${\iftagin@\DN@{\rmultline@@#1$$}\else
 \DN@{\setbox\mtagbox@\null\ltwidth@\z@\rtwidth@\z@
  \rmultline@@@#1$$}\fi\next@}
\def\rmultline@@#1\endmultline\tag#2$${\ltwidth@\z@
 \setbox\mtagbox@\hbox{\maketag@#2\maketag@}%
 \rmmeasure@#1\endmultline\dimen@\mwidth@\advance\dimen@\wd\mtagbox@
 \advance\dimen@\multlinetaggap@
 \ifdim\dimen@>\displaywidth\rtwidth@\z@\else\rtwidth@\wd\mtagbox@\fi
 \rmultline@@@#1\endmultline$$}
\def\rmultline@@@{\displ@y
 \def\shoveright##1{##1\hfilneg\iftagin@\ifdim\rtwidth@>\z@
  \hskip\rtwidth@\hskip\multlinetaggap@\fi\else\hskip\multlinegap@\fi}%
 \def\shoveleft##1{\setboxz@h{$\m@th\displaystyle{}##1$}%
  \setbox@ne\hbox{$\m@th\displaystyle##1$}%
  \hfilneg\hskip\multlinegap@\hskip.5\wd@ne\hskip-.5\wdz@##1}%
 \halign\bgroup\Let@\hbox to\displaywidth
  {\strut@$\m@th\displaystyle\hfil{}##\hfil$}\crcr
 \hfilneg\hskip\multlinegap@}
\def\endmultline{\iftagsleft@\expandafter\lendmultline@\else
 \expandafter\rendmultline@\fi}
\def\lendmultline@{\hfilneg\hskip\multlinegap@\crcr\egroup}
\def\rendmultline@{\iftagin@                                                
 \ifdim\rtwidth@>\z@                                                        
  \hskip\multlinetaggap@\box\mtagbox@                                       
 \else\llap{\vtop{\normalbaselines\null\hbox{\strut@\box\mtagbox@}}}\fi     
 \else\hskip\multlinegap@\fi                                                
 \hfilneg\crcr\egroup}
\def\bmod{\mskip-\medmuskip\mkern5mu\mathbin{\fam\z@ mod}\penalty900
 \mkern5mu\mskip-\medmuskip}
\def\pmod#1{\allowbreak\ifinner\mkern8mu\else\mkern18mu\fi
 ({\fam\z@ mod}\,\,#1)}
\def\pod#1{\allowbreak\ifinner\mkern8mu\else\mkern18mu\fi(#1)}
\def\mod#1{\allowbreak\ifinner\mkern12mu\else\mkern18mu\fi{\fam\z@ mod}\,\,#1}
\message{continued fractions,}
\newcount\cfraccount@
\def\cfrac{\bgroup\bgroup\advance\cfraccount@\@ne\strut
 \iffalse{\fi\def\\{\over\displaystyle}\iffalse}\fi}
\def\lcfrac{\bgroup\bgroup\advance\cfraccount@\@ne\strut
 \iffalse{\fi\def\\{\hfill\over\displaystyle}\iffalse}\fi}
\def\rcfrac{\bgroup\bgroup\advance\cfraccount@\@ne\strut\hfill
 \iffalse{\fi\def\\{\over\displaystyle}\iffalse}\fi}
\def\gloop@#1\repeat{\gdef\body{#1}\iterate}
\def\endcfrac{\gloop@\ifnum\cfraccount@>\z@\global\advance\cfraccount@\m@ne
 \egroup\hskip-\nulldelimiterspace\egroup\repeat}
\message{compound symbols,}
\def\binrel@#1{\setboxz@h{\thinmuskip0mu
  \medmuskip\m@ne mu\thickmuskip\@ne mu$#1\m@th$}%
 \setbox@ne\hbox{\thinmuskip0mu\medmuskip\m@ne mu\thickmuskip
  \@ne mu${}#1{}\m@th$}%
 \setbox\tw@\hbox{\hskip\wd@ne\hskip-\wdz@}}
\def\overset#1\to#2{\binrel@{#2}\ifdim\wd\tw@<\z@
 \mathbin{\mathop{\kern\z@#2}\limits^{#1}}\else\ifdim\wd\tw@>\z@
 \mathrel{\mathop{\kern\z@#2}\limits^{#1}}\else
 {\mathop{\kern\z@#2}\limits^{#1}}{}\fi\fi}
\def\underset#1\to#2{\binrel@{#2}\ifdim\wd\tw@<\z@
 \mathbin{\mathop{\kern\z@#2}\limits_{#1}}\else\ifdim\wd\tw@>\z@
 \mathrel{\mathop{\kern\z@#2}\limits_{#1}}\else
 {\mathop{\kern\z@#2}\limits_{#1}}{}\fi\fi}
\def\oversetbrace#1\to#2{\overbrace{#2}^{#1}}
\def\undersetbrace#1\to#2{\underbrace{#2}_{#1}}
\def\sideset#1\and#2\to#3{%
 \setbox@ne\hbox{$\dsize{\vphantom{#3}}#1{#3}\m@th$}%
 \setbox\tw@\hbox{$\dsize{#3}#2\m@th$}%
 \hskip\wd@ne\hskip-\wd\tw@\mathop{\hskip\wd\tw@\hskip-\wd@ne
  {\vphantom{#3}}#1{#3}#2}}
\def\rightarrowfill@#1{\setboxz@h{$#1-\m@th$}\ht\z@\z@
  $#1\m@th\copy\z@\mkern-6mu\cleaders
  \hbox{$#1\mkern-2mu\box\z@\mkern-2mu$}\hfill
  \mkern-6mu\mathord\rightarrow$}
\def\leftarrowfill@#1{\setboxz@h{$#1-\m@th$}\ht\z@\z@
  $#1\m@th\mathord\leftarrow\mkern-6mu\cleaders
  \hbox{$#1\mkern-2mu\copy\z@\mkern-2mu$}\hfill
  \mkern-6mu\box\z@$}
\def\leftrightarrowfill@#1{\setboxz@h{$#1-\m@th$}\ht\z@\z@
  $#1\m@th\mathord\leftarrow\mkern-6mu\cleaders
  \hbox{$#1\mkern-2mu\box\z@\mkern-2mu$}\hfill
  \mkern-6mu\mathord\rightarrow$}
\def\overrightarrow{\mathpalette\overrightarrow@}
\def\overrightarrow@#1#2{\vbox{\ialign{##\crcr\rightarrowfill@#1\crcr
 \noalign{\kern-\ex@\nointerlineskip}$\m@th\hfil#1#2\hfil$\crcr}}}

\def\overleftarrow{\mathpalette\overleftarrow@}
\def\overleftarrow@#1#2{\vbox{\ialign{##\crcr\leftarrowfill@#1\crcr
 \noalign{\kern-\ex@\nointerlineskip}$\m@th\hfil#1#2\hfil$\crcr}}}
\def\overleftrightarrow{\mathpalette\overleftrightarrow@}
\def\overleftrightarrow@#1#2{\vbox{\ialign{##\crcr\leftrightarrowfill@#1\crcr
 \noalign{\kern-\ex@\nointerlineskip}$\m@th\hfil#1#2\hfil$\crcr}}}
\def\underrightarrow{\mathpalette\underrightarrow@}
\def\underrightarrow@#1#2{\vtop{\ialign{##\crcr$\m@th\hfil#1#2\hfil$\crcr
 \noalign{\nointerlineskip}\rightarrowfill@#1\crcr}}}

\def\underleftarrow{\mathpalette\underleftarrow@}
\def\underleftarrow@#1#2{\vtop{\ialign{##\crcr$\m@th\hfil#1#2\hfil$\crcr
 \noalign{\nointerlineskip}\leftarrowfill@#1\crcr}}}
\def\underleftrightarrow{\mathpalette\underleftrightarrow@}
\def\underleftrightarrow@#1#2{\vtop{\ialign{##\crcr$\m@th\hfil#1#2\hfil$\crcr
 \noalign{\nointerlineskip}\leftrightarrowfill@#1\crcr}}}
\message{various kinds of dots,}
\let\DOTSI\relax
\let\DOTSB\relax

\newif\ifmath@
{\uccode`7=`\\ \uccode`8=`m \uccode`9=`a \uccode`0=`t \uccode`!=`h
 \uppercase{\gdef\math@#1#2#3#4#5#6\math@{\global\math@false\ifx 7#1\ifx 8#2%
 \ifx 9#3\ifx 0#4\ifx !#5\xdef\meaning@{#6}\global\math@true\fi\fi\fi\fi\fi}}}
\newif\ifmathch@
{\uccode`7=`c \uccode`8=`h \uccode`9=`\"
 \uppercase{\gdef\mathch@#1#2#3#4#5#6\mathch@{\global\mathch@false
  \ifx 7#1\ifx 8#2\ifx 9#5\global\mathch@true\xdef\meaning@{9#6}\fi\fi\fi}}}
\newcount\classnum@
\def\getmathch@#1.#2\getmathch@{\classnum@#1 \divide\classnum@4096
 \ifcase\number\classnum@\or\or\gdef\thedots@{\dotsb@}\or
 \gdef\thedots@{\dotsb@}\fi}
\newif\ifmathbin@
{\uccode`4=`b \uccode`5=`i \uccode`6=`n
 \uppercase{\gdef\mathbin@#1#2#3{\relaxnext@
  \DNii@##1\mathbin@{\ifx\space@\next\global\mathbin@true\fi}%
 \global\mathbin@false\DN@##1\mathbin@{}%
 \ifx 4#1\ifx 5#2\ifx 6#3\DN@{\FN@\nextii@}\fi\fi\fi\next@}}}
\newif\ifmathrel@
{\uccode`4=`r \uccode`5=`e \uccode`6=`l
 \uppercase{\gdef\mathrel@#1#2#3{\relaxnext@
  \DNii@##1\mathrel@{\ifx\space@\next\global\mathrel@true\fi}%
 \global\mathrel@false\DN@##1\mathrel@{}%
 \ifx 4#1\ifx 5#2\ifx 6#3\DN@{\FN@\nextii@}\fi\fi\fi\next@}}}
\newif\ifmacro@
{\uccode`5=`m \uccode`6=`a \uccode`7=`c
 \uppercase{\gdef\macro@#1#2#3#4\macro@{\global\macro@false
  \ifx 5#1\ifx 6#2\ifx 7#3\global\macro@true
  \xdef\meaning@{\macro@@#4\macro@@}\fi\fi\fi}}}
\def\macro@@#1->#2\macro@@{#2}
\newif\ifDOTS@
\newcount\DOTSCASE@
{\uccode`6=`\\ \uccode`7=`D \uccode`8=`O \uccode`9=`T \uccode`0=`S
 \uppercase{\gdef\DOTS@#1#2#3#4#5{\global\DOTS@false\DN@##1\DOTS@{}%
  \ifx 6#1\ifx 7#2\ifx 8#3\ifx 9#4\ifx 0#5\let\next@\DOTS@@\fi\fi\fi\fi\fi
  \next@}}}
{\uccode`3=`B \uccode`4=`I \uccode`5=`X
 \uppercase{\gdef\DOTS@@#1{\relaxnext@
  \DNii@##1\DOTS@{\ifx\space@\next\global\DOTS@true\fi}%
  \DN@{\FN@\nextii@}%
  \ifx 3#1\global\DOTSCASE@\z@\else
  \ifx 4#1\global\DOTSCASE@\@ne\else
  \ifx 5#1\global\DOTSCASE@\tw@\else\DN@##1\DOTS@{}%
  \fi\fi\fi\next@}}}
\newif\ifnot@
{\uccode`5=`\\ \uccode`6=`n \uccode`7=`o \uccode`8=`t
 \uppercase{\gdef\not@#1#2#3#4{\relaxnext@
  \DNii@##1\not@{\ifx\space@\next\global\not@true\fi}%
 \global\not@false\DN@##1\not@{}%
 \ifx 5#1\ifx 6#2\ifx 7#3\ifx 8#4\DN@{\FN@\nextii@}\fi\fi\fi
 \fi\next@}}}
\newif\ifkeybin@
\def\keybin@{\keybin@true
 \ifx\next+\else\ifx\next=\else\ifx\next<\else\ifx\next>\else\ifx\next-\else
 \ifx\next*\else\ifx\next:\else\keybin@false\fi\fi\fi\fi\fi\fi\fi}
\def\dots{\RIfM@\expandafter\mdots@\else\expandafter\tdots@\fi}
\def\tdots@{\unskip\relaxnext@
 \DN@{$\m@th\mathinner{\ldotp\ldotp\ldotp}\,
   \ifx\next,\,$\else\ifx\next.\,$\else\ifx\next;\,$\else\ifx\next:\,$\else
   \ifx\next?\,$\else\ifx\next!\,$\else$ \fi\fi\fi\fi\fi\fi}%
 \ \FN@\next@}
\def\mdots@{\FN@\mdots@@}
\def\mdots@@{\gdef\thedots@{\dotso@}
 \ifx\next\boldkey\gdef\thedots@\boldkey{\boldkeydots@}\else                
 \ifx\next\boldsymbol\gdef\thedots@\boldsymbol{\boldsymboldots@}\else       
 \ifx,\next\gdef\thedots@{\dotsc}
 \else\ifx\not\next\gdef\thedots@{\dotsb@}
 \else\keybin@
 \ifkeybin@\gdef\thedots@{\dotsb@}
 \else\xdef\meaning@{\meaning\next..........}\xdef\meaning@@{\meaning@}
  \expandafter\math@\meaning@\math@
  \ifmath@
   \expandafter\mathch@\meaning@\mathch@
   \ifmathch@\expandafter\getmathch@\meaning@\getmathch@\fi                 
  \else\expandafter\macro@\meaning@@\macro@                                 
  \ifmacro@                                                                
   \expandafter\not@\meaning@\not@\ifnot@\gdef\thedots@{\dotsb@}
  \else\expandafter\DOTS@\meaning@\DOTS@
  \ifDOTS@
   \ifcase\number\DOTSCASE@\gdef\thedots@{\dotsb@}%
    \or\gdef\thedots@{\dotsi}\else\fi                                      
  \else\expandafter\math@\meaning@\math@                                   
  \ifmath@\expandafter\mathbin@\meaning@\mathbin@
  \ifmathbin@\gdef\thedots@{\dotsb@}
  \else\expandafter\mathrel@\meaning@\mathrel@
  \ifmathrel@\gdef\thedots@{\dotsb@}
  \fi\fi\fi\fi\fi\fi\fi\fi\fi\fi\fi\fi
 \thedots@}
\def\plainldots@{\mathinner{\ldotp\ldotp\ldotp}}
\def\plaincdots@{\mathinner{\cdotp\cdotp\cdotp}}
\def\dotsi{\!\plaincdots@}
\let\dotsb@\plaincdots@
\newif\ifextra@
\newif\ifrightdelim@
\def\rightdelim@{\global\rightdelim@true                                    
 \ifx\next)\else                                                            
 \ifx\next]\else
 \ifx\next\rbrack\else
 \ifx\next\}\else
 \ifx\next\rbrace\else
 \ifx\next\rangle\else
 \ifx\next\rceil\else
 \ifx\next\rfloor\else
 \ifx\next\rgroup\else
 \ifx\next\rmoustache\else
 \ifx\next\right\else
 \ifx\next\bigr\else
 \ifx\next\biggr\else
 \ifx\next\Bigr\else                                                        
 \ifx\next\Biggr\else\global\rightdelim@false
 \fi\fi\fi\fi\fi\fi\fi\fi\fi\fi\fi\fi\fi\fi\fi}
\def\extra@{%
 \global\extra@false\rightdelim@\ifrightdelim@\global\extra@true            
 \else\ifx\next$\global\extra@true                                          
 \else\xdef\meaning@{\meaning\next..........}
 \expandafter\macro@\meaning@\macro@\ifmacro@                               
 \expandafter\DOTS@\meaning@\DOTS@
 \ifDOTS@
 \ifnum\DOTSCASE@=\tw@\global\extra@true                                    
 \fi\fi\fi\fi\fi}
\newif\ifbold@
\def\dotso@{\relaxnext@
 \ifbold@
  \let\next\delayed@
  \DNii@{\extra@\plainldots@\ifextra@\,\fi}%
 \else
  \DNii@{\DN@{\extra@\plainldots@\ifextra@\,\fi}\FN@\next@}%
 \fi
 \nextii@}
\def\extrap@#1{%
 \ifx\next,\DN@{#1\,}\else
 \ifx\next;\DN@{#1\,}\else
 \ifx\next.\DN@{#1\,}\else\extra@
 \ifextra@\DN@{#1\,}\else
 \let\next@#1\fi\fi\fi\fi\next@}
\def\ldots{\DN@{\extrap@\plainldots@}%
 \FN@\next@}
\def\cdots{\DN@{\extrap@\plaincdots@}%
 \FN@\next@}

\def\dotsc{\relaxnext@
 \DN@{\ifx\next;\plainldots@\,\else
  \ifx\next.\plainldots@\,\else\extra@\plainldots@
  \ifextra@\,\fi\fi\fi}%
 \FN@\next@}
\def\cdot{\mathchar"2201 }

\def\mapsto{\DOTSB\mapstochar\rightarrow}

\message{special superscripts,}
\def\dddot#1{{\mathop{#1}\limits^{\vbox to-1.4\ex@{\kern-\tw@\ex@
 \hbox{\rm...}\vss}}}}
\def\ddddot#1{{\mathop{#1}\limits^{\vbox to-1.4\ex@{\kern-\tw@\ex@
 \hbox{\rm....}\vss}}}}
\def\sphat{^{\mathchoice{}{}%
 {\,\,\botsmash{\hbox{\lower4\ex@\hbox{$\m@th\widehat{\null}$}}}}%
 {\,\botsmash{\hbox{\lower3\ex@\hbox{$\m@th\hat{\null}$}}}}}}

\def\spacute{^{\!\botsmash{\hbox{\lower\@ne ex\hbox{\'{}}}}}}
\def\spgrave{^{\mathchoice{}{}{}{\!}%
 \botsmash{\hbox{\lower\@ne ex\hbox{\`{}}}}}}
\def\spdot{^{\hbox{\raise\ex@\hbox{\rm.}}}}
\def\spddot{^{\hbox{\raise\ex@\hbox{\rm..}}}}
\def\spdddot{^{\hbox{\raise\ex@\hbox{\rm...}}}}
\def\spddddot{^{\hbox{\raise\ex@\hbox{\rm....}}}}
\def\spbreve{^{\!\botsmash{\hbox{\lower4\ex@\hbox{\u{}}}}}}

\message{\string\text,}
\def\textonlyfont@#1#2{\def#1{\RIfM@
 \Err@{Use \string#1\space only in text}\else#2\fi}}
\textonlyfont@\rm\tenrm
\textonlyfont@\it\tenit
\textonlyfont@\sl\tensl
\textonlyfont@\bf\tenbf
\def\oldnos#1{\RIfM@{\mathcode`\,="013B \fam\@ne#1}\else
 \leavevmode\hbox{$\m@th\mathcode`\,="013B \fam\@ne#1$}\fi}
\def\text{\RIfM@\expandafter\text@\else\expandafter\text@@\fi}
\def\text@@#1{\leavevmode\hbox{#1}}
\def\mathhexbox@#1#2#3{\text{$\m@th\mathchar"#1#2#3$}}
\def\dag{{\mathhexbox@279}}
\def\ddag{{\mathhexbox@27A}}
\def\S{{\mathhexbox@278}}
\def\P{{\mathhexbox@27B}}
\newif\iffirstchoice@
\firstchoice@true
\def\text@#1{\mathchoice
 {\hbox{\everymath{\displaystyle}\def\textfonti{\the\textfont\@ne}%
  \def\textfontii{\the\textfont\tw@}\textdef@@ T#1}}
 {\hbox{\firstchoice@false
  \everymath{\textstyle}\def\textfonti{\the\textfont\@ne}%
  \def\textfontii{\the\textfont\tw@}\textdef@@ T#1}}
 {\hbox{\firstchoice@false
  \everymath{\scriptstyle}\def\textfonti{\the\scriptfont\@ne}%
  \def\textfontii{\the\scriptfont\tw@}\textdef@@ S\rm#1}}
 {\hbox{\firstchoice@false
  \everymath{\scriptscriptstyle}\def\textfonti
  {\the\scriptscriptfont\@ne}%
  \def\textfontii{\the\scriptscriptfont\tw@}\textdef@@ s\rm#1}}}
\def\textdef@@#1{\textdef@#1\rm\textdef@#1\bf\textdef@#1\sl\textdef@#1\it}
\def\rmfam{0}
\def\textdef@#1#2{%
 \DN@{\csname\expandafter\eat@\string#2fam\endcsname}%
 \if S#1\edef#2{\the\scriptfont\next@\relax}%
 \else\if s#1\edef#2{\the\scriptscriptfont\next@\relax}%
 \else\edef#2{\the\textfont\next@\relax}\fi\fi}
\scriptfont\itfam\tenit \scriptscriptfont\itfam\tenit
\scriptfont\slfam\tensl \scriptscriptfont\slfam\tensl
\newif\iftopfolded@
\newif\ifbotfolded@
\def\topfoldedtext{\topfolded@true\botfolded@false\foldedtext@}
\def\botfoldedtext{\botfolded@true\topfolded@false\foldedtext@}
\def\foldedtext{\topfolded@false\botfolded@false\foldedtext@}
\Invalid@\foldedwidth
\def\foldedtext@{\relaxnext@
 \DN@{\ifx\next\foldedwidth\let\next@\nextii@\else
  \DN@{\nextii@\foldedwidth{.3\hsize}}\fi\next@}%
 \DNii@\foldedwidth##1##2{\setbox\z@\vbox
  {\normalbaselines\hsize##1\relax
  \tolerance1600 \noindent\ignorespaces##2}\ifbotfolded@\boxz@\else
  \iftopfolded@\vtop{\unvbox\z@}\else\vcenter{\boxz@}\fi\fi}%
 \FN@\next@}
\message{math font commands,}
\def\bold{\RIfM@\expandafter\bold@\else
 \expandafter\nonmatherr@\expandafter\bold\fi}
\def\bold@#1{{\bold@@{#1}}}
\def\bold@@#1{\fam\bffam\relax#1}
\def\slanted{\RIfM@\expandafter\slanted@\else
 \expandafter\nonmatherr@\expandafter\slanted\fi}
\def\slanted@#1{{\slanted@@{#1}}}
\def\slanted@@#1{\fam\slfam\relax#1}
\def\roman{\RIfM@\expandafter\roman@\else
 \expandafter\nonmatherr@\expandafter\roman\fi}
\def\roman@#1{{\roman@@{#1}}}
\def\roman@@#1{\fam\rmfam\relax#1}
\def\italic{\RIfM@\expandafter\italic@\else
 \expandafter\nonmatherr@\expandafter\italic\fi}
\def\italic@#1{{\italic@@{#1}}}
\def\italic@@#1{\fam\itfam\relax#1}
\def\Cal{\RIfM@\expandafter\Cal@\else
 \expandafter\nonmatherr@\expandafter\Cal\fi}
\def\Cal@#1{{\Cal@@{#1}}}
\def\Cal@@#1{\noaccents@\fam\tw@#1}
\mathchardef\Gamma="0000
\mathchardef\Delta="0001
\mathchardef\Theta="0002
\mathchardef\Lambda="0003
\mathchardef\Xi="0004
\mathchardef\Pi="0005
\mathchardef\Sigma="0006
\mathchardef\Upsilon="0007
\mathchardef\Phi="0008
\mathchardef\Psi="0009
\mathchardef\Omega="000A
\mathchardef\varGamma="0100
\mathchardef\varDelta="0101
\mathchardef\varTheta="0102
\mathchardef\varLambda="0103
\mathchardef\varXi="0104
\mathchardef\varPi="0105
\mathchardef\varSigma="0106
\mathchardef\varUpsilon="0107
\mathchardef\varPhi="0108
\mathchardef\varPsi="0109
\mathchardef\varOmega="010A
\let\alloc@@\alloc@
\def\hexnumber@#1{\ifcase#1 0\or 1\or 2\or 3\or 4\or 5\or 6\or 7\or 8\or
 9\or A\or B\or C\or D\or E\or F\fi}
\def\loadmsam{%
 \font@\tenmsa=msam10
 \font@\sevenmsa=msam7
 \font@\fivemsa=msam5
 \alloc@@8\fam\chardef\sixt@@n\msafam
 \textfont\msafam=\tenmsa
 \scriptfont\msafam=\sevenmsa
 \scriptscriptfont\msafam=\fivemsa
 \edef\next{\hexnumber@\msafam}%
 \mathchardef\dabar@"0\next39
 \edef\dashrightarrow{\mathrel{\dabar@\dabar@\mathchar"0\next4B}}%
 \edef\dashleftarrow{\mathrel{\mathchar"0\next4C\dabar@\dabar@}}%
 \let\dasharrow\dashrightarrow
 \edef\ulcorner{\delimiter"4\next70\next70 }%
 \edef\urcorner{\delimiter"5\next71\next71 }%
 \edef\llcorner{\delimiter"4\next78\next78 }%
 \edef\lrcorner{\delimiter"5\next79\next79 }%
 \edef\yen{{\noexpand\mathhexbox@\next55}}%
 \edef\checkmark{{\noexpand\mathhexbox@\next58}}%
 \edef\circledR{{\noexpand\mathhexbox@\next72}}%
 \edef\maltese{{\noexpand\mathhexbox@\next7A}}%
 \global\let\loadmsam\empty}%
\def\loadmsbm{%
 \font@\tenmsb=msbm10 \font@\sevenmsb=msbm7 \font@\fivemsb=msbm5
 \alloc@@8\fam\chardef\sixt@@n\msbfam
 \textfont\msbfam=\tenmsb
 \scriptfont\msbfam=\sevenmsb \scriptscriptfont\msbfam=\fivemsb
 \global\let\loadmsbm\empty
 }
\def\widehat#1{\ifx\undefined\msbfam \DN@{362}%
  \else \setboxz@h{$\m@th#1$}%
    \edef\next@{\ifdim\wdz@>\tw@ em%
        \hexnumber@\msbfam 5B%
      \else 362\fi}\fi
  \mathaccent"0\next@{#1}}
\def\widetilde#1{\ifx\undefined\msbfam \DN@{365}%
  \else \setboxz@h{$\m@th#1$}%
    \edef\next@{\ifdim\wdz@>\tw@ em%
        \hexnumber@\msbfam 5D%
      \else 365\fi}\fi
  \mathaccent"0\next@{#1}}
\message{\string\newsymbol,}
\def\newsymbol#1#2#3#4#5{\define#1{}%
  \count@#2\relax \advance\count@\m@ne 
 \ifcase\count@
   \ifx\undefined\msafam\loadmsam\fi \let\next@\msafam
 \or \ifx\undefined\msbfam\loadmsbm\fi \let\next@\msbfam
 \else  \Err@{\Invalid@@\string\newsymbol}\let\next@\tw@\fi
 \mathchardef#1="#3\hexnumber@\next@#4#5\space}

\def\Bbb{\RIfM@\expandafter\Bbb@\else
 \expandafter\nonmatherr@\expandafter\Bbb\fi}
\def\Bbb@#1{{\Bbb@@{#1}}}
\def\Bbb@@#1{\noaccents@\fam\msbfam\relax#1}
\message{bold Greek and bold symbols,}
\def\loadbold{%
 \font@\tencmmib=cmmib10 \font@\sevencmmib=cmmib7 \font@\fivecmmib=cmmib5
 \skewchar\tencmmib'177 \skewchar\sevencmmib'177 \skewchar\fivecmmib'177
 \alloc@@8\fam\chardef\sixt@@n\cmmibfam
 \textfont\cmmibfam\tencmmib
 \scriptfont\cmmibfam\sevencmmib \scriptscriptfont\cmmibfam\fivecmmib
 \font@\tencmbsy=cmbsy10 \font@\sevencmbsy=cmbsy7 \font@\fivecmbsy=cmbsy5
 \skewchar\tencmbsy'60 \skewchar\sevencmbsy'60 \skewchar\fivecmbsy'60
 \alloc@@8\fam\chardef\sixt@@n\cmbsyfam
 \textfont\cmbsyfam\tencmbsy
 \scriptfont\cmbsyfam\sevencmbsy \scriptscriptfont\cmbsyfam\fivecmbsy
 \let\loadbold\empty
}
\def\boldnotloaded#1{\Err@{\ifcase#1\or First\else Second\fi
       bold symbol font not loaded}}
\def\mathchari@#1#2#3{\ifx\undefined\cmmibfam
    \boldnotloaded@\@ne
  \else\mathchar"#1\hexnumber@\cmmibfam#2#3\space \fi}
\def\mathcharii@#1#2#3{\ifx\undefined\cmbsyfam
    \boldnotloaded\tw@
  \else \mathchar"#1\hexnumber@\cmbsyfam#2#3\space\fi}
\edef\bffam@{\hexnumber@\bffam}
\def\boldkey#1{\ifcat\noexpand#1A%
  \ifx\undefined\cmmibfam \boldnotloaded\@ne
  \else {\fam\cmmibfam#1}\fi
 \else
 \ifx#1!\mathchar"5\bffam@21 \else
 \ifx#1(\mathchar"4\bffam@28 \else\ifx#1)\mathchar"5\bffam@29 \else
 \ifx#1+\mathchar"2\bffam@2B \else\ifx#1:\mathchar"3\bffam@3A \else
 \ifx#1;\mathchar"6\bffam@3B \else\ifx#1=\mathchar"3\bffam@3D \else
 \ifx#1?\mathchar"5\bffam@3F \else\ifx#1[\mathchar"4\bffam@5B \else
 \ifx#1]\mathchar"5\bffam@5D \else
 \ifx#1,\mathchari@63B \else
 \ifx#1-\mathcharii@200 \else
 \ifx#1.\mathchari@03A \else
 \ifx#1/\mathchari@03D \else
 \ifx#1<\mathchari@33C \else
 \ifx#1>\mathchari@33E \else
 \ifx#1*\mathcharii@203 \else
 \ifx#1|\mathcharii@06A \else
 \ifx#10\bold0\else\ifx#11\bold1\else\ifx#12\bold2\else\ifx#13\bold3\else
 \ifx#14\bold4\else\ifx#15\bold5\else\ifx#16\bold6\else\ifx#17\bold7\else
 \ifx#18\bold8\else\ifx#19\bold9\else
  \Err@{\string\boldkey\space can't be used with #1}%
 \fi\fi\fi\fi\fi\fi\fi\fi\fi\fi\fi\fi\fi\fi\fi
 \fi\fi\fi\fi\fi\fi\fi\fi\fi\fi\fi\fi\fi\fi}
\def\boldsymbol#1{%
 \DN@{\Err@{You can't use \string\boldsymbol\space with \string#1}#1}%
 \ifcat\noexpand#1A%
   \let\next@\relax
   \ifx\undefined\cmmibfam \boldnotloaded\@ne
   \else {\fam\cmmibfam#1}\fi
 \else
  \xdef\meaning@{\meaning#1.........}%
  \expandafter\math@\meaning@\math@
  \ifmath@
   \expandafter\mathch@\meaning@\mathch@
   \ifmathch@
    \expandafter\boldsymbol@@\meaning@\boldsymbol@@
   \fi
  \else
   \expandafter\macro@\meaning@\macro@
   \expandafter\delim@\meaning@\delim@
   \ifdelim@
    \expandafter\delim@@\meaning@\delim@@
   \else
    \boldsymbol@{#1}%
   \fi
  \fi
 \fi
 \next@}
\def\mathhexboxii@#1#2{\ifx\undefined\cmbsyfam
    \boldnotloaded\tw@
  \else \mathhexbox@{\hexnumber@\cmbsyfam}{#1}{#2}\fi}
\def\boldsymbol@#1{\let\next@\relax\let\next#1%
 \ifx\next\cdot\mathcharii@201 \else
 \ifx\next\prime{{\null\mathcharii@030 \null}}\else
 \ifx\next\lbrack\mathchar"4\bffam@5B \else
 \ifx\next\rbrack\mathchar"5\bffam@5D \else
 \ifx\next\{\mathcharii@466 \else
 \ifx\next\lbrace\mathcharii@466 \else
 \ifx\next\}\mathcharii@567 \else
 \ifx\next\rbrace\mathcharii@567 \else
 \ifx\next\surd{{\mathcharii@170}}\else
 \ifx\next\S{{\mathhexboxii@78}}\else
 \ifx\next\P{{\mathhexboxii@7B}}\else
 \ifx\next\dag{{\mathhexboxii@79}}\else
 \ifx\next\ddag{{\mathhexboxii@7A}}\else
 \DN@{\Err@{You can't use \string\boldsymbol\space with \string#1}#1}%
 \fi\fi\fi\fi\fi\fi\fi\fi\fi\fi\fi\fi\fi}
\def\boldsymbol@@#1.#2\boldsymbol@@{\classnum@#1 \count@@@\classnum@        
 \divide\classnum@4096 \count@\classnum@                                    
 \multiply\count@4096 \advance\count@@@-\count@ \count@@\count@@@           
 \divide\count@@@\@cclvi \count@\count@@                                    
 \multiply\count@@@\@cclvi \advance\count@@-\count@@@                       
 \divide\count@@@\@cclvi                                                    
 \multiply\classnum@4096 \advance\classnum@\count@@                         
 \ifnum\count@@@=\z@                                                        
  \count@"\bffam@ \multiply\count@\@cclvi
  \advance\classnum@\count@
  \DN@{\mathchar\number\classnum@}%
 \else
  \ifnum\count@@@=\@ne                                                      
   \ifx\undefined\cmmibfam \DN@{\boldnotloaded\@ne}%
   \else \count@\cmmibfam \multiply\count@\@cclvi
     \advance\classnum@\count@
     \DN@{\mathchar\number\classnum@}\fi
  \else
   \ifnum\count@@@=\tw@                                                    
     \ifx\undefined\cmbsyfam
       \DN@{\boldnotloaded\tw@}%
     \else
       \count@\cmbsyfam \multiply\count@\@cclvi
       \advance\classnum@\count@
       \DN@{\mathchar\number\classnum@}%
     \fi
  \fi
 \fi
\fi}
\newif\ifdelim@
\newcount\delimcount@
{\uccode`6=`\\ \uccode`7=`d \uccode`8=`e \uccode`9=`l
 \uppercase{\gdef\delim@#1#2#3#4#5\delim@
  {\delim@false\ifx 6#1\ifx 7#2\ifx 8#3\ifx 9#4\delim@true
   \xdef\meaning@{#5}\fi\fi\fi\fi}}}
\def\delim@@#1"#2#3#4#5#6\delim@@{\if#32%
\let\next@\relax
 \ifx\undefined\cmbsyfam \boldnotloaded\@ne
 \else \mathcharii@#2#4#5\space \fi\fi}
\def\vert{\delimiter"026A30C }
\def\Vert{\delimiter"026B30D }
\let\|\Vert

\def\boldkeydots@#1{\bold@true\let\next=#1\let\delayed@=#1\mdots@@
 \boldkey#1\bold@false}  
\def\boldsymboldots@#1{\bold@true\let\next#1\let\delayed@#1\mdots@@
 \boldsymbol#1\bold@false}
\message{Euler fonts,}

\def\frak{\mathfont@\frak}

\def\loadmathfont#1{%
   \expandafter\font@\csname ten#1\endcsname=#110
   \expandafter\font@\csname seven#1\endcsname=#17
   \expandafter\font@\csname five#1\endcsname=#15
   \edef\next{\noexpand\alloc@@8\fam\chardef\sixt@@n
     \expandafter\noexpand\csname#1fam\endcsname}%
   \next
   \textfont\csname#1fam\endcsname \csname ten#1\endcsname
   \scriptfont\csname#1fam\endcsname \csname seven#1\endcsname
   \scriptscriptfont\csname#1fam\endcsname \csname five#1\endcsname
   \expandafter\def\csname #1\expandafter\endcsname\expandafter{%
      \expandafter\mathfont@\csname#1\endcsname}%
 \expandafter\gdef\csname load#1\endcsname{}%
}
\def\mathfont@#1{\RIfM@\expandafter\mathfont@@\expandafter#1\else
  \expandafter\nonmatherr@\expandafter#1\fi}
\def\mathfont@@#1#2{{\mathfont@@@#1{#2}}}
\def\mathfont@@@#1#2{\noaccents@
   \fam\csname\expandafter\eat@\string#1fam\endcsname
   \relax#2}
\message{math accents,}
\def\accentclass@{7}
\def\noaccents@{\def\accentclass@{0}}
\def\makeacc@#1#2{\def#1{\mathaccent"\accentclass@#2 }}
\makeacc@\hat{05E}
\makeacc@\check{014}
\makeacc@\tilde{07E}
\makeacc@\acute{013}
\makeacc@\grave{012}
\makeacc@\dot{05F}
\makeacc@\ddot{07F}
\makeacc@\breve{015}
\makeacc@\bar{016}

\newcount\skewcharcount@
\newcount\familycount@
\def\theskewchar@{\familycount@\@ne
 \global\skewcharcount@\the\skewchar\textfont\@ne                           
 \ifnum\fam>\m@ne\ifnum\fam<16
  \global\familycount@\the\fam\relax
  \global\skewcharcount@\the\skewchar\textfont\the\fam\relax\fi\fi          
 \ifnum\skewcharcount@>\m@ne
  \ifnum\skewcharcount@<128
  \multiply\familycount@256
  \global\advance\skewcharcount@\familycount@
  \global\advance\skewcharcount@28672
  \mathchar\skewcharcount@\else
  \global\skewcharcount@\m@ne\fi\else
 \global\skewcharcount@\m@ne\fi}                                            
\newcount\pointcount@
\def\getpoints@#1.#2\getpoints@{\pointcount@#1 }
\newdimen\accentdimen@
\newcount\accentmu@
\def\dimentomu@{\multiply\accentdimen@ 100
 \expandafter\getpoints@\the\accentdimen@\getpoints@
 \multiply\pointcount@18
 \divide\pointcount@\@m
 \global\accentmu@\pointcount@}
\def\Makeacc@#1#2{\def#1{\RIfM@\DN@{\mathaccent@
 {"\accentclass@#2 }}\else\DN@{\nonmatherr@{#1}}\fi\next@}}
\def\unbracefonts@{\let\Cal@\Cal@@\let\roman@\roman@@\let\bold@\bold@@
 \let\slanted@\slanted@@}
\def\mathaccent@#1#2{\ifnum\fam=\m@ne\xdef\thefam@{1}\else
 \xdef\thefam@{\the\fam}\fi                                                 
 \accentdimen@\z@                                                           
 \setboxz@h{\unbracefonts@$\m@th\fam\thefam@\relax#2$}
 \ifdim\accentdimen@=\z@\DN@{\mathaccent#1{#2}}
  \setbox@ne\hbox{\unbracefonts@$\m@th\fam\thefam@\relax#2\theskewchar@$}
  \setbox\tw@\hbox{$\m@th\ifnum\skewcharcount@=\m@ne\else
   \mathchar\skewcharcount@\fi$}
  \global\accentdimen@\wd@ne\global\advance\accentdimen@-\wdz@
  \global\advance\accentdimen@-\wd\tw@                                     
  \global\multiply\accentdimen@\tw@
  \dimentomu@\global\advance\accentmu@\@ne                                 
 \else\DN@{{\mathaccent#1{#2\mkern\accentmu@ mu}%
    \mkern-\accentmu@ mu}{}}\fi                                             
 \next@}\Makeacc@\Hat{05E}
\Makeacc@\Check{014}
\Makeacc@\Tilde{07E}
\Makeacc@\Acute{013}
\Makeacc@\Grave{012}
\Makeacc@\Dot{05F}
\Makeacc@\Ddot{07F}
\Makeacc@\Breve{015}
\Makeacc@\Bar{016}
\def\Vec{\RIfM@\DN@{\mathaccent@{"017E }}\else
 \DN@{\nonmatherr@\Vec}\fi\next@}
\def\accentedsymbol#1#2{\csname newbox\expandafter\endcsname
  \csname\expandafter\eat@\string#1@box\endcsname
 \expandafter\setbox\csname\expandafter\eat@
  \string#1@box\endcsname\hbox{$\m@th#2$}\define
  #1{\copy\csname\expandafter\eat@\string#1@box\endcsname{}}}
\message{roots,}
\def\sqrt#1{\radical"270370 {#1}}
\let\underline@\underline
\let\overline@\overline
\def\underline#1{\underline@{#1}}
\def\overline#1{\overline@{#1}}
\Invalid@\leftroot
\Invalid@\uproot
\newcount\uproot@
\newcount\leftroot@
\def\root{\relaxnext@
  \DN@{\ifx\next\uproot\let\next@\nextii@\else
   \ifx\next\leftroot\let\next@\nextiii@\else
   \let\next@\plainroot@\fi\fi\next@}%
  \DNii@\uproot##1{\uproot@##1\relax\FN@\nextiv@}%
  \def\nextiv@{\ifx\next\space@\DN@. {\FN@\nextv@}\else
   \DN@.{\FN@\nextv@}\fi\next@.}%
  \def\nextv@{\ifx\next\leftroot\let\next@\nextvi@\else
   \let\next@\plainroot@\fi\next@}%
  \def\nextvi@\leftroot##1{\leftroot@##1\relax\plainroot@}%
   \def\nextiii@\leftroot##1{\leftroot@##1\relax\FN@\nextvii@}%
  \def\nextvii@{\ifx\next\space@
   \DN@. {\FN@\nextviii@}\else
   \DN@.{\FN@\nextviii@}\fi\next@.}%
  \def\nextviii@{\ifx\next\uproot\let\next@\nextix@\else
   \let\next@\plainroot@\fi\next@}%
  \def\nextix@\uproot##1{\uproot@##1\relax\plainroot@}%
  \bgroup\uproot@\z@\leftroot@\z@\FN@\next@}
\def\plainroot@#1\of#2{\setbox\rootbox\hbox{$\m@th\scriptscriptstyle{#1}$}%
 \mathchoice{\r@@t\displaystyle{#2}}{\r@@t\textstyle{#2}}
 {\r@@t\scriptstyle{#2}}{\r@@t\scriptscriptstyle{#2}}\egroup}
\def\r@@t#1#2{\setboxz@h{$\m@th#1\sqrt{#2}$}%
 \dimen@\ht\z@\advance\dimen@-\dp\z@
 \setbox@ne\hbox{$\m@th#1\mskip\uproot@ mu$}\advance\dimen@ 1.667\wd@ne
 \mkern-\leftroot@ mu\mkern5mu\raise.6\dimen@\copy\rootbox
 \mkern-10mu\mkern\leftroot@ mu\boxz@}
\def\boxed#1{\setboxz@h{$\m@th\displaystyle{#1}$}\dimen@.4\ex@
 \advance\dimen@3\ex@\advance\dimen@\dp\z@
 \hbox{\lower\dimen@\hbox{%
 \vbox{\hrule height.4\ex@
 \hbox{\vrule width.4\ex@\hskip3\ex@\vbox{\vskip3\ex@\boxz@\vskip3\ex@}%
 \hskip3\ex@\vrule width.4\ex@}\hrule height.4\ex@}%
 }}}
\message{commutative diagrams,}
\let\ampersand@\relax
\newdimen\minaw@
\minaw@11.11128\ex@
\newdimen\minCDaw@
\minCDaw@2.5pc
\def\minCDarrowwidth#1{\RIfMIfI@\onlydmatherr@\minCDarrowwidth
 \else\minCDaw@#1\relax\fi\else\onlydmatherr@\minCDarrowwidth\fi}
\newif\ifCD@
\def\CD{\bgroup\vspace@\relax\let\ampersand@&\iffalse}\fi
 \CD@true\vcenter\bgroup\Let@\tabskip\z@skip\baselineskip20\ex@
 \lineskip3\ex@\lineskiplimit3\ex@\halign\bgroup
 &\hfill$\m@th##$\hfill\crcr}
\def\endCD{\crcr\egroup\egroup\egroup}
\newdimen\bigaw@
\atdef@>#1>#2>{\ampersand@                                                  
 \setboxz@h{$\m@th\ssize\;{#1}\;\;$}
 \setbox@ne\hbox{$\m@th\ssize\;{#2}\;\;$}
 \setbox\tw@\hbox{$\m@th#2$}
 \ifCD@\global\bigaw@\minCDaw@\else\global\bigaw@\minaw@\fi                 
 \ifdim\wdz@>\bigaw@\global\bigaw@\wdz@\fi
 \ifdim\wd@ne>\bigaw@\global\bigaw@\wd@ne\fi                                
 \ifCD@\enskip\fi                                                           
 \ifdim\wd\tw@>\z@
  \mathrel{\mathop{\hbox to\bigaw@{\rightarrowfill@\displaystyle}}%
    \limits^{#1}_{#2}}
 \else\mathrel{\mathop{\hbox to\bigaw@{\rightarrowfill@\displaystyle}}%
    \limits^{#1}}\fi                                                        
 \ifCD@\enskip\fi                                                          
 \ampersand@}                                                              
\atdef@<#1<#2<{\ampersand@\setboxz@h{$\m@th\ssize\;\;{#1}\;$}%
 \setbox@ne\hbox{$\m@th\ssize\;\;{#2}\;$}\setbox\tw@\hbox{$\m@th#2$}%
 \ifCD@\global\bigaw@\minCDaw@\else\global\bigaw@\minaw@\fi
 \ifdim\wdz@>\bigaw@\global\bigaw@\wdz@\fi
 \ifdim\wd@ne>\bigaw@\global\bigaw@\wd@ne\fi
 \ifCD@\enskip\fi
 \ifdim\wd\tw@>\z@
  \mathrel{\mathop{\hbox to\bigaw@{\leftarrowfill@\displaystyle}}%
       \limits^{#1}_{#2}}\else
  \mathrel{\mathop{\hbox to\bigaw@{\leftarrowfill@\displaystyle}}%
       \limits^{#1}}\fi
 \ifCD@\enskip\fi\ampersand@}
\begingroup
 \catcode`\~=\active \lccode`\~=`\@
 \lowercase{%
  \global\atdef@)#1)#2){~>#1>#2>}
  \global\atdef@(#1(#2({~<#1<#2<}}
\endgroup
\atdef@ A#1A#2A{\llap{$\m@th\vcenter{\hbox
 {$\ssize#1$}}$}\Big\uparrow\rlap{$\m@th\vcenter{\hbox{$\ssize#2$}}$}&&}
\atdef@ V#1V#2V{\llap{$\m@th\vcenter{\hbox
 {$\ssize#1$}}$}\Big\downarrow\rlap{$\m@th\vcenter{\hbox{$\ssize#2$}}$}&&}
\atdef@={&\enskip\mathrel
 {\vbox{\hrule width\minCDaw@\vskip3\ex@\hrule width
 \minCDaw@}}\enskip&}
\atdef@|{\Big\Vert&&}
\atdef@\vert{\Big\Vert&&}
\def\pretend#1\haswidth#2{\setboxz@h{$\m@th\scriptstyle{#2}$}\hbox
 to\wdz@{\hfill$\m@th\scriptstyle{#1}$\hfill}}
\message{poor man's bold,}
\def\pmb{\RIfM@\expandafter\mathpalette\expandafter\pmb@\else
 \expandafter\pmb@@\fi}
\def\pmb@@#1{\leavevmode\setboxz@h{#1}%
   \dimen@-\wdz@
   \kern-.5\ex@\copy\z@
   \kern\dimen@\kern.25\ex@\raise.4\ex@\copy\z@
   \kern\dimen@\kern.25\ex@\box\z@
}
\def\binrel@@#1{\ifdim\wd2<\z@\mathbin{#1}\else\ifdim\wd\tw@>\z@
 \mathrel{#1}\else{#1}\fi\fi}
\newdimen\pmbraise@
\def\pmb@#1#2{\setbox\thr@@\hbox{$\m@th#1{#2}$}%
 \setbox4\hbox{$\m@th#1\mkern.5mu$}\pmbraise@\wd4\relax
 \binrel@{#2}%
 \dimen@-\wd\thr@@
   \binrel@@{%
   \mkern-.8mu\copy\thr@@
   \kern\dimen@\mkern.4mu\raise\pmbraise@\copy\thr@@
   \kern\dimen@\mkern.4mu\box\thr@@
}}
\def\documentstyle#1{\W@{}\input #1.sty\relax}
\message{syntax check,}
\font\dummyft@=dummy
\fontdimen1 \dummyft@=\z@
\fontdimen2 \dummyft@=\z@
\fontdimen3 \dummyft@=\z@
\fontdimen4 \dummyft@=\z@
\fontdimen5 \dummyft@=\z@
\fontdimen6 \dummyft@=\z@
\fontdimen7 \dummyft@=\z@
\fontdimen8 \dummyft@=\z@
\fontdimen9 \dummyft@=\z@
\fontdimen10 \dummyft@=\z@
\fontdimen11 \dummyft@=\z@
\fontdimen12 \dummyft@=\z@
\fontdimen13 \dummyft@=\z@
\fontdimen14 \dummyft@=\z@
\fontdimen15 \dummyft@=\z@
\fontdimen16 \dummyft@=\z@
\fontdimen17 \dummyft@=\z@
\fontdimen18 \dummyft@=\z@
\fontdimen19 \dummyft@=\z@
\fontdimen20 \dummyft@=\z@
\fontdimen21 \dummyft@=\z@
\fontdimen22 \dummyft@=\z@
\def\fontlist@{\\{\tenrm}\\{\sevenrm}\\{\fiverm}\\{\teni}\\{\seveni}%
 \\{\fivei}\\{\tensy}\\{\sevensy}\\{\fivesy}\\{\tenex}\\{\tenbf}\\{\sevenbf}%
 \\{\fivebf}\\{\tensl}\\{\tenit}}
\def\font@#1=#2 {\rightappend@#1\to\fontlist@\font#1=#2 }
\def\dodummy@{{\def\\##1{\global\let##1\dummyft@}\fontlist@}}
\def\nopages@{\output{\setbox\z@\box\@cclv \deadcycles\z@}%
 \alloc@5\toks\toksdef\@cclvi\output}
\let\galleys\nopages@
\newif\ifsyntax@
\newcount\countxviii@
\def\syntax{\syntax@true\dodummy@\countxviii@\count18
 \loop\ifnum\countxviii@>\m@ne\textfont\countxviii@=\dummyft@
 \scriptfont\countxviii@=\dummyft@\scriptscriptfont\countxviii@=\dummyft@
 \advance\countxviii@\m@ne\repeat                                           
 \dummyft@\tracinglostchars\z@\nopages@\frenchspacing\hbadness\@M}
\def\first@#1#2\end{#1}
\def\printoptions{\W@{Do you want S(yntax check),
  G(alleys) or P(ages)?}%
 \message{Type S, G or P, followed by <return>: }%
 \begingroup 
 \endlinechar\m@ne 
 \read\m@ne to\ans@
 \edef\ans@{\uppercase{\def\noexpand\ans@{%
   \expandafter\first@\ans@ P\end}}}%
 \expandafter\endgroup\ans@
 \if\ans@ P
 \else \if\ans@ S\syntax
 \else \if\ans@ G\galleys
 \else\message{? Unknown option: \ans@; using the `pages' option.}%
 \fi\fi\fi}
\def\alloc@#1#2#3#4#5{\global\advance\count1#1by\@ne
 \ch@ck#1#4#2\allocationnumber=\count1#1
 \global#3#5=\allocationnumber
 \ifalloc@\wlog{\string#5=\string#2\the\allocationnumber}\fi}
\def\document{\def\alloclist@{}\def\fontlist@{}}
\let\enddocument\bye

\let\proclaim\undefined
\let\footnote\undefined
\let\=\undefined
\let\>\undefined

\catcode`\@=\active
\message{... finished}

\expandafter\ifx\csname mathdefs.tex\endcsname\relax
  \expandafter\gdef\csname mathdefs.tex\endcsname{}
\else \message{Hey!  Apparently you were trying to
  \string\input{mathdefs.tex} twice.   This does not make sense.} 
\errmessage{Please edit your file (probably \jobname.tex) and remove
any duplicate ``\string\input'' lines}\endinput\fi




\catcode`\X=12\catcode`\@=11

\def\n@wcount{\alloc@0\count\countdef\insc@unt}
\def\n@wwrite{\alloc@7\write\chardef\sixt@@n}
\def\n@wread{\alloc@6\read\chardef\sixt@@n}
\def\r@s@t{\relax}\def\v@idline{\par}\def\@mputate#1/{#1}
\def\l@c@l#1X{\firstpart.#1}\def\gl@b@l#1X{#1}\def\t@d@l#1X{{}}

\def\crossrefs#1{\ifx\all#1\let\tr@ce=\all\else\def\tr@ce{#1,}\fi
   \n@wwrite\cit@tionsout\openout\cit@tionsout=\jobname.cit 
   \write\cit@tionsout{\tr@ce}\expandafter\setfl@gs\tr@ce,}
\def\setfl@gs#1,{\def\@{#1}\ifx\@\empty\let\next=\relax
   \else\let\next=\setfl@gs\expandafter\xdef
   \csname#1tr@cetrue\endcsname{}\fi\next}
\def\m@ketag#1#2{\expandafter\n@wcount\csname#2tagno\endcsname
     \csname#2tagno\endcsname=0\let\tail=\all\xdef\all{\tail#2,}
   \ifx#1\l@c@l\let\tail=\r@s@t\xdef\r@s@t{\csname#2tagno\endcsname=0\tail}\fi
   \expandafter\gdef\csname#2cite\endcsname##1{\expandafter
     \ifx\csname#2tag##1\endcsname\relax?\else\csname#2tag##1\endcsname\fi
     \expandafter\ifx\csname#2tr@cetrue\endcsname\relax\else
     \write\cit@tionsout{#2tag ##1 cited on page \folio.}\fi}
   \expandafter\gdef\csname#2page\endcsname##1{\expandafter
     \ifx\csname#2page##1\endcsname\relax?\else\csname#2page##1\endcsname\fi
     \expandafter\ifx\csname#2tr@cetrue\endcsname\relax\else
     \write\cit@tionsout{#2tag ##1 cited on page \folio.}\fi}
   \expandafter\gdef\csname#2tag\endcsname##1{\expandafter
      \ifx\csname#2check##1\endcsname\relax
      \expandafter\xdef\csname#2check##1\endcsname{}%
      \else\immediate\write16{Warning: #2tag ##1 used more than once.}\fi
      \multit@g{#1}{#2}##1/X%
      \write\t@gsout{#2tag ##1 assigned number \csname#2tag##1\endcsname\space
      on page \number\count0.}%
   \csname#2tag##1\endcsname}}

\def\multit@g#1#2#3/#4X{\def\t@mp{#4}\ifx\t@mp\empty%
      \global\advance\csname#2tagno\endcsname by 1 
      \expandafter\xdef\csname#2tag#3\endcsname
      {#1\number\csname#2tagno\endcsnameX}%
   \else\expandafter\ifx\csname#2last#3\endcsname\relax
      \expandafter\n@wcount\csname#2last#3\endcsname
      \global\advance\csname#2tagno\endcsname by 1 
      \expandafter\xdef\csname#2tag#3\endcsname
      {#1\number\csname#2tagno\endcsnameX}
      \write\t@gsout{#2tag #3 assigned number \csname#2tag#3\endcsname\space
      on page \number\count0.}\fi
   \global\advance\csname#2last#3\endcsname by 1
   \def\t@mp{\expandafter\xdef\csname#2tag#3/}%
   \expandafter\t@mp\@mputate#4\endcsname
   {\csname#2tag#3\endcsname\lastpart{\csname#2last#3\endcsname}}\fi}
\def\t@gs#1{\def\all{}\m@ketag#1e\m@ketag#1s\m@ketag\t@d@l p
\let\realscite\scite
\let\realstag\stag
   \m@ketag\gl@b@l r \n@wread\t@gsin
   \openin\t@gsin=\jobname.tgs \re@der \closein\t@gsin
   \n@wwrite\t@gsout\openout\t@gsout=\jobname.tgs }
\outer\def\localtags{\t@gs\l@c@l}
\outer\def\globaltags{\t@gs\gl@b@l}
\outer\def\newlocaltag#1{\m@ketag\l@c@l{#1}}
\outer\def\newglobaltag#1{\m@ketag\gl@b@l{#1}}

\newif\ifpr@ 
\def\m@kecs #1tag #2 assigned number #3 on page #4.%
   {\expandafter\gdef\csname#1tag#2\endcsname{#3}
   \expandafter\gdef\csname#1page#2\endcsname{#4}
   \ifpr@\expandafter\xdef\csname#1check#2\endcsname{}\fi}
\def\re@der{\ifeof\t@gsin\let\next=\relax\else
   \read\t@gsin to\t@gline\ifx\t@gline\v@idline\else
   \expandafter\m@kecs \t@gline\fi\let \next=\re@der\fi\next}
\def\pretags#1{\pr@true\pret@gs#1,,}
\def\pret@gs#1,{\def\@{#1}\ifx\@\empty\let\n@xtfile=\relax
   \else\let\n@xtfile=\pret@gs \openin\t@gsin=#1.tgs \message{#1} \re@der 
   \closein\t@gsin\fi \n@xtfile}

\newcount\sectno\sectno=0\newcount\subsectno\subsectno=0
\newif\ifultr@local \def\ultralocal{\ultr@localtrue}
\def\firstpart{\number\sectno}
\def\lastpart#1{\ifcase#1 \or a\or b\or c\or d\or e\or f\or g\or h\or 
   i\or k\or l\or m\or n\or o\or p\or q\or r\or s\or t\or u\or v\or w\or 
   x\or y\or z \fi}

\def\resetall{\global\advance\sectno by 1\subsectno=0
   \gdef\firstpart{\number\sectno}\r@s@t}
\def\resetsub{\global\advance\subsectno by 1
   \gdef\firstpart{\number\sectno.\number\subsectno}\r@s@t}
\def\newsection#1\par{\resetall\vskip0pt plus.3\vsize\penalty-250
   \vskip0pt plus-.3\vsize\bigskip\bigskip
   \message{#1}\leftline{\bf#1}\nobreak\bigskip}
\def\subsection#1\par{\ifultr@local\resetsub\fi
   \vskip0pt plus.2\vsize\penalty-250\vskip0pt plus-.2\vsize
   \bigskip\smallskip\message{#1}\leftline{\bf#1}\nobreak\medskip}


\newdimen\marginshift

\newdimen\margindelta
\newdimen\marginmax
\newdimen\marginmin

\def\margininit{       
\marginmax=3 true cm                  
				      
\margindelta=0.1 true cm              
\marginmin=0.1true cm                 
\marginshift=\marginmin
}    

\def\t@gsjj#1,{\def\@{#1}\ifx\@\empty\let\next=\relax\else\let\next=\t@gsjj
   \def\@@{p}\ifx\@\@@\else
   \expandafter\gdef\csname#1cite\endcsname##1{\citejj{##1}}
   \expandafter\gdef\csname#1page\endcsname##1{?}
   \expandafter\gdef\csname#1tag\endcsname##1{\tagjj{##1}}\fi\fi\next}
\newif\ifshowstuffinmargin
\showstuffinmarginfalse
\def\jjtags{\ifx\shlhetal\relax 
  \else
\ifx\shlhetal\undefinedcontrolseq
\else
\showstuffinmargintrue
\ifx\all\relax\else\expandafter\t@gsjj\all,\fi\fi \fi
}

\def\tagjj#1{\realstag{#1}\mginpar{\zeigen{#1}}}
\def\citejj#1{\rechnen{#1}\mginpar{\zeigen{#1}}}     

\def\rechnen#1{\expandafter\ifx\csname stag#1\endcsname\relax ??\else
                           \csname stag#1\endcsname\fi}

\newdimen\theight

\def\marginfont{\sevenrm}

\def\trymarginbox#1{\setbox0=\hbox{\marginfont\hskip\marginshift #1}%
		\global\marginshift\wd0 
		\global\advance\marginshift\margindelta}

\def \mginpar#1{%
\ifvmode\setbox0\hbox to \hsize{\hfill\rlap{\marginfont\quad#1}}%
\ht0 0cm
\dp0 0cm
\box0\vskip-\baselineskip
\else 
             \vadjust{\trymarginbox{#1}%
		\ifdim\marginshift>\marginmax \global\marginshift\marginmin
			\trymarginbox{#1}%
                \fi
             \theight=\ht0
             \advance\theight by \dp0    \advance\theight by \lineskip
             \kern -\theight \vbox to \theight{\rightline{\rlap{\box0}}%
\vss}}\fi}


\def\t@gsoff#1,{\def\@{#1}\ifx\@\empty\let\next=\relax\else\let\next=\t@gsoff
   \def\@@{p}\ifx\@\@@\else
   \expandafter\gdef\csname#1cite\endcsname##1{\zeigen{##1}}
   \expandafter\gdef\csname#1page\endcsname##1{?}
   \expandafter\gdef\csname#1tag\endcsname##1{\zeigen{##1}}\fi\fi\next}
\def\verbatimtags{\showstuffinmarginfalse
\ifx\all\relax\else\expandafter\t@gsoff\all,\fi}
\def\zeigen#1{\hbox{$\langle$}#1\hbox{$\rangle$}}

\def\(#1){\edef\dot@g{\ifmmode\ifinner(\hbox{\noexpand\etag{#1}})
   \else\noexpand\eqno(\hbox{\noexpand\etag{#1}})\fi
   \else(\noexpand\ecite{#1})\fi}\dot@g}

\newif\ifbr@ck
\def\eat#1{}
\def\[#1]{\br@cktrue[\br@cket#1'X]}
\def\br@cket#1'#2X{\def\temp{#2}\ifx\temp\empty\let\next\eat
   \else\let\next\br@cket\fi
   \ifbr@ck\br@ckfalse\br@ck@t#1,X\else\br@cktrue#1\fi\next#2X}
\def\br@ck@t#1,#2X{\def\temp{#2}\ifx\temp\empty\let\neext\eat
   \else\let\neext\br@ck@t\def\temp{,}\fi
   \def\teemp{#1}\ifx\teemp\empty\else\rcite{#1}\fi\temp\neext#2X}
\def\resetbr@cket{\gdef\[##1]{[\rtag{##1}]}}
\def\references{\resetbr@cket\newsection References\par}

\newtoks\symb@ls\newtoks\s@mb@ls\newtoks\p@gelist\n@wcount\ftn@mber
    \ftn@mber=1\newif\ifftn@mbers\ftn@mbersfalse\newif\ifbyp@ge\byp@gefalse
\def\defm@rk{\ifftn@mbers\n@mberm@rk\else\symb@lm@rk\fi}
\def\n@mberm@rk{\xdef\m@rk{{\the\ftn@mber}}%
    \global\advance\ftn@mber by 1 }
\def\rot@te#1{\let\temp=#1\global#1=\expandafter\r@t@te\the\temp,X}
\def\r@t@te#1,#2X{{#2#1}\xdef\m@rk{{#1}}}
\def\b@@st#1{{$^{#1}$}}\def\str@p#1{#1}
\def\symb@lm@rk{\ifbyp@ge\rot@te\p@gelist\ifnum\expandafter\str@p\m@rk=1 
    \s@mb@ls=\symb@ls\fi\write\f@nsout{\number\count0}\fi \rot@te\s@mb@ls}
\def\byp@ge{\byp@getrue\n@wwrite\f@nsin\openin\f@nsin=\jobname.fns 
    \n@wcount\currentp@ge\currentp@ge=0\p@gelist={0}
    \re@dfns\closein\f@nsin\rot@te\p@gelist
    \n@wread\f@nsout\openout\f@nsout=\jobname.fns }
\def\m@kelist#1X#2{{#1,#2}}
\def\re@dfns{\ifeof\f@nsin\let\next=\relax\else\read\f@nsin to \f@nline
    \ifx\f@nline\v@idline\else\let\t@mplist=\p@gelist
    \ifnum\currentp@ge=\f@nline
    \global\p@gelist=\expandafter\m@kelist\the\t@mplistX0
    \else\currentp@ge=\f@nline
    \global\p@gelist=\expandafter\m@kelist\the\t@mplistX1\fi\fi
    \let\next=\re@dfns\fi\next}
\def\symbols#1{\symb@ls={#1}\s@mb@ls=\symb@ls} 
\def\bigsymbol{\textstyle}
\symbols{\bigsymbol\ast,\dagger,\ddagger,\sharp,\flat,\natural,\star}
\def\ftnumbers{\ftn@mberstrue} \def\ftsymbols{\ftn@mbersfalse}
\def\paginal{\byp@ge} \def\resetftnumbers{\ftn@mber=1}
\def\ftnote#1{\defm@rk\expandafter\expandafter\expandafter\footnote
    \expandafter\b@@st\m@rk{#1}}

\long\def\jump#1\endjump{}
\def\ssum{\mathop{\lower .1em\hbox{$\textstyle\Sigma$}}\nolimits}

\def\qed{\nobreak\kern 1em \vrule height .5em width .5em depth 0em}
\def\newneq{\hbox{\rlap{\hbox to 1\wd9{\hss$=$\hss}}\raise .1em 
   \hbox to 1\wd9{\hss$\scriptscriptstyle/$\hss}}}
\def\subsetne{\setbox9 = \hbox{$\subset$}\mathrel{\hbox{\rlap
   {\lower .4em \newneq}\raise .13em \hbox{$\subset$}}}}
\def\supsetne{\setbox9 = \hbox{$\subset$}\mathrel{\hbox{\rlap
   {\lower .4em \newneq}\raise .13em \hbox{$\supset$}}}}

\def\vbar{\mathchoice{\vrule height6.3ptdepth-.5ptwidth.8pt\kern-.8pt}
   {\vrule height6.3ptdepth-.5ptwidth.8pt\kern-.8pt}
   {\vrule height4.1ptdepth-.35ptwidth.6pt\kern-.6pt}
   {\vrule height3.1ptdepth-.25ptwidth.5pt\kern-.5pt}}
\def\f@dge{\mathchoice{}{}{\mkern.5mu}{\mkern.8mu}}
\def\b@c#1#2{{\rm \mkern#2mu\vbar\mkern-#2mu#1}}
\def\b@b#1{{\rm I\mkern-3.5mu #1}}
\def\b@a#1#2{{\rm #1\mkern-#2mu\f@dge #1}}
\def\bb#1{{\count4=`#1 \advance\count4by-64 \ifcase\count4\or\b@a A{11.5}\or
   \b@b B\or\b@c C{5}\or\b@b D\or\b@b E\or\b@b F \or\b@c G{5}\or\b@b H\or
   \b@b I\or\b@c J{3}\or\b@b K\or\b@b L \or\b@b M\or\b@b N\or\b@c O{5} \or
   \b@b P\or\b@c Q{5}\or\b@b R\or\b@a S{8}\or\b@a T{10.5}\or\b@c U{5}\or
   \b@a V{12}\or\b@a W{16.5}\or\b@a X{11}\or\b@a Y{11.7}\or\b@a Z{7.5}\fi}}

\catcode`\X=11 \catcode`\@=12


\expandafter\ifx\csname citeadd.tex\endcsname\relax
\expandafter\gdef\csname citeadd.tex\endcsname{}
\else \message{Hey!  Apparently you were trying to
\string\input{citeadd.tex} twice.   This does not make sense.} 
\errmessage{Please edit your file (probably \jobname.tex) and remove
any duplicate ``\string\input'' lines}\endinput\fi

\sectno=-1   
\localtags
\jjtags
\NoBlackBoxes
\define\mr{\medskip\roster}
\define\sn{\smallskip\noindent}

\define\bn{\bigskip\noindent}
\define\ub{\underbar}

\define\ermn{\endroster\medskip\noindent}

\define \nl{\newline}
\magnification=\magstep 1
\documentstyle{amsppt}

{    
\catcode`@11

\ifx\alicetwothousandloaded@\relax
  \endinput\else\global\let\alicetwothousandloaded@\relax\fi

\gdef\subjclass{\let\savedef@\subjclass
 \def\subjclass##1\endsubjclass{\let\subjclass\savedef@
   \toks@{\def\usualspace{{\rm\enspace}}\eightpoint}%
   \toks@@{##1\unskip.}%
   \edef\thesubjclass@{\the\toks@
     \frills@{{\noexpand\rm2000 {\noexpand\it Mathematics Subject
       Classification}.\noexpand\enspace}}%
     \the\toks@@}}%
  \nofrillscheck\subjclass}
} 


\expandafter\ifx\csname alice2jlem.tex\endcsname\relax
  \expandafter\xdef\csname alice2jlem.tex\endcsname{\the\catcode`@}
\else \message{Hey!  Apparently you were trying to
\string\input{alice2jlem.tex}  twice.   This does not make sense.}
\errmessage{Please edit your file (probably \jobname.tex) and remove
any duplicate ``\string\input'' lines}\endinput\fi

\expandafter\ifx\csname bib4plain.tex\endcsname\relax
  \expandafter\gdef\csname bib4plain.tex\endcsname{}
\else \message{Hey!  Apparently you were trying to \string\input
  bib4plain.tex twice.   This does not make sense.}
\errmessage{Please edit your file (probably \jobname.tex) and remove
any duplicate ``\string\input'' lines}\endinput\fi

\def\renewcommand{\newcommand}	       
\edef\cite{\the\catcode`@}%
\catcode`@ = 11
\let\@oldatcatcode = \cite
\chardef\@letter = 11
\chardef\@other = 12
%
%
%
%
\def\@innerdef#1#2{\edef#1{\expandafter\noexpand\csname #2\endcsname}}%
%
%
\@innerdef\@innernewcount{newcount}%
\@innerdef\@innernewdimen{newdimen}%
\@innerdef\@innernewif{newif}%
\@innerdef\@innernewwrite{newwrite}%
%
%
%
\def\@gobble#1{}%
%
%
%
\ifx\inputlineno\@undefined
   \let\@linenumber = \empty 
\else
   \def\@linenumber{\the\inputlineno:\space}%
\fi
%
%
%
\def\@futurenonspacelet#1{\def\cs{#1}%
   \afterassignment\@stepone\let\@nexttoken=
}%
\begingroup 
\def\\{\global\let\@stoken= }%
\\ 
\endgroup
\def\@stepone{\expandafter\futurelet\cs\@steptwo}%
\def\@steptwo{\expandafter\ifx\cs\@stoken\let\@@next=\@stepthree
   \else\let\@@next=\@nexttoken\fi \@@next}%
\def\@stepthree{\afterassignment\@stepone\let\@@next= }%
%
%
%
\def\@getoptionalarg#1{%
   \let\@optionaltemp = #1%
   \let\@optionalnext = \relax
   \@futurenonspacelet\@optionalnext\@bracketcheck
}%
%
%
\def\@bracketcheck{%
   \ifx [\@optionalnext
      \expandafter\@@getoptionalarg
   \else
      \let\@optionalarg = \empty
      \expandafter\@optionaltemp
   \fi
}%
\def\@@getoptionalarg[#1]{%
   \def\@optionalarg{#1}%
   \@optionaltemp
}%
%
%
%
\def\@nnil{\@nil}%
\def\@fornoop#1\@@#2#3{}%
\def\@for#1:=#2\do#3{%
   \edef\@fortmp{#2}%
   \ifx\@fortmp\empty \else
      \expandafter\@forloop#2,\@nil,\@nil\@@#1{#3}%
   \fi
}%
\def\@forloop#1,#2,#3\@@#4#5{\def#4{#1}\ifx #4\@nnil \else
       #5\def#4{#2}\ifx #4\@nnil \else#5\@iforloop #3\@@#4{#5}\fi\fi
}%
\def\@iforloop#1,#2\@@#3#4{\def#3{#1}\ifx #3\@nnil
       \let\@nextwhile=\@fornoop \else
      #4\relax\let\@nextwhile=\@iforloop\fi\@nextwhile#2\@@#3{#4}%
}%
%
%
%
\@innernewif\if@fileexists
\def\@testfileexistence{\@getoptionalarg\@finishtestfileexistence}%
\def\@finishtestfileexistence#1{%
   \begingroup
      \def\extension{#1}%
      \immediate\openin0 =
         \ifx\@optionalarg\empty\jobname\else\@optionalarg\fi
         \ifx\extension\empty \else .#1\fi
         \space
      \ifeof 0
         \global\@fileexistsfalse
      \else
         \global\@fileexiststrue
      \fi
      \immediate\closein0
   \endgroup
}%
%
%
%
%
\def\bibliographystyle#1{%
   \@readauxfile
   \@writeaux{\string\bibstyle{#1}}%
}%
\let\bibstyle = \@gobble
%
%
\let\bblfilebasename = \jobname
\def\bibliography#1{%
   \@readauxfile
   \@writeaux{\string\bibdata{#1}}%
   \@testfileexistence[\bblfilebasename]{bbl}%
   \if@fileexists
      \nobreak
      \@readbblfile
   \fi
}%
\let\bibdata = \@gobble
%
%
\def\nocite#1{%
   \@readauxfile
   \@writeaux{\string\citation{#1}}%
}%
\@innernewif\if@notfirstcitation
%
%
\def\cite{\@getoptionalarg\@cite}%
%
%
\def\@cite#1{%
   \let\@citenotetext = \@optionalarg
   \printcitestart
   \nocite{#1}%
   \@notfirstcitationfalse
   \@for \@citation :=#1\do
   {%
      \expandafter\@onecitation\@citation\@@
   }%
   \ifx\empty\@citenotetext\else
      \printcitenote{\@citenotetext}%
   \fi
   \printcitefinish
}%
\def\@onecitation#1\@@{%
   \if@notfirstcitation
      \printbetweencitations
   \fi
   \expandafter \ifx \csname\@citelabel{#1}\endcsname \relax
      \if@citewarning
         \message{\@linenumber Undefined citation `#1'.}%
      \fi
      \expandafter\gdef\csname\@citelabel{#1}\endcsname{%
\strut
\vadjust{\vskip-\dp\strutbox
\vbox to 0pt{\vss\parindent0cm \leftskip=\hsize 
\advance\leftskip3mm
\advance\hsize 4cm\strut\openup-4pt 
\rightskip 0cm plus 1cm minus 0.5cm ?  #1 ?\strut}}
         {\tt
            \escapechar = -1
            \nobreak\hskip0pt
            \expandafter\string\csname#1\endcsname
            \nobreak\hskip0pt
         }%
      }%
   \fi
   \csname\@citelabel{#1}\endcsname
   \@notfirstcitationtrue
}%
%
%
\def\@citelabel#1{b@#1}%
%
%
\def\@citedef#1#2{\expandafter\gdef\csname\@citelabel{#1}\endcsname{#2}}%
%
%
%
\def\@readbblfile{%
   \ifx\@itemnum\@undefined
      \@innernewcount\@itemnum
   \fi
   \begingroup
      \def\begin##1##2{%
         \setbox0 = \hbox{\biblabelcontents{##2}}%
         \biblabelwidth = \wd0
      }%
      \def\end##1{}
      %
      %
      \@itemnum = 0
      \def\bibitem{\@getoptionalarg\@bibitem}%
      \def\@bibitem{%
         \ifx\@optionalarg\empty
            \expandafter\@numberedbibitem
         \else
            \expandafter\@alphabibitem
         \fi
      }%
      \def\@alphabibitem##1{%
         \expandafter \xdef\csname\@citelabel{##1}\endcsname {\@optionalarg}%
         \ifx\biblabelprecontents\@undefined
            \let\biblabelprecontents = \relax
         \fi
         \ifx\biblabelpostcontents\@undefined
            \let\biblabelpostcontents = \hss
         \fi
         \@finishbibitem{##1}%
      }%
      \def\@numberedbibitem##1{%
         \advance\@itemnum by 1
         \expandafter \xdef\csname\@citelabel{##1}\endcsname{\number\@itemnum}%
         \ifx\biblabelprecontents\@undefined
            \let\biblabelprecontents = \hss
         \fi
         \ifx\biblabelpostcontents\@undefined
            \let\biblabelpostcontents = \relax
         \fi
         \@finishbibitem{##1}%
      }%
      \def\@finishbibitem##1{%
         \biblabelprint{\csname\@citelabel{##1}\endcsname}%
         \@writeaux{\string\@citedef{##1}{\csname\@citelabel{##1}\endcsname}}%
         \ignorespaces
      }%
      %
      %
      \let\em = \bblem
      \let\newblock = \bblnewblock
      \let\sc = \bblsc
      \frenchspacing
      \clubpenalty = 4000 \widowpenalty = 4000
      \tolerance = 10000 \hfuzz = .5pt
      \everypar = {\hangindent = \biblabelwidth
                      \advance\hangindent by \biblabelextraspace}%
      \bblrm
      \parskip = 1.5ex plus .5ex minus .5ex
      \biblabelextraspace = .5em
      \bblhook
      \input \bblfilebasename.bbl
   \endgroup
}%
%
%
\@innernewdimen\biblabelwidth
\@innernewdimen\biblabelextraspace
%
%
%
\def\biblabelprint#1{%
   \noindent
   \hbox to \biblabelwidth{%
      \biblabelprecontents
      \biblabelcontents{#1}%
      \biblabelpostcontents
   }%
   \kern\biblabelextraspace
}%
%
%
%
\def\biblabelcontents#1{{\bblrm [#1]}}%
%
%
\def\bblrm{\rm}%
%
%
\def\bblem{\it}%
%
%
\def\bblsc{\ifx\@scfont\@undefined
              \font\@scfont = cmcsc10
           \fi
           \@scfont
}%
%
%
\def\bblnewblock{\hskip .11em plus .33em minus .07em }%
%
%
\let\bblhook = \empty
%
%
%
\def\printcitestart{[}
\def\printcitefinish{]}
\def\printbetweencitations{, }
\def\printcitenote#1{, #1}
%
%
%
\let\citation = \@gobble
%
%
%
\@innernewcount\@numparams
%
%
\def\newcommand#1{%
   \def\@commandname{#1}%
   \@getoptionalarg\@continuenewcommand
}%
%
%
\def\@continuenewcommand{%
   \@numparams = \ifx\@optionalarg\empty 0\else\@optionalarg \fi \relax
   \@newcommand
}%
%
%
\def\@newcommand#1{%
   \def\@startdef{\expandafter\edef\@commandname}%
   \ifnum\@numparams=0
      \let\@paramdef = \empty
   \else
      \ifnum\@numparams>9
         \errmessage{\the\@numparams\space is too many parameters}%
      \else
         \ifnum\@numparams<0
            \errmessage{\the\@numparams\space is too few parameters}%
         \else
            \edef\@paramdef{%
               \ifcase\@numparams
                  \empty  No arguments.
               \or ####1%
               \or ####1####2%
               \or ####1####2####3%
               \or ####1####2####3####4%
               \or ####1####2####3####4####5%
               \or ####1####2####3####4####5####6%
               \or ####1####2####3####4####5####6####7%
               \or ####1####2####3####4####5####6####7####8%
               \or ####1####2####3####4####5####6####7####8####9%
               \fi
            }%
         \fi
      \fi
   \fi
   \expandafter\@startdef\@paramdef{#1}%
}%
%
%
%
%
\def\@readauxfile{%
   \if@auxfiledone \else 
      \global\@auxfiledonetrue
      \@testfileexistence{aux}%
      \if@fileexists
         \begingroup
            \endlinechar = -1
            \catcode`@ = 11
            \input \jobname.aux
         \endgroup
      \else
         \message{\@undefinedmessage}%
         \global\@citewarningfalse
      \fi
      \immediate\openout\@auxfile = \jobname.aux
   \fi
}%
%
%
\newif\if@auxfiledone
\ifx\noauxfile\@undefined \else \@auxfiledonetrue\fi
%
%
%
%
\@innernewwrite\@auxfile
\def\@writeaux#1{\ifx\noauxfile\@undefined \write\@auxfile{#1}\fi}%
%
%
%
\ifx\@undefinedmessage\@undefined
   \def\@undefinedmessage{No .aux file; I won't give you warnings about
                          undefined citations.}%
\fi
%
%
\@innernewif\if@citewarning
\ifx\noauxfile\@undefined \@citewarningtrue\fi
%
%
%
\catcode`@ = \@oldatcatcode


\def\widestnumber#1#2{}

\def\citewarning#1{\ifx\shlhetal\relax 
    \else
    \par{#1}\par
    \fi
}

\def\rm{\fam0 \tenrm}

\def\fakesubhead#1\endsubhead{\bigskip\noindent{\bf#1}\par}



%
%
%

%

\font\textrsfs=rsfs10
\font\scriptrsfs=rsfs7
\font\scriptscriptrsfs=rsfs5

\newfam\rsfsfam
\textfont\rsfsfam=\textrsfs
\scriptfont\rsfsfam=\scriptrsfs
\scriptscriptfont\rsfsfam=\scriptscriptrsfs

\edef\oldcatcodeofat{\the\catcode`\@}
\catcode`\@11

\def\Cal@@#1{\noaccents@ \fam \rsfsfam #1}

\catcode`\@\oldcatcodeofat


\expandafter\ifx \csname margininit\endcsname \relax\else\margininit\fi

\pageheight{8.5truein}
\topmatter
\title{Logical Dreams} \endtitle
\author {Saharon Shelah \thanks {\null\newline I would like to thank 
Alice Leonhardt for the beautiful typing. \null\newline
This paper is based on my lecture (and the preparations to the lecture)
during the conference {\bf Mathematical Challenges of the 21$^{th}$
Century}. Publication E23 } \endthanks} \endauthor

\affil{Institute of Mathematics\\
 The Hebrew University\\
 Jerusalem, Israel
 \medskip
 Rutgers University\\
 Mathematics Department\\
 New Brunswick, NJ  USA} \endaffil

\subjclass{{{03-02, 03Exx}}} \endsubjclass

\abstract  We discuss the past and future of set theory,
axiom systems and independence results.  We deal in particular with
cardinal arithmetic. \endabstract

\toc
\subhead ${{}}$ {\bf Reading instructions} \endsubhead
\subhead 0. {\bf Introductory remarks} \endsubhead
\subhead 1  {\bf What does mathematical logic do for you?} \endsubhead
\subhead 2  {\bf The glory of proven ignorance} \endsubhead
\subhead 3  {\bf Set theory and additional axioms} \endsubhead
\subhead 4  {\bf ZFC indecisiveness} \endsubhead
\subhead 5  {\bf Is ZFC really so weak?} \endsubhead
\subhead 6  {\bf Concluding remarks} \endsubhead
\endtoc
\endtopmatter
\document  
 
\newpage 

\head {Annotated Content} \endhead  \resetall 
\bn
\ \ \ \ $\quad$ Reading instructions.
\bn
\S0 $\quad$ Introductory remarks.
\mr
\item "{${{}}$}"  [How to read; what this article tries to do; what is
cardinal arithmetic.]
\ermn
\S1 $\quad$ What does mathematical logic do for me?
\mr
\item "{${{}}$}"  [Why is it needed; what is this animal ZFC; thesis:
general results; the Jeffersonian thesis for the 
best framework; the scale thesis.]
\ermn
\S2 $\quad$ The glory of proven ignorance.
\mr
\item "{${{}}$}" [G\"odel's diet, Cohen's fattening; independence in number
theory; primes are random; the Riemann hypothesis, consistency strength.]
\ermn
\S3 $\quad$ Set theory and additional axioms.
\mr
\item "{${{}}$}"  [G\"odel sentences and large cardinals; semi-axioms;
new axioms; semi-axioms on $AD_{{\bold L}[{\Bbb R}]}$ positive but not true;  
is there just one nice descriptive set theory? cardinal invariants for
the continuum.]
\ermn
\S4 $\quad$ ZFC indecisiveness.
\mr
\item "{${{}}$}"  [Large cardinals, equi-consistency;
supercompacts; significance of values of $2^{\aleph_0}$; 
dividing lines/watershed lines; independence over $\bold L$.]
\ermn
\S5 $\quad$ Is ZFC really so weak?
\mr
\item "{${{}}$}"  [A new perspective on cardinal arithmetic; rubble 
removal thesis; treasures thesis; pcf; relatives of GCH, dividing
lines/watershed lines.]
\ermn
\S6 $\quad$ Concluding Remarks.
\newpage

\noindent {\bf Reading instructions}: The real instructions, i.e., my hopes,
are that you start in the beginning and read untill the end. But if you look
for advice on how much not to read and still get what you like, note that
the sections are essentially independent. If you wonder what is the
relevance of mathematical logic to the rest of mathematics and whether you
implicitly accept the usual axioms of set theory (that is, the axioms of
E.~{\bf Z}ermelo and A.~{\bf F}raenkel with the Axiom of {\bf C}hoice,
further referred to as ZFC), you should read mainly \S1.

If you are excited about the independence phenomenon you should read \S2.
If you are interested in considerations on new additional axioms, and/or in
looking at definable sets of reals (as in descriptive set theory) or general
sets of reals, and in looking more deeply into independence, you should go to
\S3. For more on independence in set theory see \S4.  If you really are
interested in cardinal arithmetic, i.e., the arithmetic of infinite numbers
of G.~Cantor, go to \S5.

I had intended to also write a section on ``Pure model theory", ``Applied
model theory", and ``On Speculations and Nonsense" but only \S6
materialized. 

I thank J. Baldwin, A. Blass, G. Cherlin, P. Eklof, U. Hrushovski,
J. Kennedy, C. Laskowski, A. Levy, A. Ros{\l}anowski, J. V\"a\"an\"anen and
A. Villaveces for many comments.

\bn
\head {\S0 Introductory remarks} \endhead  \resetall \sectno=0
\medskip

The intention to set the agenda for the $21^{\text{th}}$ century for
mathematical logic is certainly over-ambitious, if not to say megalomaniac.

Unavoidably, I will speak mostly on directions I am interested in and/or
relatively knowledgeable in (which are quite close to those I have worked
on), so this selection will be riddled with prejudices but at least they are
mine; hopefully some others will be infuriated enough to offer differing
opinions.
\medskip

I will try to make this article accessible to any mathematician (if (s)he
ignores some more specialized parts), as this is intended for a general
audience (of mathematicians); this differs from my recent papers on open
problems, \cite{Sh:666}, \cite{Sh:702}.  Also, as readers tend to get lost,
each section makes a fresh start and even subsections can be read separately
(so from a negative point of view, there are repetitions and a somewhat
fragmented character).

Usually we try to appeal to non-expert readers, do our best, fail, and
barely make it interesting for the professional, but I am doomed to make the
attempt. {\it Many of the issues discussed below have been very popular and
everyone has an opinion; I will mention some other opinions in order to
disagree with them}.
\medskip

\ub{Cardinal Arithmetic}:  We shall deal several times with \ub{cardinal
arithmetic}, so recall that two (possibly infinite) sets
$A,B$ are called equinumerous (in short: $A\sim B$), if there is a
one-to-one mapping from one onto the other.  Essentially, a cardinal number
$\lambda$ is $A/{\sim}$ written as $|A|$.  We define addition,
multiplication and exponentiation of those numbers by 
\mr
\item "{}" $|A| + |B| = |A\cup B|$, when $A,B$ are disjoint, 
\item "{}" $|A| \times |B| = |A \times B|$,
\item "{}" $|A|^{|B|} = |{}^B A| = |\{f:f$ a function from $B$ into $A\}|$,
also 
\item "{}" $\prod\limits_{i\in I} |A_i| = |\prod\limits_{i \in I} A_i|$.
\endroster
Then all the usual equalities and weak $(\le)$ inequalities hold, but for
infinite cardinals $\lambda,\mu$ we have $\lambda +1=\lambda$, $\lambda+\mu=
\max\{\lambda,\mu\}=\lambda\times\mu$, but $2^\lambda>\lambda$ for every $\lambda$.

Let $\aleph_0$ be the number of natural numbers, $\aleph_{n+1}$ the
successor of $\aleph_n$, $\aleph_\omega = \sum\{\aleph_n:n \in \Bbb 
N\}$.  

The {\bf C}ontinuum {\bf H}ypothesis (CH) is the statement $2^{\aleph_0} =
\aleph_1$, and the {\bf G}eneralized {\bf C}ontinuum {\bf H}ypothesis (GCH)
says that $2^\lambda$ is $\lambda^+$ (the successor of $\lambda$) for every
infinite $\lambda$. In fact, GCH says that not only the rules of addition
and multiplication of infinite numbers are simple, but also the rules of
exponentiation.

An ordinal is the order type of a linear order which is well ordered,
i.e., whose every non empty subset has a first element.  In fact we can
choose a representative; for ordinals $\alpha,\beta$ the meaning of $\alpha <
\beta$ is clear and, in fact, $(\{\alpha:\alpha < \beta\},<)$ is of
order type $\alpha$.  A cardinal number $\lambda$ is identified with
$\{\alpha:\alpha$ is an ordinal and $|\{\beta:\beta < \alpha\}| <
\lambda\}$, a set of cardinality $\lambda$; this is the representative.

We denote by $\aleph_\alpha$ the $\alpha$-th infinite cardinal.
\bigskip\bigskip

\noindent
\head{\S1 What does mathematical logic do for you?} 
\endhead  \resetall \sectno=1
\medskip

What does {\it mathematical logic} deal with?

It is {\bf mathematics applied to mathematics itself} (and to some problems
in philosophy).

But what does it do for you, a mathematician from another field?  It does
not help you to solve a thorny differential equation or anything like
that. But if you suspect that the Riemann Hypothesis or P=NP is {\it
undecidable}, or, say, cannot be decided by your present methods by whatever
ingenuity which can be mastered, then you are on the hook: how can you even
phrase this coherently, let alone {\it prove} anything like this?

For such questions you need to phrase a general framework for doing
mathematics (this is set theory -- usually axiomatized in the ``ZFC set of
axioms\footnote{axioms of Zermelo and Fraenkel with Axiom of Choice}'' or
relatives of it).

You also need to {\it define} what a mathematical proof is (G\"odel's
completeness theorem as well as G\"odel's incompleteness theorem speak about
the relationship between provability and truth; more on this than I would
like to know is the subject of proof theory).

You may develop mathematical theories well enough without looking for
mathematical logic, but if you like a formalization of ``a mathematical
theory'' and consider investigation at such a level of generality, you
arrive at model theory. You may wonder whether anything interesting can be
said about such an arbitrary theory; does model theory just grind water or
does it have theorems with meat?  Does model theory pass this test?
Certainly, I definitely think that it does.  Does model theory have
relevance to other parts of mathematics?  I think that for a general theory
to give interesting results when specialized to older contexts is strong
evidence of its being deep (though certainly not a necessary condition).
Naturally, I think that model theory passes this test too, and notions which
arise from investigating structure of a general class of structures, in
arbitrary cardinality later serve in investigating classical objects.
Though I have much to say on model theory, I have not managed to write what
I would like to say on it, so except for \S6 the present paper does not deal
with model theory, but only with set theory.

Similarly if you like to know about inherent limitations of algorithms and the 
semi-lattice of Turing degrees, you turn to
computability (of old called recursion theory).
\medskip

\ub{What is this animal ``ZFC''?}
\smallskip

Fundamental for us is ZFC. I feel it is not something esoteric for
you. Rather, most mathematicians ``do not know'' of ZFC just as Monsieur
Jourdain in Moli\`ere's {\it Bourgeois Gentilhomme} ``did not know'' he was
speaking in prose; or, in the same way that most of us use freely the law of
the excluded middle, ``$A$ or non-$A$'', without any scruples, because we do
not even notice that we are using it.  Many times, we can read but not
understand (mathematical writings).  However, I guarantee you will
understand, though not necessarily be able to read (including most of my
audience in LA and unfortunately myself) the following paragraph:

\block
{\it 
 \noindent  MONSIEUR JOURDAIN: Quoi? quand je dis: ``Nicole, apportez-moi mes
 pantoufles, et me donnez mon bonnet de nuit" , c'est de la prose? 

 \noindent  MAITRE DE PHILOSOPHIE: Oui, Monsieur.

 \noindent  MONSIEUR JOURDAIN: Par ma foi! il y a plus de quarante ans que
je dis de la prose sans que j'en susse rien, et je vous suis le plus
oblig\'e du monde de m'avoir appris cela.}
\endblock
When I explain the axioms of ZFC, the usual response is ``fine, but by what
right do you assume CH?" (the Continuum Hypothesis, see below, which is not
included).
\medskip

In non-technical terms, 
\medskip

\ub{{\stag{1.0} The Jourdain Thesis}}:\quad  If you deal freely with the set
of reals, or the vector space of functions from the reals to the reals or,
e.g., with Hom$(G,H)$ for groups $G,H$, then you -- for all practical
purposes -- work within and accept the axioms of ZFC, that is they formalize
whatever proofs you accept. Even other positions (like omitting the axiom of
choice or replacing it by weak versions) are best dealt within this
framework.    
\medskip

\ub{{\stag{1.1} The Generality Thesis}}: \quad If you would like to have
{\it general} results, you have to use a set theoretic framework.  While by
now we know well how to generate ``generalized nonsense" which grinds water
and tells us nothing new, many times a general framework shows you that
isolated claims are parts of general phenomena; the real test is whether you
discover deep beautiful mathematics.
\medskip

That is, if you want to know something about {\bf ALL} structures of some
kind (all groups, all manifolds, etc), then you need to be able to deal with
infinite unions or infinite products of sets, which are inherently
set-theoretical concepts.

Moreover, even if your main interest is in, say, finitely generated groups,
you will be drawn into more general ones, e.g., taking some
compactification, or using infinite products. 

How far should you go? When to stop? I suggest that we adopt the following
position: 
\medskip

\ub{{\stag{1.2} The Jeffersonian Thesis}}:\quad The best framework, the best
foundation, is the one that governs you  least; that is you do not notice
its restrictions (except when really necessary, like arriving at a
contradiction). 
\medskip

For instance, look at the following scheme of a proof: given a colouring 
of the natural numbers with red and green, we look at the Stone-\v Cech
compactification $\beta({\Bbb N})$ of ${\Bbb N}$, which is just the set of 
ultrafilters on ${\Bbb N}$ (the principal ultrafilters are identified
with the natural numbers); we define on it an operation $+$:
$$p+q=\{A\subseteq{\Bbb N}:\{n:n+A \in p\}\in q\},$$ 
and using an easy fixed point theorem in the semi-group $(\beta({\Bbb
N}),+)$ we get an infinite monochromatic set such that its members and all
(finite non-trivial) sums of members (with no repetition) have the same
color (this is a well known, very elegant alternative proof of  Hindman's
theorem). Would you object to such a proof? Or, would you stop at the power
set ${\Cal P}({\Bbb N})$ of $\Bbb N$?  First it seems to me unnatural not to
have ${\Cal P}({\Cal P}({\Bbb N}))$ and second you will not gain much, e.g.,
the problems of cardinal arithmetic are replaced by relatives even if you
consider only cardinalities of sets of reals (on pcf theory see \S5; also
definable sets may give such phenomena). 

There is a natural scale of theories, some stronger than ZFC (large
cardinals), some weaker (e.g., PA = Peano arithmetic = a set of axioms for
the natural numbers in first order logic, describing addition,
multiplication, and the scheme of mathematical induction).

PA already tells us that the universe is infinite, but PA ``stops'' after we
have all the natural numbers. ZFC goes beyond the natural numbers; in ZFC we
can distinguish different infinite cardinalities such as ``countable'' and
``uncountable'', and we can show that there are infinitely many
cardinalities, uncountably many, etc. 

But there are also set theories stronger than ZFC, which are as high above
ZFC as ZFC is above PA, and even higher.
\medskip

\ub{{\stag{1.3} The Scale Thesis}}:\quad Even if you feel ZFC assumes too
much or too little (and you do not work artificially), you will end up
somewhere along this scale, going from PA to the large cardinals.  \nl 
(What does ``artificial'' mean?  For Example, there are 17 strongly
inaccessible \footnote{that is, an uncountable cardinal $\lambda$ satisfying
$\mu<\lambda\Rightarrow 2^\mu<\lambda$ and: if $\lambda_t < \lambda$ for
$t\in I$ where $|I|<\lambda$, then $\sum\limits_{t\in I}\lambda_t<\lambda$;
in other words is strong limit and regular} cardinals, the theory ZFC +
``there are 84 strongly inaccessible cardinals" is contradictory and the
theory ZFC + ``there are 49 strongly inaccessible cardinals" is consistent
but has no well found model.) 
\medskip

An extreme skeptic goes ``below PA", e.g., (s)he may doubt not only whether 
$2^n$  (for every natural number $n$) necessarily exists but even whether 
$n^{[\log n]}$ exists. [In the latter case (s)he still has a chance to prove
``there are infinitely many primes".]  The difference between two such
positions will be just where they put their belief; so the theory is quite 
translatable, just a matter of stress.  For instance, by one we know that
there are infinitely many primes, by the other we have an implication.
There is a body of work supporting this, the so called equi-consistency
results (e.g., on real valued measurable cardinals, see later).

So far I have mainly defended accepting ZFC, as for believing in more,
see later.
\bn

\head {\S2 The glory of proven ignorance: To {\it show} that we cannot
know!} \endhead  \resetall 
\medskip

In short: The Continuum Problem asks: 

How many real numbers are there?
\medskip

G.~Cantor proved: There are  {\it more} reals than rationals. (In a
technical sense: ``uncountable'', ``there is no bijection from $\Bbb R$ into
$\Bbb Q$''). 

The Continuum Hypothesis (CH) says: yes, more, but barely so. Every set $A
\subseteq {\Bbb R}$ is either countable or equinumerous with $\Bbb R$.
\medskip

K.~G\"odel proved:  Perhaps CH holds.

P.~Cohen proved:   Perhaps CH does not hold.
\medskip

What is a better starting point, than Hilbert's first problem, {\it the
Continuum Problem}?  It deals with the arithmetic of infinite numbers,
called cardinals, discovered by Cantor; they are just the equivalence
classes of the relation ``there is a one-to-one mapping from $A$ onto $B$"
with natural operations.  Specifically, it asks whether the continuum, the
number of reals, is the successor of the number of natural numbers. In
other words, can we for any infinite set of real numbers find a one-to-one
map from it onto the set of all reals or onto the set of natural
numbers. {\it I find this a great problem, in particular it has induced
much of the 20th century achievements in set theory.}

G\"odel has shown that the Generalized Continuum Hypothesis (GCH) may hold.
In fact, it holds if we restrict ourselves to the class of the constructible
sets: the universe $\bold L$  consisting of the constructible sets, which 
satisfies all the axioms of set theory and also GCH (see \S3).  The
class $\bold L$ can be described as the minimal family of sets you have to
have as long as you have the same ordinals (that is, order types of linear
orders which are well founded, i.e., linear orders for which every nonempty
set has a first element). On the other hand Cohen has proved that 
you cannot prove the Continuum Hypothesis. Whereas G\"odel has ``shrunk" the
universe, Cohen has ``extended" it, adding ``generic" subsets to old
partially ordered sets (which from this perspective are called forcing
notions). 
\medskip

{\bf G\"odel}:  CH cannot be refuted.  Moreover, the {\it Generalized
Continuum Hypothesis} may hold, in fact it holds if we restrict ourselves to
the class $\bold L$ = the class of constructible sets.
\medskip

{\bf Cohen}: You cannot prove that all sets are constructible, and you
cannot even prove the weaker statement CH.
\medskip

Cohen discovered the method of {\it forcing} and used it to prove this
``independence'' result.

So the history of set theory is like that of many people - first non
self-awareness, happiness of discovering the world (not paying too much
attention to B.~Russell), then going on a crash diet for a look of thinness
($\bold L$ of G\"odel) and then discovering the comforts of growing a nice
belly (after Cohen).

In this section we treat the ``independence of CH'' not for its own sake but
as a prototype for an independence result (= proof of unprovability).
\medskip

\ub{{\stag{2.1} Dream}}:\quad  1) Find a ``forcing method'' relative to PA
which shows that PA and even ZFC does not decide ``reasonable" arithmetical 
statements, just like the known forcing method works for showing that
ZFC cannot decide reasonable set theoretic questions; even showing the
unprovability of various statements in bounded arithmetic (instead of PA) is
formidable. \nl 
2) In particular, in a way parallel to forcing, find independence 
results for classical problems in number theory; {\it for example}: 
\mr
\item "{$(\boxdot)$}"  prove that various statements formalizing the
``randomness of the prime numbers" (see below, e.g., there are infinitely
many primes of the form $n^2+1$) are independent (from PA and even ZFC and
other such statements). 
\endroster
\medskip

Of course, $(\boxdot)$ is a possible materialization of
\medskip

\ub{\stag{2.2q} Dream}:\quad  Find a theory which will formalize the thesis
that the primes behave like a random set of natural numbers gotten by
tossing a coin with probability $1/\log(n)$ for deciding if $n$ is prime.  
\medskip

Of course a parallel achievement for other classical problems in number
theory will be very welcome too.  
\medskip

{\bf Warning:} There are two major known kinds of ``independence results''.
The method of forcing can only be used to make the universe ``fatter'', not
``taller''. In technical terms: if we use forcing starting from models of
ZFC to prove that ``ZFC neither proves nor refutes statement A'' (or
equivalently, each of ``ZFC + A'' and ``ZFC + non-A'' is consistent), then
the ``consistency strengths'' of ``ZFC'', ``ZFC+A'', ``ZFC + non-A'' , are
all equal.

We should not confuse forcing with another method for proving independence:
that of consistency strength.  This means sentences like those gotten by
G\"odel's incompleteness theorem and its many illustrious descendants; some
of which are discussed below, for example:
\smallskip

``in ZFC we cannot prove that ZFC is consistent".  
\smallskip

\noindent (See more in the next section after \scite{3.1}.) Here we do not
use forcing, rather we rely on the fact that the consistency strength of
``ZFC + ZFC {\it is consistent}'' is strictly higher than the consistency
strength of ZFC alone. 

Independence results (=unprovability proofs) of this second kind (comparing
consistency strengths) are possible also over much weaker theories such as
Peano Arithmetic, a common axiomatization for the natural numbers.

Note that if we take as our base PA instead ZFC, statements like: ``every
definition $\{x \in \Bbb N:\varphi(x)\}$ defines a set'' play the role of
large cardinals.

There has been much work devoted to trying to find finitary combinatorial
statements which are equiconsistent with such ``large cardinals
statements'': Paris and Harrington \cite{PH77}, Kirby and Paris \cite{KP82}
(the Hydra), Ketonen and Solovay \cite{KS81}, and many works of Harvey
Friedman, see \cite{FOM} for much information and discussion. Friedman has
been saying for some years that all this is irrelevant to ``mainstream
mathematics'', and thinks that other mathematics relevant for his
discoveries will be developed.  So from his view it follows that Dream
\scite{2.4} below is false.

Let me stress that those two methods of forcing and of comparing consistency
strength are incomparable: though both may prove sometimes the same theorem,
they say different things; the forcing one is more specific to set theory,
gives a strong ``no answer", whereas the method of comparing consistency
strength (or large cardinals) gives an answer whose meaning may naturally
lead to debate.  It depends on our degree of confidence in the consistency
of the large cardinal considered.  Also it has a limitation as a method of
proving independence: if ZFC $+ \varphi_1$ is equiconsistent with ZFC $+
``\exists$ a strongly inaccessible cardinal" whereas ZFC $+ \varphi_2$ is
equiconsistent with ZFC, then ZFC $+\varphi_2+\neg\varphi_1$ is consistent,
but we cannot do anything about ZFC $+\varphi_1+\neg\varphi_2$.  

My ``dream'' about classical problems in number theory refers to
independence results of the first kind, such as those obtained by forcing,
where consistency strength is not increased.  Concerning the other method:
\medskip

\ub{{\stag{2.4} Dream}}:\quad  Prove that the Riemann Hypothesis is
unprovable in PA, but is provable in  some higher theory.
\medskip

What basis does my hope for this dream have? First, the solution of
Hilbert's 10th problem tells us that each problem of the form ``is the
theory ZFC $+ \varphi$ consistent" can be translated to a (specific)
Diophantine equation being unsolvable in the integers, moreover the
translation is uniform (this works for any reasonable (defined) theory,
where consistent means that no contradiction can be proved from it).
Second, we may look at parallel development ``higher up"; as the world is
quite ordered and reasonable.

Note that there is a significant difference between $\Pi_2$ sentences (which
say, e.g., for a given polynomial $f$, the sentence $\varphi_f$ saying that
for all natural numbers $x_0,\dotsc,x_{n-1}$ there are natural numbers
$y_0,\dotsc,y_m$ such that $f(x_0,\dotsc,y_0,\ldots) = 0$) and $\Pi_1$
sentences saying just that, e.g., a certain Diophantine equation is
unsolvable. The first ones can be proved not to follow from PA by
restricting ourselves to a proper initial ``segment" of a nonstandard model
of PA.  For $\Pi_1$ sentences, in some sense proving their consistency show
they are true (as otherwise PA is inconsistent).  Naturally, concerning
statements in set theory, models of ZFC are more malleable, as the method of
forcing shows.
\medskip

\ub{\stag{2.4m} Problem}:\quad Show that forcing is the unique method in a
non-trivial sense. 
\medskip

I have the impression that number theorists were generally not so excited
about the unsolvability of Diophantine equations.  Probably they reasoned
that anyway they did not hope for such a grand solution and all this has no
direct bearing on their work.  So far they are right, but several times
researchers in other fields have felt similarly (true, in a posteriori
wisdom we can explain the difference but we all do not lack this kind of
wisdom), clearly my dream contradicts this feeling.
\bigskip\bigskip

\head {\S3 Set theory and additional axioms} \endhead  \resetall \sectno=3
\medskip

\ub{G\"odel's sentences and large cardinals}:

So the sentence saying `` there is no proof of contradiction in PA (Peano
Arithmetic)" or similarly in ZFC is undecidable in this theory ( but in ZFC
we can prove the consistency of PA). Now it seems strange to believe in PA
while not believing in ``PA+ CON(PA)'' where CON(T) is the assertion that
the theory $T$ is consistent, (but where to stop?), so though we prove
independence, this is not like the proof of the independence of CH, which
leaves us with no indication what is true; a related difference is also
expressed by ``those sentences are metamathematical''.  Related to this
incompleteness method are large cardinals.  The so called ``large cardinal
axioms'' states that a ``large cardinal $\kappa$" exists which means it
resembles $\aleph_0$ in some senses in particular with properties implying
that ${\Cal H}(\kappa)$, the family of sets with hereditary closure of
cardinality $< \kappa$, form a model of ZFC and more.  Noticeable among them
is ``$\kappa$ is a measurable cardinal", which says that there is a 0-1
measure on the family of subsets of $\kappa$, which gives singletons measure
zero and is ${<}\kappa$-complete (that is, the union of $<\kappa$ sets of
measure zero is of measure zero). The first large cardinal property is
``$\kappa$ is a strongly inaccessible cardinal" which means that it is
strong limit (i.e., $\mu < \kappa \Rightarrow 2^\mu <\kappa$ ), regular 
(i.e., is not the sum of $<\kappa$ cardinals each $<\kappa$) and uncountable
(the other two properties are enjoyed by $\aleph_0$).  I feel those are
mathematical statements in set theory; some others will call them not
mathematical, but logical or set-theoretical. They reserve ``mathematical'' 
for number theory or generally what they call ``mainstream mathematics''.
\medskip

\ub{Strengthening ZFC, new axioms} 

Should we add more axioms to ZFC?  There were some options on various
grounds.  For some time using GCH was quite popular; it gives a coherent
theory for set theory of the real line (e.g., there is an uncountable set
$S$ of reals such that every uncountable subset of $S$ is not of Lebesgue
measure zero (or is not meagre), see, e.g., the writing of W.~Sierpi\'nski).
It was also used extensively in the partition calculus investigation see
Erd\"os, Hajnal and Rado \cite{EHR} and in model theory in the 1950's and
60's, e.g., the use of saturated models and in Keisler's theorem that every
two elementary equivalent models have isomorphic ultrapowers. 

It became reasonable to assume GCH (and clearly if you prove a statement
using ZFC + GCH, then you cannot refute in ZFC that statement, by G\"odel's
work on $\bold L$).  Some people look at this behaviour as adopting GCH, and
it could be argued that you do not lose much (e.g., ZFC and ZFC + GCH have
the same arithmetical consequences). Still most mathematicians, even those
who have worked with GCH do it because they like to prove theorems and they
could not otherwise solve their problems (or get a reasonable picture),
i.e., they have no alternative in the short run.  Clearly, even after
forcing was found, it seems better to prove that something follows from GCH
than just proving it is consistent; statements which we treat like this we
shall call semi-axioms (later we explicate more).  Of course, the extent to
which we consider a statement a semi-axiom is open to opinion and may change
in time.  I give statements in cardinal arithmetic a high score in this
respect.

Note that a semi-axiom may be (consistent with ZFC and) very atypical
(= the family of universes satisfying it is ``negligible") but still very
interesting as for some set of problems it gives a coherent nice
picture.

After the works of K.~G\"odel and R.~B.~Jensen on $\bold L$, when it became
clear that there is very little independence over ZFC + $\bold V = \bold L$, 
adopting $\bold V = \bold L$ as an axiom became an issue. It is useful, it
decides many problems (in particular it implies GCH so we understand
cardinal arithmetic), it also is a natural statement. Just one problem: why
the hell should it be true?   (What does true mean here? See below.)  Jensen
has thought that to prove a statement holds in $\bold L$  or just in a
universe which is canonical, has a structure (a parallel of his fine
structure theory) is better than just proving consistency.

I think so too and give $\bold V = \bold L$ and other inner models
semi-axiom status, though probably less high than Jensen does.
\medskip

J.~Kennedy wondered what I mean by ``$\bold V = \bold L$ is not true".

Some believe that compelling, additional axioms for set theory which settle
problems of real interest will be found or even have been found.  It is hard
to argue with hope and problematic to consider arguments which have not yet
been suggested.  However, I do not agree with the pure Platonic view that
the interesting problems in set theory can be decided, we just have to
discover the additional axiom.  My mental picture is that we have many
possible set theories, all conforming to ZFC.  I do not feel ``a universe of
ZFC" is like ``the sun", it is rather like ``a human being" or ``a human
being of some fixed nationality", see more in \cite{Sh:E16}.

So my meaning in saying ``why the hell should it be true", is not that it is
provably false, just as ``the national lottery in the last ten years was won
successively in turn by the nephews of the manager, so we know that 
there was cheating" is mathematically not proved.  Clearly $\bold L$ is
very special, to some extent unique, thus, the statement $\bold V = \bold L$
should get probability zero (thought not being impossible).
So $\bold L$ is certainly a citizen with full rights but a very atypical
one.  Also a typical citizen will not satisfy $(\forall
\alpha)[2^{\aleph_\alpha} = \aleph_{\alpha + \alpha +7}]$ but probably 
will satisfy $(\exists \alpha)(2^{\aleph_\alpha} = \aleph_{\alpha +
\alpha + 7})$.  However, some statements do not seem to me clearly
classified as typical or atypical.  You may think ``does CH, i.e.,
$2^{\aleph_0} = \aleph_1$ hold?" being like ``can a typical American be
Catholic". More reasonably CH has a small measure, still much much more than 
$\bold V = \bold L$.  For set theorists I will add that $\exists 0^\#$
is for me a candidate for a statement with positive measure and with a
positive measure for its negation.

What about ZFC + ``ZFC is inconsistent"?  Clearly we would not consider such
theories except that G\"odel incompleteness theorems forces us to do it; so
if $\bold V = \bold L$ is an atypical citizen, such a theory should perhaps
be a permanent illegal immigrant.
\medskip

But I am very interested in
\medskip

\ub{\stag{3.1} Dream}:\quad Find statements which yield a wonderful theory,
so will therefore be accepted as additional semi-axioms in the sense above,
which are not just other cases of those above.
\medskip

Large cardinals are certainly natural statements, as their role in finding a
quite linear scale of consistency strength on statements (arising
independently of them) shows.  But there is an important group with stronger
beliefs - the California school of set theory which holds that AD$_{{\bold
L}[\Bbb R]}$ (and relatives which we explain below) is true.  They are
interested mainly in descriptive set theory so let me first try to explain
it. The point is that we know that some sets of reals are not Lebesgue
measurable, but Borel sets are, and even $\Sigma^1_1$-sets (projections to
$\Bbb R$ of Borel subsets of $\Bbb R \times \Bbb R$).  But this is not
necessarily true for projective sets (= the sets belonging to the closure of
the family of Borel subsets of the $\Bbb R^n$-s under projection and
complement).

In particular if $\bold V = \bold L$ this is false (i.e., not all projective 
set of reals are Lebesgue measurable) and the California school believes
that generally it gives a ``false", uninteresting picture and the answers we
get in this case are ``incidental", ``artificial".  This is very reasonable.
For them, the statement AD$_{L[\Bbb R]}$, and lately the star-axiom from
Woodin \cite{Wd00}, are the remedies they adopt.  Just as AD$_{\bold L[\Bbb
R]}$ solves ``correctly" the theory of $({\Cal H}(\aleph_1),\in)$, so does
the star-axiom solve correctly the problem of $({\Cal H}(\aleph_2),\in)$.
In particular, AD$_{\bold L[\Bbb R]}$ implies that all projective sets of
reals are Lebesgue measurable and generally gives the ``right theorems of
descriptive set theory".  The star-axiom gives much more, in particular it
implies the continuum, i.e., the number of reals, is $\aleph_2$.

The argument of the California school for AD$_{\bold L[\Bbb R]}$ is that it
gives a coherent, true picture for a group of problems on the continuum,
e.g., the behavior of the projective sets and implies that the natural
problems in the theory of projective sets of reals are decidable and in the
``right way" and it is equi-consistent with suitable large cardinals,
moreover they follow from their existence.

Now they do not have an accepted catechism, different members may differ,
and even the same reasonable person may vary according to time or place, but
I heard many times that ``AD$_{L[\Bbb R]}$ is true".  More specifically, for
definitiveness, let me quote two major opinions. H.~Woodin, in a seminar at
the Mittag-Leffler institute (Fall 2000), says that the work on problems on
implications between variants of the axiom of choice (so prominent in set
theory till the sixties) has been essentially deserted, marginalized as the
axiom of choice is just accepted as true.  Similarly, he thinks the
star-axiom will be accepted as true, and research which was done on problems
it answers by other approaches will be marginalized.  In an informal
symposium in V\"a\"an\"anen's apartment around the same time, after I
explained my position on very interesting semi-axioms, J.~Steel persuasively
argued ``so having agreed that not all set theories are equal, some are more
interesting, they are better, so make one more step, is there not the best
... (star-axiom)"?  Woodin \cite{Wd00} is an excellent presentation and is
written in a non-polemic style. 

Let me stress: the mathematics of the California school is great: deep and
interesting; the program succeeded to discover a major semi-axiom.  Unlike
Woodin I feel that probably the right analogy is between the status of
AD$_{\bold L[\Bbb R]}$ and the star-axiom now with that of GCH and $\bold V
= \bold L$ earlier.  At the time they were much more informative then any
alternative, now they are not marginalized, just not ``the favorite son" any
more.

I strongly reject the California school's position on several grounds.
\mr
\item "{$(a)$}"   Generally I do not think that the fact that a
statement solves everything
really nicely, even deeply, even being the best semi-axioms (if there is
such a thing, which I doubt) 
is a sufficient reason to say it is a ``true axiom".  In
particular I do not find it compelling at all to see it as true.
\sn
\item "{$(b)$}"  The judgments of certain semi-axioms as best is based
on the groups of problems you are interested in.  For the California
school descriptive set theory problems are central.
While I agree that they are 
important and worth investigating, for me they are not ``the
center".  Other groups of problems suggest different semi-axioms as
best, other universes may be the nicest from a different perspective.
\sn
\item "{$(c)$}"  Even for descriptive set theory the adoption of the
axioms they advocate is problematic.  
It makes many interesting distinctions disappear (see more below). 
\ermn
Now I reject also the extreme formalistic attitude which
says that we just scribble symbols on paper or all consistent set
theories are equal (see above before \scite{3.1}). 
\medskip

{\stag{3.2} \ub{ Dream}}:\quad Find universes with a different but
interesting theory of projective sets (i.e., universes in which AD$_{\bold 
L[\Bbb R]}$ fails). 
\medskip

Note that $\bold V = \bold L$ gives a very clear theory but one where, e.g.,
we have a projective well ordering of reals (hence a non-Lebesgue-measurable
set). Of course, AD$_{\bold L[\Bbb R]}$ is much more reasonable than the
star-axiom which is very special though not as special as $\bold V = \bold
L$; $\bold V = \bold L$ has probability zero whereas AD$_{\bold L[\Bbb R]}$
has higher ``probability", \ub{but} it excludes no less interesting
universes.  Even for set theory of the reals, the possibility of the
continuum being real valued measurable, is not less important (see on it
later, anyhow it implies that the continuum is large).

So according to the view presented here, AD$_{{\bold L}[\Bbb R]}$ is
certainly a semi-axiom: a beautiful theory with fascinating theorems built
on it.  I may even agree it has ``positive measure", is not atypical.  But
``there is no inner model with suitable large cardinal" seems to me also of
positive measure.  Also the large cardinals themselves seems to me of
positive, co-positive measure (decreasing with largeness, of course).

I think $\bold V = \bold L$ passes this criterion of being a great
semi-axiom for many problems. It is extremely helpful in building examples,
e.g., it gives a very coherent theory on an important group of problems in
Abelian group theory (see Eklof and Mekler \cite{EM}, \cite{EM02}).

The descriptive set theorists have reasonable reasons to reject the axiom
$\bold V = \bold L$, but are there not similar grounds for rejecting the
continuum being $\aleph_n$, $n$ a natural number, hence rejecting the star
axiom?  Various combinatorial properties of the continuum follows from
$2^{\aleph_0} = \aleph_n$ which are akin to having a definable well order
(on $(n+1)$-place function from reals to reals), see \cite{EHR}.  Moreover,
this applies even to problems with descriptive set theory flavours.  For
example I am currently interested in a very explicit definition of an
Abelian group $G_n$ (for $n$ a natural number; we can represent the set of
elements as the reals, and $x \in G$, $x + y$, $-x$ are not just Borel but
even $F_\sigma$).  Now this Abelian group is free iff $2^{\aleph_0} \ge 
\aleph_{n+1}$ (see \cite[\S5]{Sh:771}).  It seems to me that its being free
is accidental, just the continuum is not large enough for us to see how
complicated it is, to see its true nature.  So assuming $2^{\aleph_0} \le
\aleph_n$ is artificial, wrong, just as $\bold V = \bold L$ is wrong for
``are projective set Lebesgue measurable".

It may be more interesting to consider a family of problems on the
continuum: investigate cardinal invariants.  We may measure the continuum
not only just in size but in other ways, like non(null) = the minimal
cardinality of a non-null set or non(meagre) = the minimal cardinality a
non-meagre set (= not first category set), ${\frak d} = \{|{\Cal F}|:{\Cal
F}$ a family of functions from $\Bbb N$ to $\Bbb N$ such that every such
function is bounded by one of them$\}$.  There is a myriad of such measures
many of them important in many directions (see Bartoszy\'nski \cite{Baxx},
Blass \cite{Bsxx}) naturally they are uncountable but $\le 2^{\aleph_0}$. 

If the continuum is $\le \aleph_2$, we will not have (by a trivial
pigeonhole principle) many relations, as no three can be simultaneously
distinct. This seems to me not being less artificial than the answers to
descriptive set theoretic problems in $\bold L$.
\medskip

\ub{\stag{3.3} Dream}:\quad Find a consistent axiom which gives a coherent true
picture for cardinal invariants of the continuum, e.g., it implies
that any two cardinal invariants which are ``not necessarily equal" will be not
equal so necessarily the continuum is large (as there are many such
invariants, see Goldstern and Shelah \cite{GoSh:448}).
\medskip

Such an axiom is not necessarily unique since for many pairs of cardinal
invariants both inequalities are consistent.

This is naturally connected to 
\medskip

\ub{\stag{3.3a} Dream}:\quad  1) Develop a theory of iterated forcing for a
large continuum as versatile as the one we have for $2^{\aleph_0} = 
\aleph_2$, and to a lesser extent, for $2^{\aleph_0} = \aleph_1$ (see
\cite{Sh:f}, and a different approach in Woodin \cite{Woxx}); but see 
\scite{4.2a}. \nl
2) So far we have many independence results and few theorems but hopefully
(see the rubble thesis in \S5) as in the case of cardinal arithmetic a
positive theory will appear.
\medskip

Another approach, not disjoint to the descriptive set theory one, is to
adopt large cardinal axioms, the argument being that they are the natural
extension of how we arrive to ZFC.

A large cardinal axiom is one stating there is an uncountable cardinal
$\lambda$ which satisfies a property satisfied by $\aleph_0$ like ``strongly
inaccessible", see above.  A picture justifying it is defining by induction
on the ordinal $\alpha$ the set $\bold V_\alpha$ by $\bold V_0 =\emptyset$,
$\bold V_\alpha =\bigcup\{{\Cal P}(\bold V_\beta):\beta < \alpha\}$. The
intuitive argument is that if we ``wait long enough", i.e., for some
$\alpha$ large enough, whatever is not forbidden to happen will happen, so
there are $\alpha$ which reflect well what the whole universe $\bold V =
\bigcup\{\bold V_\alpha:\alpha$ an ordinal$\}$ satisfies.  This seems to me
reasonable (but see later). A motivation was the hope that this would decide
classical problems about relatively small sets; it has failed for CH (though 
above a supercompact cardinal cardinal arithmetic is simple) but succeeded
for problems like ``every projective set of reals is Lebesgue measure" (see
above). In my view the analogy of arriving at large cardinals with ZFC is
problematic:  we arrive at ZFC by considering natural formations of sets
(the set of natural numbers, taking Cartesian products and power sets); even  
the first strongly inaccessible has no parallel justification.  If you go
higher, up in the large cardinal hierarchy, the justification for their
existence is decreasing so large cardinal axioms are great semi-axioms but
not to be accepted as true.

J.~Kennedy has wondered where do we make the comparison between various set
theories, and by what criterions.  Now we can do it informally, Cantor has
understood set theory quite well and understood the Continuum Problem
without sticking to a formalization.  Alternatively, we may work within a
``bare bones set theory", just enough to formalize first order theories,
proofs, and say, having the completeness theorem.  We may well agree we are
in a universe which is set theory and discuss it without apriori having a
common agreement on all its properties.

What are our criterions for semi-axioms?  First of all, having many
consequences, rich, deep beautiful theory is important. Second, it is
preferable that it is reasonable and ``has positive measure".  Third, it
is preferred to be sure it leads to no contradiction (so lower consistency
strength is better).  

Naturally, those are conflicting hopes. So $\bold V=\bold L$ is preferable
on GCH as it has more consequences but as GCH is much more reasonable and
still has a very large set of consequences, it is worthwhile to use it, and
to try to prove from it what earlier was proved from $\bold V = \bold L$. 

Is the position presented here consistent?  Can I on the one hand be opposed
to ``the true, unique, set theory" and on the other hand not give equal
place to all consistent set theories (measure zero to $\bold V = \bold L$)?
Having pointed out a case where we accept such a position (typical
American), I think I have shown the consistency of this position; in fact,
considering any biography of a person is no better than the ``typical
American", as any interesting statement about a person normally has such
character.

It may help the reader to note that my interest in mathematical logic has
been mathematical, seeking generality, rather than philosophical.  So we
have not seriously addressed issues which do not seem to make a difference
for mathematics: you may be a pure Platonist, telling me that there is a
unique universe of set theory about which we know no more than ZFC, \ub{or}
you believe in going to the sources, ``would Cantor have such a thing in
mind", \ub{or} believe that the axioms of set theory reflect the
relationship between the two hemispheres of the human brain \ub{or} that we 
discover (rather than invent it) by our research, \ub{or} you may think that
starting with some mental picture, we change it as our knowledge increases
and converge to a ``natural theory" (up to presentation).  Very important
and fascinating intellectual issues but will not make a real difference for
the discussion here.
\bigskip \bigskip 

\head {\S4 ZFC indecisiveness} \endhead  \resetall \sectno=4
\medskip

Again, I certainly think large cardinals are natural notions, an
indispensable part of our knowledge.  For me the most 
compelling reason for their being important is the linear order
phenomena and their roles in, see below.

Two statements $\varphi,\psi$ are equi-consistent if we can prove (even
in PA) that ZFC $+ \varphi$ is consistent iff ZFC $+ \psi$ is
consistent,  i.e., ZFC $\models \neg \varphi$ iff ZFC $\models \neg
\psi$.  Now:
\mr
\item "{$(A)$}"  For a plethora of classical problems $\varphi$ in set
theory, large cardinal statements $\varphi_{LC}(x)$
were found such that $\varphi$ and $\exists \kappa \varphi_{LC}(\kappa)$ are
equi-consistent or at least were sandwiched between two large cardinal
statements; in fact, for most of the cases, if set theorists were really
interested, they get equi-consistency results; for examples
see below.
\sn
\item "{$(B)$}"  The large cardinals are linearly ordered, i.e., for two
such properties $\varphi^1_{LC}(\kappa)$, $\varphi^2_{LC}(\kappa)$ we can
almost always prove for some $\ell \in \{1,2\}$ that
{\roster
\itemitem{ $(\boxdot_{\varphi^1_{LC},\varphi^2_{LC}})$ }   
CON(ZFC $+ \exists \kappa \varphi^\ell_{LC}(\kappa))$ implies
CON(ZFC $+ \exists \kappa\varphi^{3 - \ell}_{LC}(\kappa))$
\nl
that is, relative consistency gives a linear order (transitivity holds
trivially, the interesting phenomena is the comparability).
\endroster}
\ermn
In fact, usually when $(\boxdot_{\varphi^1_{LC},\varphi^2_{LC}})$ holds,
then the result is more explicit: if $\kappa_2$ satisfies
$\varphi^2_{LC}(-)$, then we can find $\kappa_1 < \kappa_2$ such that
$\varphi^1_{LC}(\kappa_1)$ (at least in a smaller universe).

The way we prove from the consistency of ZFC $+ \varphi_{LC}(\kappa)$ the
consistency of ZFC $+ \varphi$ is usually by forcing.  The way we prove
from the consistency of ZFC $+ \varphi$ the consistency ZFC $+
\varphi_{LC}(\kappa)$ is by inner models like $\bold L$ (i.e., we shrink the
universe of sets putting in only ``necessary sets").
\medskip

\ub{\stag{3.4} Dream}:\quad  Explain the phenomena of linearity of
consistency strength (and of large cardinal properties): to what extent is
it the outcome of the history of set theory, the history of mathematics and
human perceptions; of course, we hope this will lead to exciting
mathematical discoveries.
\medskip

We may give a naive answer along the following lines:

For a large cardinal statement $\varphi = \varphi_{LC}(\kappa)$, let
$\alpha_\varphi = \text{ Min}\{\alpha:\alpha$ a countable ordinal and
there is a countable set $\bold V^*$ which is transitive (i.e., $x\in y\in  
\bold V^* \Rightarrow x \in \bold V^*$) and an ordinal is in $\bold V^*$
iff it is $<\alpha$, and $\bold V^*$ is a model of ZFC $+\exists\kappa
\varphi(\kappa)\}$.

So we actually investigate $\alpha_{\varphi^1} \le \alpha_{\varphi^2}$,
which is linear.  But this does not explain why the statements we are
interested in so far tend to be equi-consistent with such statements; we can
produce counterexamples but they are artificial; (see \cite{Sh:170} and
lately I heard on recent lectures of Woodin).  Like several other dreams
here it is foggy, being an answer may well be disputable, may have no
solution or several.

As an illustration consider the classical case of real-valued measurable
cardinals. 

We know that, unfortunately, we cannot have a measure on all sets of
reals as good as Lebesgue measure, but mathematicians tried to remedy
this.  One way is to omit the requirement that the measure $\mu$ is
preserved by translations (so restricting the attention to the unit
interval, any $A\subseteq [0,1]_{\Bbb R}$ is given a measure which is a
real number $\mu(A) \in [0,1]_{\Bbb R}$, the measure is $\sigma$-additive
(if $A_n$ are pairwise disjoint, then $\mu(\bigcup\{A_n:n \in \Bbb N\}) =
\sum\{\mu(A_n):n \in \Bbb N\}$), $\mu([0,1]_{\Bbb R})=1$, $\mu(\emptyset)
= 0$ and $\mu(\{x\}) = 0$). If $\lambda = \text{ Min}\{|A|:\mu(A)=1\}$ we
say that $\lambda$ is real-valued measurable since we can copy the measure
from the set $A$ to $\lambda$. 

Investigations on such measures on bigger sets lead to measurable cardinals;
a cardinal $\lambda$ is called measurable if on some set $A^*$ of
cardinality $\lambda$ there is a $\{0,1\}$-measure, i.e., $\mu:{\Cal P}(A^*) 
\rightarrow \{0,1\}$ with $\mu(\emptyset)=0,\mu(A^*)=1,\mu(\{x\})=0$, which
is even $(< \lambda)$-additive (i.e., if $A_i \subseteq A^*$ for $i \in I$
are pairwise disjoint and $|I| < \lambda$ then $M(\bigcup\{A_i:i \in I\})= 
\text{ max}\{\mu(A_i):i\in I\}$). This is one of the most important large
cardinal properties.  In particular, its existence is unprovable in ZFC.

R.~Solovay has shown that if there is a measurable cardinal then in some
universe (obtained by a forcing extension) the continuum is real-valued
measurable, thus showing the relative consistency.  The reader may hope
for a proof without a large cardinal or at least with a smaller large
cardinal, but (s)he will hope in vain.  Solovay proves that if there is a
measure as above on the family of all subsets of $[0,1]_{\Bbb R}$,
then if we shrink the universe (getting a so-called inner model) we
can get a universe with a measurable cardinal.  True, the two
properties have some affinity to begin with, but this serves well as a
measuring stick.  To have AD$_{\bold L[\Bbb R]}$ requires much larger
cardinals, to have ``all projective sets are Lebesgue measurable"
requires much less.  Moreover, this analysis leads sometimes to ZFC
results, e.g., close to my heart is: if $A \subseteq [0,1]_{\Bbb R}$
has positive outer Lebesgue measure, then we can find pairwise disjoint
sets $A_n \subseteq A$ for $n \in \Bbb N$, each with the same outer
Lebesgue measure as $A$, see Gitik and Shelah \cite{GiSh:582}.

Note that consistency strength gives us independence and more, hence
necessarily is less versatile than forcing; {\it consistency strength can
give one side of independence}, forcing usually can give both; but, of
course, the forcing may well start with large cardinals.

The notable hole in the program above is the supercompact, a very high large
cardinal which has been quite widely used in consistency proofs.  The hole
is that we do not know how to prove the consistency of ZFC + ``there is a
supercompact cardinal" from ``low level statements".
\medskip

\ub{\stag{3.5} Problem}:\quad Find an inner model for 
supercompactness.  Without this we may always hope to replace it as an
assumption for consistency results by smaller large cardinals (as had been
done in some important cases).
\medskip

\centerline {$* \qquad * \qquad *$}
\medskip

The theory of forcing and, in particular, of iterated forcing seems
important to me.  There is a reasonable amount of such theory for forcing
notions for the continuum; in the beginning we had one kind of such forcing:
Cohen forcing.  Now we know more, though far from what we like to know as
mentioned above.  The case of continuum $> \aleph_2$ is a case where we miss
much (in proof or in forcing, probably in both).  For larger cardinals we
know less, e.g., in work on the $2^\mu$ for $\mu$ singular, very
sophisticated forcings were discovered (see writing of M.~Gitik, M.~Magidor,
W.~Mitchell and H.~Woodin).
\medskip

\ub{\stag{3.7} Question}:\quad  Can we find a theory of iterated forcing
parallel to the one we have for the continuum, for the following purposes.
\mr
\item "{$(a)$}"  For forcing as developed for the singular cardinal
problem. 
\sn
\item "{$(b)$}"  For other cases see more in \cite{Sh:666}, on several
missing theories of this kind. 
\endroster
\medskip

\ub{\stag{3.8} Dream}:\quad  Find additional methods for independence (in
addition to forcing and large cardinals/consistency strength), or prove the
uniqueness of these methods.  
\medskip

\ub{Dividing lines}:
  
``Dividing line" or ``watershed line" means here a property such that both
it and its negation has strong consequences; hence it helps in proofs by
cases; dividing lines were fruitful in my work in model theory.

So $2^{\aleph_0} = \aleph_{\omega^3 + \omega+1}$ does not look like a
good dividing line (or even a property).  What about CH, i.e., $2^{\aleph_0}
= \aleph_1$?  This statement has many consequences and hence it is an
impressive assertion but its negation seemingly does not, so it is an
important semi-axiom but not a good dividing line.
\medskip

\ub{\stag{dl.1} Dream}:
\mr
\item "{$(a)$}" Find a real significance for $2^{\aleph_0}=\aleph_{753}$,
\sn
\item "{$(b)$}" or for $2^{\aleph_0} = \aleph_{\omega^3 + \omega +5}$,
\sn
\item "{$(c)$}"  show that all values of $2^{\aleph_0}$ which are $>
\aleph_2$ are similar in some sense (or at least all values $\aleph_n>
\aleph_2$, all regular $\aleph_\alpha>\aleph_{\omega_1}$ or whatever). 
\endroster
\medskip

\ub{\stag{dl.2} Dream}:\quad Can we find important dividing lines  and
develop a theory for combinatorial set theory?
\medskip
 
Now Jensen has a different dream (I do not believe that it will
materialize). 
\medskip

\ub{\stag{dl.3} Dream}:
Find a super-duper ``inner model'', so large that basically it always
behaves like $\bold L$ in a universe with no $0^{\#}$; so we can answer
questions by translating our problems to it were we will have fine structure
to help us.
\medskip

Of course, his optimism has something to do with his success with $0^{\#}$;
if there is no such real then the universe is ``close to $\bold L$",
otherwise it is very far from it.

Considering inner models there are good dividing lines for descriptive set
theory. Consider the statement `` for every real $r$, $r^{\#}$ exists" (see,
e.g., Jech \cite{J}).  Why is it a good dividing line? If $r$ is a
counterexample then descriptive set theory, the case where real parameters
are allowed or just with parameter $r$, is very much like the one in $\bold
L$, so all the questions traditionally asked are answerable, though, many
say, in the ``wrong, uninteresting'' way. If it holds, we have determinacy
for $\Sigma^1_1$ games (which gives nice consequences for ``low level''
projective sets of reals).

Let us consider two candidates for being dividing lines, unfortunately not
impressive ones.  The statement ``there is a nice normal filter on
$\omega_1$", see \cite[Ch.V]{Sh:g}, helps to prove cardinal arithmetic
inequalities; if it holds it is used to define rank on functions from
$\omega_1$ to the ordinals, if it fails the universe is similar enough to
$\bold L$ hence cardinal arithmetic is trivial by Dodd and Jensen
(\cite{DJ1}); see a use in \cite{ShSt:419}.

So far values of the continuum have not appeared to be good dividing lines,
but we may look at another candidate.  (Whereas $2^{\aleph_0} = \aleph_1$
has many consequences, $2^{\aleph_0} = \aleph_2$ has few consequences but
many consistency results.)   

The second candidate we consider is the cardinal arithmetic equality $2^\mu
= \mu^+$, where $\mu$ is a strong limit singular cardinal, e.g., $\mu =
\beth_\omega$, see \S5. 
\medskip

\ub{Forcing for $\bold L$: base theory = ZFC:\quad highly developed
``forcing'' technology}: 
  
Forcing (as has been developed) was very successful to show that there are
many problems which (like the continuum hypothesis) cannot be decided in set
theory; moreover, this has become a method, which for a typical set
theoretic problem, gives us a reasonable way to prove its consistency and
its independence.  This indicates that the axioms of set theory (ZFC) are
weak, usually not able to decide given questions.

However, these problems (such as CH) which can be attacked by forcing can
typically be decided by the single additional axiom that {\it all sets are
constructible}, $\bold V = \bold L$.

(G\"odel proved that $\bold L$ is a model of ZFC+GCH, Jensen solved many
specific problems and developed general methods under $\bold V = \bold L$).

So this leads naturally to
\medskip

\ub{\stag{3.10} Dream}:\quad Can we find a method parallel to forcing for
$\bold L$ (i.e., for the usual axioms of set theory + every set is
constructible).  
\medskip

Such a method would enable us to prove that certain statements are
independent of ZFC + $\bold V = \bold L$; the current forcing results only
give independence over ZFC.

A major preliminary obstacle to this dream is the lack of a good candidate
to be a test problem, since so many questions have already been settled
under the assumption $\bold V = \bold L$.

A negative answer, explaining why ``a large body of set theoretic problems
is decidable'' would be marvelous too (this would give ``quasi-decidability"
for $\bold L$).
\bigskip\bigskip

\head{\S5 Is ZFC really so weak?} \endhead  \resetall \sectno=5
\medskip 

In this section we deal with a (relatively) new perspective on cardinal
arithmetic.  
\medskip 

\ub{\stag{4.1} The not so poor Thesis}:\quad The view that {\it ZFC is a
deficient theory, since it does not decide so many basic questions, in
contrast to classical theories} is one sided; in fact,  a ``random''
question in, say, number theory is similarly hopeless as far as
answerability is concerned, certainly practically and, by G\"odel's 
incompleteness theorem + the undecidability of solvability of
Diophantine equations, I think also fundamentally. 
\medskip 

Set theory has actually an {\it advantage} over many other fields of
mathematics. When we seem to be unable to prove or refute something, we have
strong methods to try to show that a proof or refutation may be {\it
impossible}; such results, in addition to being interesting in their own
right, also help to clear the air, directing us to what actually {\it is}
decidable, discarding the undecidable ones.

It also tells us, what kind of problems are, at least relatively, of the
answerable kind.
\medskip 

\ub{\stag{4.2} The Rubble Removal Thesis}:\quad  
Methods for proving independence, in addition to their intrinsic value, work
for us like a sieve: when we have a myriad of problems in some directions,
and we have tried to prove independence (and many times this results in
discarding most of them), we are left with strong candidates for theorems of 
ZFC.
\medskip 

In this connection we may hope (compare with \scite{3.3}):
\medskip 

\ub{\stag{4.2a} Dream}:\quad  Our problem in proving consistency results for
continuum $> \aleph_2$ comes from the existence of some positive
theory which will become trivialized if the continuum is too
small, i.e., most of the theorems which we will prove become trivialized.
\medskip 

The discussion in sections 2,4 may support the impression that ``all is
independent in ZFC", this is not groundless but also not the whole truth.
We shall now concentrate on cardinal arithmetic not assuming specialized
knowledge (see \S0 for basic definitions).

For an even more leisurely explanation on cardinal arithmetic and pcf
for the general mathematical audience see \cite{Sh:E25}; a book devoted to
this subject is \cite{Sh:g}. 

Cardinal arithmetic is a good example for \scite{4.2}; after proving there
is nothing more to say on $2^{\aleph_0},2^{\aleph_1},2^{\aleph_{\alpha
+1}}$, we had found that we can say something about $2^{\aleph_{\omega_1}}$
and even on $2^{\aleph_\omega}$ and even on $\dsize \prod_{n < \omega}
\aleph_n$.  In fact, the Thesis of \cite{Sh:g} is that there are two
separate phenomena. The first one is the behaviour of $2^\lambda$ for
$\lambda$ regular (mainly $\lambda=\aleph_0,\aleph_{\alpha +1}$) for which
everything is independent. The second one is the cofinality problem, the
domain of pcf theory, which appears later.  

In other words, looking at the bright side, we know all
the true rules  we can know on $2^\lambda$ for regular $\lambda$ (no more
rules than the classical ones: it is non-decreasing, i.e., $\lambda\le\mu 
\Rightarrow 2^\lambda \le 2^\mu$ and cf$(2^\lambda)>\lambda$; on cf see
below). This leaves us with the singular cardinals like $\aleph_\omega =
\sum\{\aleph_n:n$ a natural number$\}$ 

First it has been ``clear" that the case of singular cardinals would be
similar to that of the regular cardinals just more complicated, a
``technical problem"; in fact independence results were found for singular
cardinals, however using large cardinals. Second, it was proved that there
are some limitations. Third, it was proved that large cardinals are
necessary, using the theory of inner models, (see Magidor \cite{Mg1},
\cite{Mg2} for independence results; Silver \cite{Si}, Galvin and Hajnal
\cite{GH}, \cite{Sh:111}, \cite[Ch.XIII]{Sh:b} on limitations on
$2^\lambda$; Devlin and Jensen \cite{DeJ}, Dodd and Jensen \cite{DJ1} on
inner models; for more on forcing see in the writings of M.~Gitik,
M.~Magidor, W.~Mitchell, H.~Woodin; for more on inner models see in the
writings of M.~Gitik, W.~Mitchell, J.~Steel).  
 
A thesis of \cite{Sh:g} is that
\medskip

\ub{\stag{4.3} Thesis}:  [``Treasures are waiting for you'']\quad 
There are many laws of (infinite) cardinal arithmetic concerning
exponentiation; in the past there seemed to be few and scattered ones
because we have concentrated on $2^\lambda$, but if we deal with relatively
small exponent and large base, there is much to be discovered; see more
below.
\medskip

Cardinal arithmetic investigations have concentrated on the function
$\lambda \mapsto 2^\lambda$ which had good reasons but made us ignore other
directions. Even after the independence results, researchers tended to be
influenced by remnants of GCH, e.g., the concentration on $2^\mu$ for $\mu$
a strong limit singular cardinal. 
\medskip

\ub{\stag{4.4} Dream}:\quad Find all the laws of (infinite) cardinal
exponentiation. 
\medskip

\ub{pcf theory}:

Close to my heart is
\medskip

\ub{\stag{4.5} Thesis}:\quad Cardinal arithmetic is loaded with consistency
results because we ask the wrong questions; the ``treasures'' thesis  is not
enough; we should replace cardinality by cofinality as explained below (pcf
theory).  
\medskip

\definition{\stag{4.6} Definition} 1) For a partially ordered set $P$ let
cf$(P)$, the cofinality of $P$, be 
$$
\min\{|Q|:Q \subseteq P \text{ satisfies } 
(\forall x \in P)(\exists y \in Q)[x \le_P y]\}.
$$
2) cf$([\lambda]^\kappa) = \text{cf}([\lambda]^\kappa,\subseteq)$ is
cf$({\Cal P})$ when for some set $A$ of cardinality $\lambda$, ${\Cal P}$  
is the family of subsets of $A$ of cardinality $\le \kappa$ partially
ordered by inclusion. \nl 
3)  We usually identify any cardinal $\lambda$ with a linear order, in fact
a well ordering  of this cardinality; every initial segment of it has
smaller cardinality. \nl 
4) For a cardinal $\lambda$, let cf($\lambda$) be Min$\{|J|:\lambda =
\sum\limits_{t\in J}\lambda_t$ for some $\lambda_t<\lambda$ for $t\in
J\}$, this is compatible with the definition of cf in part (1); we call
$\lambda$ regular if cf$(\lambda) = \lambda$ and singular otherwise; recall
that $\aleph_0$ and every successor cardinal are regular whereas
$\aleph_\omega=\sum\{\aleph_n:n$ a natural number$\}$ is the first
singular cardinal, and cf$(\lambda)$ is always a regular cardinal.\nl
5) A partially ordered set $P$ is said to have true cofinality
$\lambda$ and we write tcf$(P) = \lambda$ \ub{if} ($\lambda$ is a regular
cardinal and) there is a $\le_P$-increasing sequence $\langle p_i:i <
\lambda\rangle$ such that $(\forall q \in P)(\exists i)(q \le_P p_i)$; note 
that $\lambda$ is unique.  (Can $P$ fail to have true cofinality?  Yes,
e.g., if it is $P_1 \times P_2$, ordered coordinate-wise, $P_1,P_2$
have true cofinalities but different ones).
\enddefinition
\medskip

\ub{\stag{4.7} Convention}:\quad Let ${\frak a}$ denote a set of regular
cardinals $> |{\frak a}|$, also let ${\frak b},{\frak c}$ denote such sets.   
\medskip

\definition{\stag{4.8} Definition}\quad  1) For an ideal $J$ on ${\frak
a},<_J$ is the partial order on $\prod{\frak a}$ defined by 
$$f<_J g\ \ \Leftrightarrow\ \ \{\theta \in {\frak a}:f(\theta)<
g(\theta)\}={\frak a}\mod J.$$  
If $J = \{\emptyset\}$, then we write $<$. \nl
2)  pcf$({\frak a}) = \{\text{tcf}(\prod{\frak a},<_J):J$ an ideal on
${\frak a}$ and $(\prod{\frak a},<_J)$ has true cofinality$\}$\nl
(it is enough to consider maximal ideals). \nl
3) pcf$_{\theta\text{-complete}}({\frak a}) = \{\text{tcf}(\prod{\frak
a},<_J):J$ is a $\theta$-complete ideal on ${\frak a}$ and 
$(\prod{\frak a},<_J)$ has true cofinality$\}$. \nl
4) If $\mu$ is singular, i.e., $\mu>\text{cf}(\mu)$, let\nl
pp$(\mu)=\sup\{\text{tcf}(\prod {\frak a},<_J):{\frak a} \subseteq \text{
Reg} \cap \mu,\sup({\frak a}) = \mu,J$ an ideal on ${\frak a}$ such that
$\lambda<\mu\Rightarrow {\frak a}\cap\lambda\in J$ and $(\prod {\frak
a},<_J)$ has true cofinality$\}$. \nl
5) pp$_{\theta\text{-complete}}(\mu)$ is defined similarly
restricting ourselves to $\theta$-complete ideals. 
\enddefinition
\medskip

\ub{\stag{4.9} Thesis}:  [pp vs power set]\quad 
The power set function $\lambda \rightarrow 2^\lambda$ on regular cardinals
is totally independent (the only rules are that it is not decreasing and
cf$(2^\lambda ) > \lambda$).  It is like hair colour today, easily
manipulated, whereas  pp($\lambda$), and pcf(${\frak a}$) are like the
skeleton of set theory, not totally immune to ``plastic surgery'' (i.e.,
forcing starting with large cardinals) but at great price (and pains). 
\medskip

So the chaotic behaviour of cardinal arithmetic comes from the static noise
of the interference of two different phenomena:

One, the mapping $\lambda \rightarrow 2^\lambda$ for regular cardinals which
actually is very well understood: we know all the rules; anything fulfilling
them is permissible. The second speaks about pcf theory, there are many
mysteries that have not disappeared but much is decided in ZFC. The claim
that almost all is independent was wrong; the picture is more balanced.

So ``the armies of god (fighting for resolution in ZFC)" and ``the armies of
the devil (trying to prove independence)" have advanced much and arrive to a
new stand-off.  So we should reformulate \scite{4.4}.
\medskip

\ub{\stag{4.10} Dream}:\quad Find the laws of (infinite) cardinal
exponentiation, under pcf interpretation.  
\medskip

Note that pcf$({\frak a})$ replaces the cardinality product $\prod {\frak
a}$ by a spectrum of possible cofinalities, a phenomena which has many
honourable precedents (e.g., in decomposition into primes in algebraic
extensions of $\Bbb Z$).  On the calculus of pcf we know some rules:
\medskip

\proclaim{\stag{4.11} Theorem}\quad  1) {\rm pcf}$({\frak a})$ includes
${\frak a}$ and has cardinality $\le 2^{|{\frak a}|}$ (and not merely
$2^{2^{|{\frak a}|}}$ which is the obvious upper bound being the number of
ultrafilters on  ${\frak a}$).  \nl
2) {\rm pcf}$({\frak a})$ has a maximal member {\rm max pcf}$({\frak a})$
which is equal to {\rm cf}$(\prod{\frak a})$. \nl 
3) If ${\frak b} \subseteq {\text{\rm pcf\/}}({\frak a})$ and $|{\frak b}| 
< {\text{\rm Min\/}}({\frak b})$, \ub{then} {\rm pcf}$({\frak b})\subseteq
{\text{\rm pcf\/}}({\frak a})$.  
\endproclaim
\medskip

In short, {\rm pcf}$({\frak a})$ is not as large as we may suspect, it has a
last element which is a reasonable measure of $\prod{\frak a}$, in fact it
is cf$(\prod{\frak a},<_J)$ for $J$ the trivial ideal $\{\emptyset\}$.
Moreover, {\rm pcf} essentially acts like a closure operation, e.g.,
it is increasing.

Moreover
\medskip

\proclaim{\stag{4.12} Theorem}\quad  1) If a sequence $\langle \lambda_i:i<  
\aleph_1\rangle$ is increasing continuous with limit $\lambda$, \ub{then}
for some closed unbounded set $C \subseteq \aleph_1$ we have {\rm max} 
${\text{\rm pcf\/}}\{\lambda^+_i:i \in C\}=\lambda^+$. \nl
2) [\ub{Locality}]:  If ${\frak b} \subseteq \text{\rm pcf\/}({\frak
a}),|{\frak b}| < {\text{\rm Min\/}}({\frak b}),\theta \in \text{\rm
pcf\/}({\frak b})$, \ub{then} for some ${\frak c} \subseteq {\frak b}$ of
cardinality $\le |{\frak a}|$ we have $\theta \in \text{\rm pcf}({\frak
c})$. \nl 
3) [\ub{No hole}]: If ${\frak a}$ is an interval of the class of
regular cardinals (i.e., ${\frak a}=\{\aleph_{\alpha +1}:\alpha_*\le\alpha
<\beta\}$ so necessarily $\beta<\aleph_{\alpha_*}$), \ub{then} {\rm
pcf}$({\frak a})$ is an interval too (necessarily end-extending ${\frak a}$,
i.e., {\rm pcf}$({\frak a})$ is an initial segment of $\{\aleph_{\alpha
+1}:\alpha \ge \alpha_*\}$). 
\endproclaim
\medskip

What does this mean?  The closure operation inside {\rm pcf}$({\frak a})$
has ``character" at most the cardinality of ${\frak a}$, and there are some
forms of continuity and convexity.
\medskip

\demo{\stag{4.14} Observation}\quad  $\lambda^\kappa = 2^\kappa + 
\text{ cf}([\lambda]^\kappa,\subseteq)$.
\enddemo
\medskip

In general
\medskip

\ub{\stag{4.15} Thesis}:\quad  The natural measures of $[\lambda]^\kappa$ can be
expressed by cases of pp and of $2^\lambda$ for $\lambda$ regular.
\medskip

The measures which come to my mind are $\lambda^\kappa$,
cf$([\lambda]^\kappa,\subseteq)$ from \scite{4.6}(1), and 
$$
\align
\lambda^{[\kappa]}=\text{ Min}\{|{\Cal P}|:&{\Cal P}\subseteq
[\lambda]^\kappa\text{ and every subset of }\lambda\text{ of cardinality
 }\kappa\\
  &\text{ is the union of }< \kappa\text{ members of }{\Cal P}\},
\endalign
$$
$$
\align
\lambda^{<\kappa>} = \text{ Min}\{|{\Cal P}|:&{\Cal P} \subseteq
[\lambda]^\kappa\text{ and every subset of } \lambda \text{ of
cardinality } \kappa \\
  &\text{ is included in the union of } < \kappa \text{ members of }
{\Cal P}\}
\endalign
$$

$$
\lambda^{<\kappa>_{tr}} = \sup\{|\text{lim}_\kappa(T)|:T \text{
is a tree with } \le \lambda \text{ nodes and } \kappa \text{
levels}\}.
$$

For example
\medskip

\proclaim{\stag{4.13} Theorem}\quad {\rm cf}$([\aleph_\omega]^{\aleph_0},
\subseteq) = \text{\rm pp}(\aleph_\omega) = \text{\rm max pcf}\{\aleph_n:n <
\omega\}$. 
\endproclaim
\medskip

So pcf is a closure operation with some rules; those listed above were
enough to prove a result, popular among my works: pp$(\aleph_\omega) <
\aleph_{\aleph_4}$; i.e. by \scite{4.13} it is enough to investigate
pcf$\{\aleph_n:0 < n < \omega\}$ which, if the desired inequality
$({\text{\rm pp\/}}(\aleph_\omega) \ge \aleph_{\aleph_4})$ fails, includes
$\{\aleph_{\alpha+1}:\alpha < \aleph_4\}$.

Now pcf is a closure operation with several rules, \ub{but} there is much we
do not know. Where are the new lines between what is known and what we do
not know? Locality may be a poor substitute to a positive answer to: 
\medskip

\ub{\stag{4.16} Question}:  Is pcf$({\frak a})$ always of cardinality $\le
|{\frak a}|$? 
\medskip

In the scale of the problems (in this direction) which we do not know the
answers, this question seems to me to lie in the middle.

Less hopeless for forcing are:
\medskip

\ub{\stag{4.17} Question}:\quad Show the consistency of the failure of the
Weak Hypothesis (WH), which means:
\mr
\item "{(WH)$_1$}"  for every cardinal $\lambda$ the following set is
finite: \nl 
$\{\mu < \lambda:\mu$ is singular of uncountable cofinality that is
$\aleph_0 < \text{ cf}(\mu) < \mu$ and 
pp$_{\aleph_1\text{-complete}}(\mu) \ge \lambda\}$
\sn
\item "{(WH)$_2$}"  for any cardinal $\lambda$ the following set is
finite or at most countable:\nl
$\{\mu < \lambda:\text{cf}(\mu) < \mu, {\text{\rm pp\/}}(\mu) \ge
\lambda\}$. 
\ermn
\medskip

{\stag{4.18}\ub{ Dream}}:\quad  Prove (WH), i.e., (WH)$_1$ or (WH)$_2$.
\medskip

This is really like a dream: I do not believe in it; but it is the best
substitute for GCH which has not been proved impossible (essentially).  The
dual dream is not to prove its failure but essentially to prove that there
are no more rules or at least to show that some pcf ``bizarre" structures
are possible:
\medskip

\ub{\stag{4.19} The forcer Dream}:\quad Prove the consistency of: {\it there
is a set ${\frak a}$ of regular cardinals $> |{\frak a}|$ such that for some
inaccessible (= regular limit uncountable) cardinal $\lambda$ we have
$\lambda = \sup(\lambda \cap \text{ pcf}({\frak a}))$}. This will be enough
for proving the consistent failure of pcf(pcf$({\frak a})) = \text{
pcf}({\frak a})$. 
\medskip

A quite reasonable hope is:
\medskip

\ub{\stag{4.20} Question}:\quad For every $\lambda\ge\aleph_\omega$ there is
$n<\omega$ for which we have: 
$$\text{for no }\mu < \lambda\text{ do we have cf}(\mu)>\aleph_n \and
{\text{\rm pp\/}}_{\text{cf}(\mu)\text{-complete}}(\mu) > \lambda.$$
\medskip

A related statement for $\beth_\omega$ instead of $\aleph_\omega$ was proved
in \cite{Sh:460} and put forward as a positive solution of Hilbert first
problem, i.e., GCH; we have to see whether this is justified and/or
accepted. 

Note that \scite{4.20}  will be enough for deriving consequences of
$2^{\aleph_0} > \aleph_\omega$, confirming the following
\medskip

\ub{\stag{4.21} Thesis}:\quad pcf theory will make failures of GCH
semi-axioms, i.e., from some cardinal arithmetic equation or statement 
we shall prove many consequences. Note that whereas earlier $2^{\aleph_0} =
\aleph_1$, and (in few cases) $2^{\aleph_0} < \aleph_\omega$ has been used
as an assumption now there are some cases in which, e.g., $2^{\aleph_0} >
\aleph_\omega$ is used. 
\medskip

In particular we may consider a singular strong limit cardinal $\mu$, any
$\mu=\sum\{\mu_n:n\in \Bbb N\}$ with $2^{\mu_n}<\mu_{n+1}$ will do. We
already know that $2^\mu=\mu^+$ has serious consequences. Now we can hope
that also $2^\mu>\mu^+$ has. It implies that there is a $\mu^+$-free
non-free Abelian group of cardinality $\mu^+$ (\cite[5.5]{Sh:E12}).  Also
for some increasing sequence $\langle\lambda_n:n \in \Bbb N \rangle$ of
regular cardinals with limit $\mu$, $\prod\limits_n\lambda_n$ (ordered by
$<_J$ for the ideal $J$ of bounded subsets of $\Bbb N$) has the true
cofinality $\mu^{++}$ (this is an instance of a major theme of \cite{Sh:g}
that pp$(\mu)$ is the ``true" $\mu^{\text{cf}(\mu)}$, see \cite[IX]{Sh:g}). 
So it is not unreasonable to hope that this will provide a good dividing
line. 

The proofs in pcf depend very little on the advances in set theory from the
sixties on (in \cite{Sh:g} it has been claimed that Cantor, arising from his
grave would be able to understand them; certainly he could understand the
theorems). It may well be that the next stage in the evolution will be
\medskip

\ub{\stag{4.22} Dream}:\quad Combine the methods of pcf and inner models  
to answer questions as above, see \cite{Sh:413}.
\bigskip \bigskip 

\head {\S6 Concluding remarks} \endhead  \resetall \sectno=6
\medskip

We discuss here, usually concisely and with scant references and background 
various things.  Not that they are less worthwhile than those discussed
earlier, it's just that I have not found the time to write on them in
leisure.

What is model theory?  Classically you have a theory $T$ and the class
EC$(T)$ of all its models, i.e., structures of a fixed signature (e.g., for
rings we have $+,\times,0,1$) that satisfy $T$; the main case was $T$ first
order, then the class is called first order or elementary, but other logics
have been important as well; the point has been the interplay between what
can be seen in $T$ (syntactical side) and what can be said on EC($T$),
semantical side. For example, a sentence $\psi$ (in first order logic) is
preserved by submodels (a semantical property) iff $\psi$ is equivalent to a
universal sentence (syntactical). See the introduction to the Berkeley
Symposium in 1964.

There were some other frameworks. Some have suggested to look at less
general structures.  We may look at computable (= recursive) structures, we
may look at Polish structures (i.e., for a separable complete metric space
with the universal of the model and function are continuous and relations,
e.g., closed) we may look at Borel structures, at $\Sigma^1_n$-structures or
structure from $\bold L[\Bbb R]$.

Could you choose any of them as the main framework?  If you look at
computable structures, are the atomic relations computable or the first
order definable ones are computable?  Why Borel and not closed (like Polish
structures) or $F_\sigma$?  Also dealing with Borel forces you to deal with
$\Sigma^1_1$ and hence higher. $\bold L[\Bbb R]$ is a natural stopping
point but it is highly set-theoretic sensitive -- we have very different
picture if $\bold V = \bold L$ on the one hand and if we assume AD$_{\bold
L[\Bbb R]}$ on the other hand (which requires high consistency
strength). Hence
\medskip

\ub{\stag{nf.0} Model Theory Content Thesis}:\quad  1) While all those
``restrictive frameworks" are interesting, their (absence of) closure
properties make them unnatural as the focus of model theory. \nl
2) Also none of them is ``\ub{the} effective framework". \nl
3) To advance with most of them we should better start with a general theory  
of non-elementary classes (see below). \nl
4) Also if you are interested in universal first order classes or the
existentially closed members of such classes, you better do it based on a
theory for general first order classes.  
\medskip

\ub{\stag{6.1A} ZFC vs Model Theory Thesis}:\quad  1) It is much preferable 
to have the theorem in ZFC or at least ZFC + a semi-axiom, like CH. \nl
2) If (1) fails, using a simple division of the $\lambda$'s or proving
one side in ZFC, the other by consistency are reasonable substitutes.
\medskip

We may strengthen first order logic.  We may allow infinite conjunction,
i.e. $\bigwedge\limits_n \varphi_n$, while still every formula have finitely
many free variables, this is $\Bbb L_{\aleph_1,\aleph_0}$ (similarly $\Bbb
L_{\lambda^+,\aleph_0}$ if we allow $\bigwedge\limits_{\alpha<\lambda}
\varphi_\alpha$). Now $\Bbb L_{\aleph_1,\aleph_0}$ has many good properties;
every sentence with an infinite model has a countable one (the D.L.S.,
downward L\"owenheim--Skolem theorem), a completeness theorem (with an
infinitary rule), interpolation theorem (hence implicit definability implies
explicit definability).  $\Bbb L_{\lambda^+,\aleph_0}$ still has the
D.L.S. but usually not the others. 

For me a central topic of model theory is classification\footnote{Note: it
does not necessarily generalize stability theory.  Probably better to say
taxonomy theory, as some people mistakenly interpret my intention in
``classification theory" as finding the first order classes for which every
model can be characterized up to isomorphism by (general) cardinal
invariants (an important case, the major gap), harmless to let it have also
this meaning.}  theory for elementary classes, i.e., those defined as the
class of models of a first order logic theory.  The idea being finding
natural dividing lines: positive theory for the lower halves (understanding
the structure) and on the upper part (proving they have complicated models).
Much has been done in the investigation of the theory of low classes
(particularly stable but also simple and little on others). There is much to
be done on them. 
\medskip

\ub{\stag{nf.1} High Taxonomy Dream}:\quad  Find a good dividing line (or
lines), which is much higher than those mentioned above (see more in
\cite{Sh:702}).
\medskip

\ub{\stag{nf.2} The Mountain Air Thesis}:\quad The air on high mountains is 
clearer, many aspects are more transparent when we work in a more general
thesis.  In particular: 
\mr
\item "{$(a)$}"   Even if you are interested just in the model theory
of specific structures (like the real and the $p$-adic field), general model
theory will help you.
\sn
\item "{$(b)$}"   Even if you are interested just in countable models of
first order classes EC$(T)$ you will be helped by looking at
$\kappa$-saturated models of $T$ of cardinality $> \kappa$.
\sn
\item "{$(c)$}"   It is good to prove an ``outside" 
property Pr$_T(\lambda)$ which speaks on $\{M:M \models T$ has cardinality
$\lambda\}$ of a first order $T$ does not depend on $\lambda$ by proving it
equivalent to an ``inside" syntactical property.  It is also an excellent
way to discover interesting syntactical properties of $T$ even if you have
no interest in what occurs in every $\lambda$.
\sn 
\item "{$(d)$}"  Model theory of non-first-order classes will help the
elementary case and may be more relevant than first order for investigating
some natural uncountable structures. 
\endroster
\medskip

\ub{\stag{nf.3} The Specific Stability Dream}:\quad  Find interesting
natural first order classes EC$(T)$ with an interesting model theory and in
particular stability theory; I mean that the general methods of model theory
are really combined with the investigation of those specific classes (not
just quoting the results).
\medskip

This has occurred for differentially closed fields (and to some
extent, separably closed not algebraically closed fields).
\medskip

\ub{\stag{nf.4} Question}:\quad Are there infinite non-separably closed
fields $F$ with stable Th$(F)$?
\medskip

\ub{\stag{nf.5} The Specific Taxonomy Dream}:\quad Carry out the
classification program for the family of interesting specific classes
EC$(T)$, i.e., on $\{EC(T'):T'\supseteq T$ is complete$\}$.  This is not the
same as \scite{nf.3}, there we are interested in $T$ in which we can say
much on $M \in EC(T)$, here we look for dividing lines in the family of
completions of $T$; e.g., $T =$ the theory of fields.
\medskip

\ub{\stag{ne.1} Problem}:\quad  What are the right contexts for stability
theory? Develop the in the various specific cases.
\medskip

Various contexts we are considering:
\mr
\item "{$(A)$}"  first order classes (with $\prec$),
\sn
\item "{$(B)$}"  EC$(T)$, for $T$ a universal first order theory with
amalgamation under $\subseteq$,
\sn
\item "{$(C)$}"  the class of existentially closed models of a first order
theory $T$ under $\subseteq$ (here naturally we do not have negation, i.e.,
the negation of a formula is an infinite disjunction of one), 
\sn
\item "{$(D)$}"  $D$-homogeneous models with $\prec$ (so we restrict EC$(T)$
to $\{M:\{\text{tp}(\bar a,\emptyset,M):\bar a \in {}^{\omega>} M\}
\subseteq D$ and assume a strong form of amalgamation),
\sn
\item "{$(E)$}"  universal classes (\cite{Sh:300}),
\sn
\item "{$(F)$}"  EC$(T),T \subseteq \Bbb L_{\lambda^+,\omega}$ and
$\prec_\Delta,\Delta$ a fragment of $\Bbb L_{\lambda^+,\omega}$,
\sn
\item "{$(G)$}"  abstract elementary class with amalgamation,
\sn
\item "{$(H)$}"  abstract elementary class with no maximal model,
\sn
\item "{$(I)$}"  abstract elementary class,
\sn
\item "{$(J)$}"  good frames and relatives (\cite{Sh:705}).
\ermn
Note that one of the desired properties is having good closure properties
(e.g., for interpreting groups); so a wider context may be of interest as it
has better closure properties.
\medskip

A hard test is the following
\medskip

\ub{\stag{ne.2} The main gap Question}:\quad  Prove a form of the main gap
for $\psi \in \Bbb L_{\lambda^+,\omega}$ (or just $\Bbb
L_{\aleph_1,\aleph_0}$); i.e., for every such $\psi$ either $I(\lambda,\psi)
> \lambda$, see below, for every $\lambda$ large enough or there is an
ordinal $\gamma$ such that for every ordinal $\alpha$, $I(\aleph_\alpha,\psi)
\le \beth_\gamma(|\alpha|)$.
\medskip 

\definition{\stag{ne.3} Definition}\quad  $I(\lambda,\psi) = \{M/\cong:M
\models \psi\}$, similarly $I(\lambda,K)$ for $K$ a class of models. 
\enddefinition
\medskip

As in the case of the first order, the intention is 
\medskip

\ub{\stag{ne.3a} Main gap a good test Thesis}:\quad  Solving \scite{ne.2}
will force you to develop a theory, find interesting definitions and
theorems on them. 
\medskip

In the seventies soft model theory (i.e., with the logic as a variable) had
been very popular (see \cite{BF}), but has gone out of fashion, probably as
it seemed that there were many counterexamples and few theorems.  I do not
see it as a final verdict. 
\medskip

\ub{\stag{ne.4} Dream}:\quad  Find natural properties of logics and
nontrivial implications between them (giving a substantial mathematical
theory, of course). 
\medskip

\ub{\stag{ne.5} Dream}:\quad  Find a new logic with good model theory (like 
compactness, completeness theorem, interpolation and those from
\scite{ne.4}) and strong expressive power preferably concerning other parts
of mathematics (see \cite{Sh:702}, possibly specifically derive for them). 
\medskip

\ub{\stag{ne.6} Problem}:\quad  Develop the model theory of first order
classes EC$(T)$ with $T$ having the finite model property (i.e., every
finite subset of $T$ has a finite model). 
\medskip

\ub{\stag{ne.7} Problem}:\quad  Develop and investigate a logic for
``polynomial invariants for graphs and general structures" (there are many
worthwhile definitions of polynomials (invariants)  for a graph, which
depend on the isomorphism type only). 
\medskip

\centerline {$ * \qquad * \qquad *$}
\medskip

In fact, I am not fond of set theory without choice, nevertheless:
\medskip

\ub{\stag{st.1} Problem}:\quad  Develop combinatorial set theory for
universes with limited amount of choice (see in \cite{Sh:666}). 
\medskip

\ub{\stag{st.2} Problem}:\quad  Develop descriptive set theory for
${}^\kappa \mu$ in particular ${}^\omega \mu,\mu$ strong limit singular of
cofinality $\aleph_0$ (see discussion in \cite{Sh:724}; probably we first of
all need good questions as straight generalization of properties may not
succeed). 
\medskip

\ub{\stag{ns.1} Dream}:\quad  Try to formalize and really say
something\footnote{here, as elsewhere, e.g., \scite{ns.2}, \scite{3.4}  
the dream is not to suggest a reasonable framework for this, but build a
mathematical theory} on mathematical beauty and depth.  Of course (length of
proof)/(length of theorem) is in the right direction, etc.
\medskip

\ub{\stag{ns.2} Dream}:\quad Make a reasonable mathematical theory when we
restrict ourselves to the natural numbers up to $n$, where $n$ is a specific
natural number (say $2^{2^{100}} + 1$) (e.g., thinking our universe is
discrete with this size).

     \shlhetal 

\newpage
    
REFERENCES.  
\bibliographystyle{lit-plain}
\bibliography{lista,listb,listx,listf,liste}

\enddocument